\documentclass{amsart}
\usepackage{amsfonts}
\usepackage{pstricks,pstcol,pst-plot,pst-node,pst-tree,amssymb}
\usepackage{pstricks-add}
\usepackage{graphicx}

\numberwithin{equation}{section}
\input{tcilatex}

\begin{document}
\title[Singular set of a $C^{1}$ smooth surface ]{A Codazzi-like equation
and the singular set for $C^{1}$ smooth surfaces in the Heisenberg group}
\author{Jih-Hsin Cheng}
\address[Jih-Hsin Cheng and Jenn-Fang Hwang]{Institute of Mathematics,
Academia Sinica, Taipei, 11529 and National Center for Theoretical Sciences,
Taipei Office, Taiwan, R.O.C.}
\email[Cheng]{cheng@math.sinica.edu.tw}
\author{Jenn-Fang Hwang }
\email[Hwang]{majfh@math.sinica.edu.tw}
\author{Andrea Malchiodi}
\address[Malchiodi]{SISSA, Via Beirut 2-4, 34014, Trieste, Italy}
\email{malchiod@sissa.it}
\author{Paul Yang}
\address[Yang]{Department of Mathematics, Princeton University, Princeton,
NJ 08544, U.S.A.}
\email{yang@Math.Princeton.EDU}
\thanks{}

\begin{abstract}
In this paper, we study the structure of the singular set for a $C^{1}$
smooth surface in the $3$-dimensional Heisenberg group $\boldsymbol{H}_{1}$.
We discover a Codazzi-like equation for the $p$-area element along the
characteristic curves on the surface. Information obtained from this
ordinary differential equation helps us to analyze the local configuration
of the singular set and the characteristic curves. In particular, we can
estimate the size and obtain the regularity of the singular set. We
understand the global structure of the singular set through a Hopf-type
index theorem. We also justify that Codazzi-like equation by proving a
fundamental theorem for local surfaces in $\boldsymbol{H}_{1}$.
\end{abstract}

\maketitle

\begin{center}
{\Large Table of Contents}
\end{center}

\ \ \ \ \ {\large 1. Introduction and statement of the results}

\ \ \ \ \ {\large 2. A Codazzi-like equation and properties:}

\ \ \ \ \ \ \ \ \ {\large proofs of Theorem A and Theorem B}

{\large \ \ \ \ \ 3. Local configuration of the singular set}

{\large \ \ \ \ \ 4. Examples}

{\large \ \ \ \ \ 5. Size and regularity of the singular set}

{\large \ \ \ \ \ 6. The local theory of surfaces with}

\ \ \ \ \ \ \ \ \ {\large prescribed }$p${\large -mean curvature}

{\large \ \ \ \ \ 7. Index of the singular set}

{\large \ \ \ \ \ 8. Generalization to pseudohermitian manifolds}

{\large \ \ \ \ \ 9. Appendix: a generalized ODE}

\bigskip

\section{Introduction and statement of the results}

In the recent years, the $p$-minimal (or $\boldsymbol{H}$-minimal) surfaces
have been studied extensively in the framework of geometric measure theory
(e.g., \cite{GN96}, \cite{FSS01}, \cite{Pau01}) and from the viewpoint of
partial differential equations and that of differential geometry (e.g., \cite%
{CHMY04}, \cite{CH04}, \cite{CHY1}, \cite{CHY2}). Motivated by the
isoperimetric problem in the Heisenberg group, one also studied nonzero
constant $p$-mean curvature surfaces and the regularity problem (e.g., \cite%
{Pan82}, \cite{CDG}, \cite{LM}, \cite{LR}, \cite{RR}, \cite{MR}, \cite{Pauls}%
). In fact, the notion of $p$-mean curvature ("$p$-" stands for
"pseudohermitian") can be defined for (hyper)surfaces in a pseudohermitian
manifold. The Heisenberg group as a (flat) pseudohermitian manifold is the
simplest model example, and represents a blow-up limit of general
pseudohermitian manifolds.

The equation of prescribed $p$-mean curvature in the $3$-dimensional
Heisenberg group $\boldsymbol{H}_{1}$ is one of few known equations having
geometric significance in $2D.$ For the Plateau or Dirichlet problem with
smooth boundary value, we have reasons to believe that the minimizer is at
least $C^{1}$ (but not $C^{2}$ in general$).$ In \cite{CHY2} three of the
authors studied the regularity of the nonsingular portion of a $C^{1}$
smooth surface in $\boldsymbol{H}_{1}$. In this paper we study the local
structure of the singular set of such a surface through an ordinary
differential equation along the characteristic curves. Results on the local
structure of the singular set will be used in studying the global structure
of the singular set later. Note that the local structure of the singular set
for $C^{2}$ smooth surfaces has been classified. Namely, on a $C^{2}$ smooth
surface, a singular point is either isolated or passed through by a $C^{1}$
smooth singular curve in a neighborhood under a mild condition on the $p$%
-mean curvature (see Theorem 3.3 in \cite{CHMY04} and remarks for the
general case in Section 7 there). On the other hand, the structure of the
singular set of a $C^{1}$ smooth surface can be complicated as we will see
in this paper. The understanding of the singular set of a $C^{1}$ smooth
surface would help to solve the isoperimetric problem for $C^{1}$ smooth
domains in $\boldsymbol{H}_{1}$ (see \cite{RR} for more details on the
isoperimetric problem).

Let $\Omega $ be a domain of $R^{2}$ (by a domain we mean an open and
connected set) and let $u\in C^{1}(\Omega )$ be a weak solution to 
\begin{equation}
\func{div}\frac{\nabla u+\vec{F}}{|\nabla u+\vec{F}|}=H  \label{1.1}
\end{equation}%
\noindent where $\vec{F}$ ($H,$ resp.) is an $L_{loc}^{1}$ vector field
(function, resp.) in $\Omega ,$ that is, for any $\varphi \in C_{0}^{\infty
}(\Omega ),$ there holds%
\begin{equation}
\int_{S_{\vec{F}}(u)}|\nabla \varphi |+\int_{\Omega \backslash S_{\vec{F}%
}(u)}\frac{\nabla u+\vec{F}}{|\nabla u+\vec{F}|}\cdot \nabla \varphi
+\int_{\Omega }H\varphi \geq 0  \label{1.1'}
\end{equation}

\noindent in which $S_{\vec{F}}(u)$ denotes the singular set of $u$,
consisting of the points (called singular points) where $\nabla u+\vec{F}$ $%
= $ $0$ (see (1.2) and (1.3) in \cite{CHY2}). A point $p$ $\in $ $\Omega $
is called nonsingular if $\nabla u+\vec{F}$ $\neq $ $0$ at $p.$ We call $%
\Omega $ nonsingular if all the points of $\Omega $ are nonsingular.

Note that a $C^{2}$ smooth solution (i.e. satisfying (\ref{1.1}) at
nonsingular points) may or may not be a $C^{1}$ weak solution. The reason is
that (\ref{1.1'}) implies some equal-angle condition along $S_{\vec{F}}(u)$
if $S_{\vec{F}}(u)$ is a $C^{1}$ smooth curve, which a $C^{2}$ smooth
solution may not satisfy (see Example 7.3 in \cite{CHY1}). From the
variational point of view, that $u$ satisfies the condition (\ref{1.1'}) is
more natural than that $u$ satisfies (\ref{1.1}) pointwise at nonsingular
points. Therefore we study the solutions satisfying (\ref{1.1'}).

Let $N$ denote the planar vector $\frac{\nabla u+\vec{F}}{|\nabla u+\vec{F}|}
$ at a nonsingular point$.$ Since $N$ is of unit length, we can write 
\begin{equation*}
N=(\cos \theta ,\sin \theta )
\end{equation*}%
\noindent locally for some angular function $\theta $ ($\in $ $C^{0}$ if we
assume $\vec{F}$ $\in $ $C^{0}).$ Let $D$ $:=$ $|\nabla u+\vec{F}|$ ($\in $ $%
C^{0}$ a priori) and let $N^{\perp }$ $:=$ $(\sin \theta ,$ $-\cos \theta ).$
Suppose further that $\vec{F}$ $\in $ $C^{1}$ and $N(H)$ exists and is
continuous. In \cite{CHY2} (see Theorem D therein), we proved that $\theta $
is in fact $C^{1}$ and $N^{\perp }D$ exists and is continuous. In this paper
we will show that $N^{\perp }(N^{\perp }D)$ exists and is continuous.
Moreover, $D$ satisfies an ordinary differential equation of second order
along any characteristic curve (i.e., integral curve of $N^{\perp })$ (see
Theorem A below). For $\vec{F}$ $=$ $(F_{1},$ $F_{2})$ $\in $ $C^{1}(\Omega
) $ we define%
\begin{equation*}
curl\vec{F}:=(F_{2})_{x}-(F_{1})_{y}
\end{equation*}

\noindent (note that in \cite{CHY2} we used $rot\vec{F}$ instead of $curl%
\vec{F}).$ Denote $N^{\perp }D$ and $N^{\perp }(N^{\perp }D)$ by $D^{\prime
} $ and $D^{\prime \prime },$ respectively.

\bigskip

\textbf{Theorem A}$\mathbf{.}$ \textit{Let }$u$\textit{\ }$\in $ $%
C^{1}(\Omega )$\textit{\ be a weak solution to (}$\ref{1.1})$\textit{\ with }%
$\vec{F}$ $\in $ $C^{1}(\Omega )$ \textit{and }$H$\textit{\ }$\in $\textit{\ 
}$C^{0}(\Omega )$ \textit{such that }$\Omega $\textit{\ is nonsingular}$.$%
\textit{\ Suppose further }$N^{\perp }(curl\vec{F})$ and $N(H)$\textit{\
exist and are continuous. Then }$D^{\prime }$\textit{\ and }$D^{\prime
\prime }$\textit{\ exist and are continuous in }$\Omega .$\textit{\
Moreover, }$D$\textit{\ satisfies the following differential equation}%
\begin{equation}
DD^{\prime \prime }=2(D^{\prime }-\frac{curl\vec{F}}{2})(D^{\prime }-curl%
\vec{F})+(N^{\perp }(curl\vec{F}))D+(H^{2}+N(H))D^{2}.  \label{1.2}
\end{equation}

\bigskip

The proof of (\ref{1.2}) is based on the construction of so called $s,$ $t$
coordinates using $N^{\perp }$ and $N$ (see (\ref{2.4}))$.$ The associated
integrability condition

\begin{equation*}
(\theta _{s})_{t}=(\theta _{t})_{s}
\end{equation*}

\noindent (see (\ref{2.5}), (\ref{2.6})) is exactly (\ref{1.2}). So (\ref%
{1.2}) is a Codazzi-like equation.

When $H$ is the $p($or $H)$-mean curvature, we have $\vec{F}$ $=$ $(-y,x)$
in this case$,$ so $curl\vec{F}$ $=$ $2.$ The $p($or $H)$-mean curvature is
a notion to measure a hypersurface, which respects the ambient
pseudohermitian structure (see \cite{CHMY04}). Equation (\ref{1.2}) can then
be reduced to%
\begin{equation}
DD^{\prime \prime }=2(D^{\prime }-1)(D^{\prime }-2)+(H^{2}+N(H))D^{2}.
\label{1.3}
\end{equation}

\noindent For a $p($or $H)$-minimal graph, we have $H$ $\equiv $ $0,$ so (%
\ref{1.3}) can be further reduced to%
\begin{equation}
DD^{\prime \prime }=2(D^{\prime }-1)(D^{\prime }-2).  \label{1.4}
\end{equation}

Equation ($\ref{1.4})$ is integrable. Namely, we observe that%
\begin{eqnarray}
\frac{2D^{\prime }}{D} &=&\frac{D^{\prime }D^{\prime \prime }}{(D^{\prime
}-1)(D^{\prime }-2)}  \label{1.5} \\
&=&\frac{-D^{\prime \prime }}{D^{\prime }-1}+\frac{2D^{\prime \prime }}{%
D^{\prime }-2}  \notag \\
&=&-\frac{(D^{\prime }-1)^{\prime }}{D^{\prime }-1}+\frac{2(D^{\prime
}-2)^{\prime }}{D^{\prime }-2}  \notag
\end{eqnarray}

\noindent at the points of a characteristic curve (line), where $D^{\prime }$
$\neq $ $1,2.$ Integrating (\ref{1.5}), we obtain%
\begin{equation}
|D^{\prime }-2|^{2}=c|D^{\prime }-1|D^{2}  \label{1.5'}
\end{equation}

\noindent for some constant $0<c<\infty .$ It is not hard to see that if $%
D^{\prime }$ $\neq $ $1$ ($\neq 2,$ resp.) at some nonsingular point $q$,
then $D^{\prime }$ $\neq $ $1$ ($\neq 2,$ resp.) on the whole characteristic
curve (line) $\Gamma $ passing through $q.$ If $D^{\prime }$ $=$ $2$ ($=$ $%
1, $ resp.) at $q,$ then $c$ $=$ $0$ ($=$ $\infty ,$ resp.$)$ and $D^{\prime
}$ $\equiv $ $2$ ($\equiv $ $1,$ resp.) on $\Gamma $ by the uniqueness of
solutions to (\ref{1.4}), an ordinary differential equation$.$ When a
nonsingular point tends to a singular point along a characteristic curve
(line), either $D^{\prime }$ goes to $2$ or $D^{\prime }$ goes to $1$ ($%
D^{\prime }\equiv 1$ in fact in this case). In general, $(N^{\perp }(curl%
\vec{F}))D$ $+$ $(H^{2}+N(H))D^{2}$ $\neq $ $0$ and hence (\ref{1.2}) is not
integrable. But near singular points (where $D$ $=$ $0),$ we consider $%
(N^{\perp }(curl\vec{F}))D$ $+$ $(H^{2}+N(H))D^{2}$ to be a small
perturbation term and obtain the following result.

\bigskip

\textbf{Theorem B}$\mathbf{.}$ \textit{Let }$u$\textit{\ }$\in $ $%
C^{1}(\Omega )$\textit{\ be a weak solution to (}$\ref{1.1})$\textit{\ with }%
$\vec{F}$ $\in $ $C^{1}(\Omega )$ \textit{and }$H$\textit{\ }$\in $\textit{\ 
}$C^{0}(\Omega ).$\textit{\ Assume further }$N^{\perp }(curl\vec{F})$ 
\textit{and} $N(H)$\textit{\ exist and are continuous (extended over
singular points) in }$\Omega .$ \textit{Let }$p$\textit{\ }$\in $\textit{\ }$%
\Omega $\textit{\ be a singular point. Let }$\Gamma $\textit{\ : }$[0,\bar{%
\rho})\rightarrow $\textit{\ }$\Omega $\textit{\ }$\in $\textit{\ }$C^{1}$%
\textit{\ be such that }$\Gamma (0)=p$\textit{\ and }$\Gamma ((0,\bar{\rho}%
)) $\textit{\ is a characteristic curve with unit-speed parameter }$\rho $ $%
\in $ $(0,\bar{\rho}).$\textit{\ Suppose }$curl\vec{F}(p)$\textit{\ }$\neq $%
\textit{\ }$0.$\textit{\ Then the following statements hold:}

\textit{(a) We have either}%
\begin{equation}
\lim_{\rho \rightarrow 0}D^{\prime }(\Gamma (\rho ))=\frac{curl\vec{F}(p)}{2}%
\text{ or }\lim_{\rho \rightarrow 0}D^{\prime }(\Gamma (\rho ))=curl\vec{F}%
(p).  \label{1.6}
\end{equation}

\textit{(b)} \textit{The sign of }$curl\vec{F}(p)$\textit{\ determines the
direction of }$N^{\perp }.$\textit{\ That is, if }$curl\vec{F}(p)$\textit{\ }%
$>$\textit{\ }$0,$\textit{\ then }$N^{\perp }$\textit{\ }$=$\textit{\ }$%
\frac{\partial }{\partial \rho }$\textit{\ while if }$curl\vec{F}(p)$\textit{%
\ }$<$\textit{\ }$0,$\textit{\ then }$N^{\perp }$\textit{\ }$=$\textit{\ }$-%
\frac{\partial }{\partial \rho }.$

\textit{(c)} \textit{Let }$p,$\textit{\ }$q$\textit{\ be two distinct
singular points in }$\Omega .$\textit{\ Suppose }$curl\vec{F}$\textit{\ }$%
\neq $\textit{\ }$0$\textit{\ in} $\Omega .$ \textit{Then there does not
exist }$\Gamma $\textit{\ }$:$\textit{\ }$[0,\bar{\rho}]\rightarrow $\textit{%
\ }$\Omega $\textit{\ }$\in $\textit{\ }$C^{1}$\textit{, a characteristic
curve on }$(0,\bar{\rho})$\textit{\ with }$\Gamma (0)$\textit{\ }$=$\textit{%
\ }$p$\textit{\ and }$\Gamma (\bar{\rho})$\textit{\ }$=$\textit{\ }$q.$

\bigskip

In the Appendix we consider a generalized version of equation (\ref{1.2})
and prove a result analogous to Theorem B (a) (see Theorem A.1). Applying
Theorem A.1 to a general situation for a surface in a pseudohermitian
3-manifold (see (\ref{O.23})), we obtain a variant of Theorem B (a) (see
Theorem B$^{\prime }$ in Section 8).

Theorem B (b) will be applied to show impossibility of some situations. For
instance, in the case of a $p($or $H)$-minimal graph, if a family of
characteristic lines converges to another line, then the limit line cannot
contain any singular point and must be a characteristic line (see Lemma 3.2).

Theorem B (c) will be used often in the study of the configuration of
singular points and characteristic curves in Section 3. Among other things,
we have a result about the local structure of singular points. Recall that
for a $C^{2}$ smooth surface, a singular point is either isolated or passed
through by a $C^{1}$ smooth singular curve in a neighborhood under a mild
condition on $H$ (see Theorem 3.3 in \cite{CHMY04}).

\bigskip

\textbf{Theorem C.} \textit{Let }$u$\textit{\ }$\in $ $C^{1}(\Omega )$%
\textit{\ be a weak solution to (}$\ref{1.1})$\textit{\ with }$\vec{F}$ $\in 
$ $C^{1}(\Omega )$ \textit{and }$H$\textit{\ }$\in $\textit{\ }$C^{0}(\Omega
).$\textit{\ Assume further }$N^{\perp }(curl\vec{F})$ \textit{and} $N(H)$%
\textit{\ exist and are continuous (extended over singular points) in }$%
\Omega .$ \textit{Let }$p$\textit{\ be a singular point in }$\Omega .$%
\textit{\ Suppose }$curl\vec{F}(p)$\textit{\ }$\neq $\textit{\ }$0.$ \textit{%
Then we have}

\textit{(a)} \textit{Either }$p$\textit{\ is an isolated singular point,
i.e.,} \textit{there exists a neighborhood }$V$\textit{\ }$\subset $\textit{%
\ }$\Omega $\textit{\ of }$p$\textit{\ such that }$V$\textit{\ contains no
other singular points except }$p,$\textit{\ or there exists at least one }$%
C^{0}$\textit{\ singular curve }$\gamma $\textit{\ }$:$\textit{\ }$[0,1]$%
\textit{\ }$\rightarrow $\textit{\ }$\Omega $\textit{\ (i.e., }$\gamma $ 
\textit{is continuous and} $\gamma (s)$\textit{\ is a singular point for
each }$s$\textit{\ }$\in $\textit{\ }$[0,1])$\textit{\ such that }$\gamma
(0) $\textit{\ }$=$\textit{\ }$p$\textit{\ and }$\gamma (1)$\textit{\ }$\neq 
$\textit{\ }$p.$

\textit{(b) Moreover, there is a neighborhood }$U$ \textit{of }$p$ \textit{%
such that for any singular point }$q$ $\in $ $U$ \textit{there exists a }$%
C^{0}$\textit{\ singular curve }$\beta $\textit{\ }$:$\textit{\ }$[0,1]$%
\textit{\ }$\rightarrow $\textit{\ }$U$\textit{\ with }$\beta (0)=p$\textit{%
\ and }$\beta (1)=q.$

\textit{\bigskip }

Theorem C (b) implies that the singular set is locally path-connected. It is
possible to construct examples having several singular curves meeting at a
singular point (see Examples 4.2, 4.3). In case a singular point is
isolated, we can describe the local configuration of characteristic curves
near such a singular point in a general situation.

\bigskip

\textbf{Theorem D}. \textit{Let }$u$\textit{\ }$\in $ $C^{1}(\Omega )$%
\textit{\ be a weak solution to (}$\ref{1.1})$\textit{\ with }$\vec{F}$ $\in 
$ $C^{1}(\Omega )$ \textit{and }$H$\textit{\ }$\in $\textit{\ }$C^{0}(\Omega
).$ \textit{Assume further }$N^{\perp }(curl\vec{F})$\textit{\ and} $N(H)$%
\textit{\ exist and are continuous (extended over singular points) in }$%
\Omega .$ \textit{Let }$p$\textit{\ be an isolated singular point in }$%
\Omega .$ \textit{Suppose }$curl\vec{F}(p)$\textit{\ }$\neq $\textit{\ }$0.$ 
\textit{Then there exists }$r_{0}$\textit{\ }$>$\textit{\ }$0$\textit{\ such
that for any }$q$\textit{\ }$\in $\textit{\ }$\bar{B}_{r_{0}}(p)\backslash
\{p\}$\textit{, }$q$\textit{\ is nonsingular and the characteristic curve }$%
\Gamma _{q}$ \textit{passing through }$q$\textit{\ has to meet }$p.$ \textit{%
Moreover, the unit tangent vector }$N^{\perp }$\textit{\ of }$\Gamma _{q}$%
\textit{\ has a limit at }$p,$\textit{\ denoted }$v(q),$\textit{\ and the
map }$\psi $\textit{\ }$:$\textit{\ }$q$\textit{\ }$\in $\textit{\ }$%
\partial B_{r_{0}}(p)$\textit{\ }$\rightarrow $\textit{\ }$v(q)$\textit{\ }$%
\in $\textit{\ }$T_{p}\Omega $\textit{\ is a homeomorphism onto the space of
unit tangent vectors at }$p.$

\bigskip

According to Theorem A in \cite{CHY2}, $H$ (plus initial condition)
determines $\theta $ ($\frac{d\theta }{d\rho }$ $=$ $-H.$ Here we use $\rho $
instead of $\sigma $ as unit-speed parameter) and hence the characteristic
curves through (\ref{1.6}) of \cite{CHY2}. In fact, the characteristic
curves $(x(\rho ),$ $y(\rho ))$ satisfy the following system of ordinary
differential equations of second order%
\begin{equation}
\frac{d^{2}x}{d\rho ^{2}}=H\frac{dy}{d\rho },\text{ \ \ }\frac{d^{2}y}{d\rho
^{2}}=-H\frac{dx}{d\rho }.  \label{1.7}
\end{equation}%
\noindent From Theorem D the value of $u$ near $p$ is completely determined
by the values of $u$ at $p,$ $\vec{F},$ and $H$ by integrating $du$ $+$ $%
F_{1}dx$ $+$ $F_{2}dy$ $=$ $0$ along the characteristic curves$.$

\bigskip

\textbf{Corollary E}. \textit{Suppose we are in the situation of Theorem D.
Then the value of }$u$\textit{\ near }$p$\textit{\ is completely determined
by the values of }$u$\textit{\ at }$p,$ $\vec{F},$ and $H.$

\bigskip

We remark that in the case of $H$ $=$ $0$ and $\vec{F}$ $=$ $(-y,x),$ the
graph defined by $u$ in Corollary E is a plane. Let $\mathcal{H}^{2}$ denote
the $2$-dimensional Hausdorff measure. For the size of the singular set, we
have the following result.

\bigskip

\textbf{Theorem F. }\textit{Let }$u$\textit{\ }$\in $ $C^{1}(\Omega )$%
\textit{\ be a weak solution to (}$\ref{1.1})$\textit{\ with }$\vec{F}$ $\in 
$ $C^{1}(\Omega )$ \textit{and }$H$\textit{\ }$\in $\textit{\ }$C^{1}(\Omega
).$\textit{\ Assume further }$N^{\perp }(curl\vec{F})$ \textit{and} $N(H)$%
\textit{\ exist and are continuous (extended over singular points) in }$%
\Omega .$ \textit{Suppose also }$curl\vec{F}$\textit{\ }$\neq $\textit{\ }$0$%
\textit{\ in} $\Omega .$ \textit{Then }$\mathcal{H}^{2}(S_{\vec{F}}(u))$%
\textit{\ }$=$\textit{\ }$0.$

\bigskip

We remark that a $C^{1}$ smooth $p$-minimal graph defined by $u$ is a
special case of Theorem F$.$ Because of Theorem F the condition (\ref{1.1'})
for a $C^{1}$ smooth $p$-minimal graph defined by $u$ is reduced to%
\begin{equation}
\int_{\Omega }N\cdot \nabla \varphi =0  \label{1.8}
\end{equation}

\noindent (note that $H$ $=$ $0)$ for any $\varphi $ $\in $ $C_{0}^{\infty
}(\Omega ),$ where $N$ $:=$ $\frac{\nabla u+\vec{F}}{|\nabla u+\vec{F}|}$
with $\vec{F}$ $=$ $(-y,x).$ To study the regularity of a singular curve
passing through $p_{0}$, we define the notion of the (inverse) expanding
rate of characteristic curves along such a singular curve, $\lambda _{\pm
}(p_{0})$ (see (\ref{5.12.10})). $p_{0}$ is nondegenerate if both $\lambda
_{+}(p_{0})$ $\neq $ $0$ and $\lambda _{-}(p_{0})$ $\neq $ $0$ (see
Definition 5.1 and (\ref{5.12.10}))$.$ It is clear from the definition that
the set of nondegenerate singular points is relatively open in $S_{\vec{F}%
}(u).$ The singular curve passing through a nondegenerate point is $C^{0}$ a
priori for $u$ $\in $ $C^{1}.$ However, we have the regularity result and
the equal-angle condition as follows.

\bigskip

\textbf{Theorem G}. \textit{Let }$u$\textit{\ }$\in $ $C^{1}(\Omega )$%
\textit{\ be a weak solution to (}$\ref{1.1})$\textit{\ with }$\vec{F}$ $\in 
$ $C^{1}(\Omega )$ \textit{and }$H$\textit{\ }$\in $\textit{\ }$C^{1}(\Omega
).$\textit{\ Assume further }$N^{\perp }(\func{curl}\vec{F})$ \textit{and }$%
N(H)$\textit{\ exist and are continuous (extended over singular points) in }$%
\Omega .$\textit{\ Suppose }$\func{curl}\vec{F}(p)$\textit{\ }$\neq $\textit{%
\ }$0$\textit{\ for any nondegenerate singular point }$p.$ \textit{Then we
have}

\textit{(a)} \textit{the set of nondegenerate singular points\ consists of }$%
C^{1}$\textit{\ smooth curves.}

\textit{(b) two characteristic curves issuing from a nondegenerate singular
point }$p_{0}$ \textit{have the same angle with the tangent line of the
singular curve through }$p_{0}.$

\bigskip

To show Theorem G, we study a more general situation. Consider $\theta $ as
an independent variable, satisfying the equation%
\begin{equation}
\func{div}(\cos \theta ,\sin \theta )\equiv (\cos \theta )_{x}+(\sin \theta
)_{y}=H.  \label{1.8.1}
\end{equation}%
\noindent Suppose $\mathcal{H}^{2}(K)$ $=$ $0$ for a subset $K$ in $\Omega .$
A solution $\theta $ $\in $ $C^{0}(\Omega \backslash K)$ is a weak solution
to (\ref{1.8.1}) (with $H$ $\in $ $L_{loc}^{1}(\Omega ),$ say) if there holds%
\begin{equation*}
\int_{\Omega }N\cdot \nabla \varphi +\int_{\Omega }H\varphi =0
\end{equation*}%
\noindent for any $\varphi $ $\in $ $C_{0}^{\infty }(\Omega ),$ where we
have written $N$ as $(\cos \theta ,\sin \theta ).$ Suppose $\gamma $ $%
\subset $ $K$ is a $C^{0}$ curve$.$ We can still define what a nondegenerate
point of $\gamma $ is (see Definition 5.1). We also define what a crack
point of $\gamma $ is (see Definition 5.2). Roughly speaking, a crack point
is a point at which $N$ has different limits along two characteristic curves.

\bigskip

\textbf{Theorem G}$^{\prime }.$ \textit{Suppose }$\mathcal{H}^{2}(K)$\textit{%
\ }$=$\textit{\ }$0$\textit{\ for a subset }$K$\textit{\ in }$\Omega .$ 
\textit{Let }$\theta $\textit{\ }$\in $\textit{\ }$C^{0}(\Omega \backslash
K) $\textit{\ be a weak solution to (\ref{1.8.1}) with }$H$ $\in $ $%
C^{1}(\Omega )$\textit{. Then we have}

\textit{(a)} \textit{the set of nondegenerate crack points\ consists of }$%
C^{1}$\textit{\ smooth curves.}

\textit{(b) two characteristic curves issuing from a nondegenerate crack
point }$p_{0}$ \textit{have the same angle with the tangent line of the
curve of nondegenerate crack points through }$p_{0}.$

\bigskip

Since a singular point is a crack point, we obtain Theorem G from Theorem G$%
^{\prime }$ (see Section 5 for more details)$.$ On the other hand, in the
situation that $N$ ($\equiv $ $(\cos \theta ,\sin \theta ))$ is defined by $%
u $ as in Theorem G, we show that a crack point is in fact a singular point
(see Theorem 5.4).

In Section 2 we give the proofs of Theorem A and Theorem B. In Section 3 we
show those of Theorem C and Theorem D. Some crucial examples are given in
Section 4. Theorems F and G are proved in Section 5.

In Riemannian geometry, we have Gauss and Codazzi equations in the
submanifold theory. The fundamental theorem for hypersurfaces in a Euclidean
space says that the equations of Gauss and Codazzi are exactly the
integrability conditions for finding an isometric imbedding with prescribed
metric and second fundamental form (see, for instance, page 47 in \cite{KN}%
). In pseudohermitian geometry, we have the analogous fundamental theorem
for surfaces in the $3$-dimensional Heisenberg group $\boldsymbol{H}_{1}$.
For simplicity we work in the $C^{\infty }$ category.

\bigskip

\textbf{Theorem H}. \textit{Given a nonzero }$C^{\infty }$\textit{\ smooth
vector field }$V$\textit{, a positive }$C^{\infty }$\textit{\ smooth
function }$D,$\textit{\ and a }$C^{\infty }$\ \textit{smooth function }$H$%
\textit{\ on an open neighborhood }$U$\textit{\ of a point }$p$\textit{\ in
the }$(\xi ,\eta )$\textit{\ plane. Suppose }$D$\textit{\ satisfies the
equations}%
\begin{equation}
DD^{\prime \prime }=2(D^{\prime }-1)(D^{\prime }-2)+(H^{2}+P(H))D^{2}
\label{1.9}
\end{equation}

\begin{equation}
L_{V}P=-(\frac{2-D^{\prime }}{D})P+HV  \label{1.9.1}
\end{equation}%
\textit{\noindent for a nonzero }$C^{\infty }$\ \textit{smooth vector field }%
$P,$ \textit{transversal to }$V$\textit{, such that }$(V,$\textit{\ }$P)$%
\textit{\ has the same orientation as} $(\partial _{\xi },$ $\partial _{\eta
}),$ \textit{where we denote }$V(D),$\textit{\ }$V(V(D))$\textit{, and the
Lie derivative in the direction }$V$ \textit{by }$D^{\prime }$\textit{, }$%
D^{\prime \prime },$\textit{\ and }$L_{V},$ \textit{resp.. Then in a perhaps
smaller neighborhood }$U^{\prime }$\textit{\ }$\subset $\textit{\ }$U$%
\textit{\ of }$p$

\textit{(1) there exist a }$C^{\infty }$\textit{\ smooth
orientation-preserving diffeomorphism: }$(\xi ,\eta )$\textit{\ }$%
\rightarrow $\textit{\ }$(x$\textit{\ }$=$\textit{\ }$x(\xi ,\eta ),y$%
\textit{\ }$=$\textit{\ }$y(\xi ,\eta ))$\textit{\ from }$U^{\prime }$%
\textit{\ onto its image in }$R^{2}$ \textit{and a }$C^{\infty }$\textit{\
smooth function }$\theta $\textit{\ }$=$\textit{\ }$\theta (\xi ,\eta )$%
\textit{\ on }$U^{\prime }$\textit{\ such that }%
\begin{eqnarray}
V(\theta ) &=&-H  \label{1.10} \\
V(x) &=&\sin \theta ,\text{ }V(y)=-\cos \theta .  \notag
\end{eqnarray}%
\textit{\noindent In addition, the vector field }$N$\textit{\ satisfying }$%
N(x)$\textit{\ }$=$\textit{\ }$\cos \theta ,$\textit{\ }$N(y)$\textit{\ }$=$%
\textit{\ }$\sin \theta $\textit{\ is equal to }$P$ and \textit{has the
property that }%
\begin{equation}
N(\theta )=\frac{2-V(D)}{D}.  \label{1.11}
\end{equation}

\textit{(2) Moreover, there exists a }$C^{\infty }$\textit{\ smooth function 
}$z$\textit{\ }$=$\textit{\ }$z(\xi ,\eta )$\textit{\ on }$U^{\prime }$%
\textit{\ to make a }$C^{\infty }$\textit{\ smooth embedding: }$(\xi ,\eta )$%
\textit{\ }$\in $\textit{\ }$U^{\prime }$\textit{\ }$\rightarrow $\textit{\ }%
$(x$\textit{\ }$=$\textit{\ }$x(\xi ,\eta ),y$\textit{\ }$=$\textit{\ }$%
y(\xi ,\eta ),$\textit{\ }$z$\textit{\ }$=$\textit{\ }$z(\xi ,\eta ))$%
\textit{\ }$\in $\textit{\ }$\boldsymbol{H}_{1}$\textit{\ such that the
image is a }$C^{\infty }$\textit{\ smooth graph }$z$\textit{\ }$=$\textit{\ }%
$u(x,y),$\textit{\ }$D$\textit{\ }$=$\textit{\ }$\sqrt{%
(u_{x}-y)^{2}+(u_{y}+x)^{2}}$\textit{, }$N=\frac{(u_{x}-y,u_{y}+x)}{D},$%
\textit{\ and}%
\begin{equation}
\func{div}N=H  \label{1.12.1}
\end{equation}

\begin{equation}
\func{div}DV=2  \label{1.12.2}
\end{equation}%
\textit{\noindent where "}$\func{div}\mathit{"}$\textit{\ denotes the
divergence operator in the }$x,$\textit{\ }$y$\textit{\ coordinates.}

\bigskip

Note that we can always solve in $P$ for equation (\ref{1.9.1})$.$ So the
real condition on the given data is (\ref{1.9}). We remark that in the $x,$ $%
y$ coordinates, $V$ is identified with $N^{\perp }$ whose integral curves
are the characteristic curves (compare (\ref{1.10}) with (2.21), (2.23) in 
\cite{CHMY04}). Note that we may consider (\ref{1.12.1}) as the (extrinsic)
Gauss-like equation in our surface theory. The Codazzi-like equation (\ref%
{1.9}) with $P$ $=$ $N$ (only involving the derivatives in the "intrinsic"
direction $V)$ can be deduced from (\ref{1.12.2}) together with (\ref{1.12.1}%
) through taking the derivative of (\ref{1.11}) in the direction $V$ and
applying (\ref{1.9.1}) to $\theta $ (also compare with the proof of Theorem
A).

\bigskip

Finally we study the global property of the singular set through a Hopf-type
index theorem. Let $u$\ $\in $\ $C^{1}(\Omega )$\ be a weak solution to ($%
\ref{1.1})$\ with $\vec{F}$\ $\in $\ $C^{1}(\Omega )$\ and $H$\ $\in $\ $%
C^{0}(\Omega ).$ Suppose $\partial \Omega $ consists of finitely many
components $C_{j},$ $j$ $=$ $1,$ $2,$ $...l,$ where each $C_{j}$ is a $C^{1}$
smooth, simple closed curve. Assume the singular set $S_{\vec{F}}(u)$ $%
\subset $ $\Omega $ is compact (which implies that $S_{\vec{F}}(u)$ does not
touch the boundary $\partial \Omega )$ and the characteristic curves hit
each $C_{j}$ in the following pattern. For $q$ $\in $ $C_{j}$ except
finitely many points$,$ there is only one characteristic curve $L_{q}$
hitting $q$ transversally (meaning that the vector $N^{\perp }$ along $L_{q}$
and the tangent vector of $C_{j}$ at $q$ are independent). Let $p$ be one of
those exceptional points. Consider the line field (1-dimensional
distribution) defined by the tangent lines (the lines having the direction $%
\pm N^{\perp }(u))$ of the characteristic curves, denoted as $\mathcal{D}.$
Denote the restriction of $\mathcal{D}$ to a small neighborhood $U$ $=$ $%
B_{\varepsilon }(p)$ $\cap $ $\Omega $ of $p$ by $\mathcal{D}_{U}.$ Take
another copy of $U$ and the corresponding line field$,$ denoted as $%
U^{\prime }$ and $\mathcal{D}_{U}^{\prime },$ resp.. We glue $U^{\prime }$
with $U$ along the boundary $C_{j}$ to get a (two-sided) neighborhood $%
\tilde{U}$ of $p.$ Denote the line field on $\tilde{U}$ obtained from $%
\mathcal{D}_{U}$ and $\mathcal{D}_{U}^{\prime }$ by $\mathcal{\tilde{D}}_{U}$%
. Define%
\begin{equation}
index(p,\mathcal{D}_{U})=\frac{1}{2}index(p,\mathcal{\tilde{D}}_{U})
\label{1.12.3}
\end{equation}

\noindent where $index(p,\mathcal{\tilde{D}}_{U})$ is the index of $p$ with
respect to the line field $\mathcal{\tilde{D}}_{U}$ (smoothing it near $%
\tilde{U}$ $\cap $ $C_{j}$ while keeping the topological type of $\mathcal{%
\tilde{D}}_{U})$ (see p.325 in \cite{Sp})$.$ Note that $index(p,\mathcal{D}%
_{U})$ is independent of the choice of small neighborhoods $U.$ See Example
7.1 in Section 7.

Let $p_{1},$..., $p_{m}$ denote those exceptional points of $C_{j}.$ Denote
the restriction of $\mathcal{D}$ to a small neighborhood $U_{k}$ of $p_{k}$
by $\mathcal{D}_{U_{k}}.$ We define the index of $C_{j}$ with respect to $u$
as follows:%
\begin{equation}
index(C_{j};u):=\dsum\limits_{k=1}^{m}index(p_{k},\mathcal{D}_{U_{k}}).
\label{1.12.4}
\end{equation}

Let $\chi (\Omega )$ denote the Euler characteristic number of $\Omega .$ We
can now formulate a Hopf-type index theorem.

\bigskip

\textbf{Theorem I.} \textit{Let }$\Omega $\textit{\ be a bounded domain of }$%
R^{2}.$\textit{\ Let }$u$\textit{\ }$\in $\textit{\ }$C^{1}(\Omega )$\textit{%
\ be a weak solution to (}$\ref{1.1})$\textit{\ with }$\vec{F}$\textit{\ }$%
\in $\textit{\ }$C^{1}(\Omega )$\textit{\ and }$H$\textit{\ }$\in $\textit{\ 
}$C^{1}(\Omega ).$\textit{\ Assume further }$N^{\perp }(curl\vec{F})$\textit{%
\ and }$N(H)$\textit{\ exist and are continuous (extended over singular
points) in }$\Omega .$\textit{\ Suppose }$\func{curl}\vec{F}$\textit{\ }$%
\neq $\textit{\ }$0$\textit{\ and }$\partial \Omega $\textit{\ consists of
finitely many components }$C_{j},$\textit{\ }$j$\textit{\ }$=$\textit{\ }$1,$%
\textit{\ }$2,$\textit{\ }$...l,$\textit{\ where each }$C_{j}$\textit{\ is a 
}$C^{1}$\textit{\ smooth, simple closed curve. Assume the singular set }$S_{%
\vec{F}}(u)$\textit{\ }$\subset $\textit{\ }$\Omega $\textit{\ is compact
and the characteristic curves hit each }$C_{j}$\textit{\ in the pattern
mentioned above. Then we have}%
\begin{equation}
\chi (\Omega )=\text{\# }\pi _{0}(S_{\vec{F}}(u))+%
\sum_{j=1}^{l}index(C_{j};u)  \label{1.13}
\end{equation}%
\textit{\noindent where \# }$\pi _{0}(S_{\vec{F}}(u))$\textit{\ denotes the
number of connected components of }$S_{\vec{F}}(u).$

\bigskip

For $u$ $\in $ $C^{1}(\bar{\Omega})$ we denote the set of singular points in 
$\bar{\Omega}$ ($\supset $ $\partial \Omega $ in particular$)$ still by $S_{%
\vec{F}}(u).$ Let $S(u)$ $:=$ $S_{\vec{F}}(u)$ for $\vec{F}$ $=$ $(-y,x).$

\bigskip

\textbf{Corollary J}. \textit{Let }$\Omega $\textit{\ be a bounded domain of 
}$R^{2}$ \textit{with }$C^{1}$\textit{\ smooth boundary. Consider a }$p$%
\textit{-minimal graph over }$\bar{\Omega},$\textit{\ defined by }$u$\textit{%
\ }$\in $\textit{\ }$C^{1}(\bar{\Omega}).$\textit{\ Suppose }$\Omega $%
\textit{\ is convex and }$S(u)$\textit{\ }$\subset \subset $\textit{\ }$%
\Omega $\textit{\ and nonempty. Then }$\#$\textit{\ }$\pi _{0}(S(u))$\textit{%
\ }$=$\textit{\ }$1.$

\bigskip

We remark that Corollary J is not the sharpest version for a convex domain
to have only one connected component of the singular set. But we won't
pursue it in this paper.

For a compact (connected) surface $\Sigma $ with no boundary, we would like
to know the configuration of its singular set $S_{\Sigma }$. When $\Sigma $
is $C^{2}$ smoothly immersed in a 3-dimensional pseudohermitian manifold
with bounded $p$-mean curvature, we learned from \cite{CHMY04} that $%
S_{\Sigma }$ consists of isolated (singular) points and closed $C^{1}$
curves. By the $C^{2}$ theory the characteristic curves meet at any singular
curve having the same tangent line , so the singular curves have no index
contribution with respect to the line field associated to the characteristic
curves. It follows that the Euler characteristic number $\chi (\Sigma )$
equals the number of isolated singular points. Therefore the genus $g(\Sigma
)$ of $\Sigma $ can only be zero or one (see Theorem E in \cite{CHMY04}) \
If we release the regularity condition, we wonder if $g(\Sigma )$ can be $%
\geq $ $2,$ say, for $\Sigma $ being $C^{1}$ smooth and of bounded $p$-mean
curvature.

\bigskip

\textbf{Acknowledgements}: J.-H. C. (P.Y., resp.) would like to thank SISSA
and Princeton University (Academia Sinica in Taiwan, resp.) for the kind
hospitality. This work has been supported by MiUR, under the project
Internazionalizzazione 2005: Equazioni alle derivate parziali nello studio
di problemi di analisi geometrica. A.M. has been supported by the project
FIRB-IDEAS "Analysis and Beyond"\emph{\ }and by the \emph{Giorgio and Elena
Petronio Fellowship} while visiting IAS in Princeton in the Fall Semester
2008-2009. He also would like to thank Academia Sinica in Taiwan and
Princeton Univesity for the kind hospitality.

\bigskip

\section{A Codazzi-like equation and properties: proofs of Theorem A and
Theorem B}

We first introduce the $s,$ $t$ coordinates. Let $N$ be a $C^{0}$ vector
field with $|N|$ $\equiv $ $1$ on a domain $\Omega $ $\subset $ $R^{2}.$ A
system of $C^{1}$ smooth local coordinates $s,$ $t$ is called a system of
characteristic coordinates if $s$ and $t$ have the property that $\nabla s$ $%
\parallel $ $N^{\perp }$ and $\nabla t$ $\parallel $ $N,$ i.e., $\nabla s$
and $\nabla t$ are parallel to $N^{\perp }$ and $N$, resp.$.$ It follows
that $\{t$ $=$ constant$\}$ are characteristic curves while $\{s$ $=$
constant$\}$ are seed curves (which are the integral curves of $N$). In \cite%
{CHY2} we proved the existence and studied the properties of such a system
of special coordinates under some mild conditions.

We now start to prove Theorem A. By Theorem C in \cite{CHY2} (since $u$%
\textit{\ }$\in $ $C^{1}(\Omega )$\textit{\ }is a weak solution to ($\ref%
{1.1})$\ with $\vec{F}$ $\in $ $C^{1}(\Omega )$ and $H$\ $\in $\ $%
C^{0}(\Omega )$ such that $\Omega $\ is nonsingular), we can find local
(near the point concerned) characteristic coordinates $s,$ $t$ and local
positive continuous functions $f,$ $g$ with the property that $Nf$ and $%
N^{\perp }g$ exist and are continuous, so that%
\begin{equation}
\frac{\partial }{\partial s}=\frac{1}{f}N^{\perp },\text{ \ \ }\frac{%
\partial }{\partial t}=\frac{1}{gD}N  \label{2.4}
\end{equation}

\noindent and%
\begin{equation}
Nf+fH=0,N^{\perp }g+\frac{(curl\vec{F})g}{D}=0  \label{2.3}
\end{equation}

\noindent (see (1.9) and (1.10) in \cite{CHY2}, resp.). From Theorem D in 
\cite{CHY2} we learn that $D^{\prime }$ ($:=$ $N^{\perp }D)$ exists and is
continuous. Moreover, $\theta $ $\in $ $C^{1}$ and (1.13) and (1.12) in \cite%
{CHY2} read (note that in \cite{CHY2} we used $rot\vec{F}$ instead of $curl%
\vec{F})$%
\begin{equation}
\theta _{t}=\frac{1}{gD^{2}}(curl\vec{F}-D^{\prime })  \label{2.1}
\end{equation}

\noindent and%
\begin{equation}
\theta _{s}=-\frac{H}{f}.  \label{2.2}
\end{equation}

\noindent Since $N(H)$ and $Nf$ exist and are continuous, ($\theta _{s})_{t}$
exists and is continuous in view of (\ref{2.4}). From (\ref{2.2}) and (\ref%
{2.4}) we compute%
\begin{eqnarray}
(\theta _{s})_{t} &=&-\frac{1}{gD}N(\frac{H}{f})  \label{2.5} \\
&=&-\frac{1}{gD}(\frac{N(H)}{f}-\frac{H(Nf)}{f^{2}})  \notag \\
&=&-\frac{1}{fgD}(N(H)+H^{2})\text{ \ (by (\ref{2.3})).}  \notag
\end{eqnarray}

The fact that ($\theta _{s})_{t}$ exists and is continuous implies that $%
(\theta _{t})_{s}$ exists and equals ($\theta _{s})_{t}$ (hence is
continuous) by a fundamental result in calculus (see Lemma 5.4 in \cite{CHY2}%
). It follows that $D^{\prime \prime }$ exists and is continuous in view of (%
\ref{2.1}) and (\ref{2.4}) since $N^{\perp }(curl\vec{F})$ exists and is
continuous by assumption.\ From (\ref{2.4}) and (\ref{2.1}) we compute%
\begin{eqnarray}
(\theta _{t})_{s} &=&\frac{1}{f}N^{\perp }(\theta _{t})  \label{2.6} \\
&=&\frac{1}{fgD^{2}}(N^{\perp }(curl\vec{F})-D^{\prime \prime })  \notag \\
&&+\frac{1}{f}(curl\vec{F}-D^{\prime })(-\frac{N^{\perp }g}{g^{2}D^{2}}-%
\frac{2D^{\prime }}{gD^{3}})  \notag \\
&=&\frac{1}{fgD}[\frac{N^{\perp }(curl\vec{F})-D^{\prime \prime }}{D}+\frac{%
(curl\vec{F}-D^{\prime })(curl\vec{F}-2D^{\prime })}{D^{2}}].  \notag
\end{eqnarray}

\noindent Here we have used (\ref{2.3}) in the last equality of (\ref{2.6}).
Finally by equating (\ref{2.5}) and (\ref{2.6}) we obtain (\ref{1.2}). We
have proved Theorem A.

\bigskip

We are going to prove Theorem B. Let $l$ $=$ $|\frac{curl\vec{F}(p)}{2}|$ $>$
$0.$ Let $m$ $=$ $\frac{curl\vec{F}(p)}{2}$ ($curl\vec{F}(p),$ resp.) if $%
curl\vec{F}(p)$ $>$ $0$ ($curl\vec{F}(p)$ $<$ $0$, resp.$)$. Let $M$ $=$ $%
curl\vec{F}(p)$ ($\frac{curl\vec{F}(p)}{2},$ resp.) if $curl\vec{F}(p)$ $>$ $%
0$ ($curl\vec{F}(p)$ $<$ $0,$ resp.). So $m$ $<$ $M.$ From equation (\ref%
{1.2}) we establish the following statement%
\begin{eqnarray}
&\text{\textit{Given} 0} \text{\TEXTsymbol{<}}\delta \text{\TEXTsymbol{<}}%
\frac{l}{3}\mathit{,}\text{ \textit{there exists} 0\TEXTsymbol{<}}%
\varepsilon \text{=}\varepsilon \text{(}\delta \text{)\TEXTsymbol{<}}\bar{%
\rho}\text{ \textit{such that} }&  \label{2.6.1} \\
&\text{\textit{for any }0} \text{\TEXTsymbol{<}}\rho \text{\TEXTsymbol{<}}%
\varepsilon \text{(}\delta \text{)}&  \notag \\
&\text{\textit{if} }D^{\prime }\text{(}\Gamma \text{(}\rho \text{))} \in 
\text{(-}\infty \text{,}m\text{-}\delta \text{)}\cup \text{(}M\text{+}\delta 
\text{,}\infty \text{)},\text{\textit{then} }D^{\prime \prime }(\Gamma (\rho
))>0\text{ \textit{while} }&  \notag \\
&\text{ \textit{if} }D^{\prime }\text{(}\Gamma \text{(}\rho \text{))} \in 
\text{(}m\text{+}\delta \text{,}M\text{-}\delta \text{)},\text{ \textit{then}
}D^{\prime \prime }(\Gamma (\rho ))<0.&  \notag
\end{eqnarray}

Next we claim that 
\begin{eqnarray}
\text{\textit{For any} }a &\in &(-\infty ,m-\delta )\cup (m+\delta ,M-\delta
)\cup (M+\delta ,\infty ),\text{ \textit{there exists}}  \label{2.6.2} \\
\text{\textit{at most one} }\rho &\in &(0,\varepsilon (\delta ))\text{ 
\textit{such that} }D^{\prime }(\Gamma (\rho ))=a.  \notag
\end{eqnarray}

\noindent Suppose there are $0$ $<$ $\rho _{1}$ $<$ $\rho _{2}$ $<$ $%
\varepsilon (\delta )$ such that $D^{\prime }(\rho _{1})$ $=$ $D^{\prime
}(\rho _{2})$ $=$ $a.$ Then there exist $\rho _{3},$ $\rho _{4},$ $\rho _{1}$
$\leq $ $\rho _{3}$ $<$ $\rho _{4}$ $\leq $ $\rho _{2},$ such that $%
D^{\prime }(\Gamma (\rho _{3}))$ $=$ $D^{\prime }(\Gamma (\rho _{1}))$ $=$ $%
D^{\prime }(\Gamma (\rho _{4}))$ $=$ $D^{\prime }(\Gamma (\rho _{2}))$ $=$ $%
a,$ and either $D^{\prime }(\Gamma (\rho ))$ $\geq $ $a$ for all $\rho $ $%
\in $ $[\rho _{3},$ $\rho _{4}]$ or $D^{\prime }(\Gamma (\rho ))$ $\leq $ $a$
for all $\rho $ $\in $ $[\rho _{3},$ $\rho _{4}].$ In both cases, we have
either $D^{\prime \prime }(\Gamma (\rho _{3}))$ $\geq $ $0$ while $D^{\prime
\prime }(\Gamma (\rho _{4}))$ $\leq $ $0$ or $D^{\prime \prime }(\Gamma
(\rho _{3}))$ $\leq $ $0$ while $D^{\prime \prime }(\Gamma (\rho _{4}))$ $%
\geq $ $0$. This contradicts (\ref{2.6.1}). We have shown (\ref{2.6.2}). (%
\ref{2.6.2}) will be used to show that $\lim_{\rho \rightarrow 0}D^{\prime
}(\Gamma (\rho ))$ exists.

Let $a_{1}$ $=$ $\lim \inf_{\rho \rightarrow 0^{+}}D^{\prime }(\Gamma (\rho
))$ and $a_{2}$ $=$ $\lim \sup_{\rho \rightarrow 0^{+}}D^{\prime }(\Gamma
(\rho )).$ Suppose $a_{1}$ $<$ $a_{2}.$ Then there exists $a$ $\in $ ($a_{1}$%
, $a_{2})$ and $a$ $\neq $ $m,$ $M.$ We can then choose small $\delta $ $\in 
$ $(0,$ $\frac{l}{3})$ such that $a$ $\in $ $(-\infty ,m-\delta )$ $\cup $ $%
(m+\delta ,M-\delta )$ $\cup $ $(M+\delta ,\infty ).$ By (\ref{2.6.2}) we
have at most one $\rho $ $\in $ $(0,\varepsilon (\delta ))$ such that $%
D^{\prime }(\Gamma (\rho ))$ $=$ $a.$ On the other hand, there are
infinitely many $\rho _{j}$ $\rightarrow $ $0$ satisfying $D^{\prime
}(\Gamma (\rho _{j}))$ $=$ $a$ by the continuity of $D^{\prime }\circ \Gamma 
$ and the definition of $a_{1}$ and $a_{2}.$ The contradiction implies $%
a_{1} $ $=$ $a_{2},$ and hence $\lim_{\rho \rightarrow 0}D^{\prime }(\Gamma
(\rho ))$ exists.

Let $\bar{a}$ $=$ $\lim_{\rho \rightarrow 0}D^{\prime }(\Gamma (\rho )).$
Note that $\bar{a}$ may be $\pm \infty .$ Suppose $\bar{a}$ $\neq $ $m,$ $M.$
Then we can find small $\delta $ $\in $ $(0,$ $\frac{l}{3})$ and associated $%
\varepsilon $ such that $D^{\prime }(\Gamma (\rho ))$ is monotonically
increasing or decreasing for $\rho \in (0,\varepsilon )$ according to $%
N^{\perp }$ $=$ $\pm \frac{\partial }{\partial \rho }$ ($N^{\perp }$ $=$ $%
\mp \frac{\partial }{\partial \rho },$ resp.) by (\ref{2.6.1}) in case $\bar{%
a}$ $\in $ $[-\infty ,m-\delta )\cup (M+\delta ,\infty ]$ ($\bar{a}$ $\in $ $%
(m+\delta ,M-\delta ),$ resp.). We then make estimate from (\ref{1.2}) and
integrate the resulting inequality to reach a contradiction.

In case $\bar{a}$ $\in $ $[-\infty ,m-\delta )\cup (M+\delta ,\infty ],$
there exists a constant $C_{1}$ $>$ $0$ such that $DD^{\prime \prime }$ $%
\geq $ $C_{1}$ for $\rho $ $\in $ $(0,\varepsilon ).$ Multiplying this
inequality by $\frac{D^{\prime }}{D},$ we obtain%
\begin{equation}
\frac{1}{2}[(D^{\prime })^{2}]^{\prime }=D^{\prime }D^{\prime \prime }\geq
C_{1}\frac{D^{\prime }}{D}\text{ (}\leq C_{1}\frac{D^{\prime }}{D},\text{
resp.)}  \label{2.6.3}
\end{equation}

\noindent if $N^{\perp }$ $=$ $\frac{\partial }{\partial \rho }$ ($N^{\perp
} $ $=$ $-\frac{\partial }{\partial \rho },$ resp.). Note that if $N^{\perp
} $ $=$ $\frac{\partial }{\partial \rho },$ we have%
\begin{eqnarray*}
D^{\prime } &:&=N^{\perp }D=\frac{\partial }{\partial \rho }D \\
&=&\lim_{\rho \rightarrow 0^{+}}\frac{D(\Gamma (\rho ))-D(\Gamma (0))}{\rho }%
=\lim_{\rho \rightarrow 0^{+}}\frac{D(\Gamma (\rho ))}{\rho }\geq 0
\end{eqnarray*}%
\noindent at $p$ $=$ $\Gamma (0),$ where $D$ $\geq $ $0$ and $D(p)$ $=$ $0$
since $p$ is a singular point. Integrating (\ref{2.6.3}) (for both cases)
over $[\rho ,$ $\rho _{0}]$ $\subset $ $(0,\varepsilon )$ gives%
\begin{eqnarray}
&&\frac{1}{2}(D^{\prime })^{2}(\Gamma (\rho _{0}))-\frac{1}{2}(D^{\prime
})^{2}(\Gamma (\rho ))  \label{2.6.4} \\
&\geq &C_{1}[\log D(\Gamma (\rho _{0}))-\log D(\Gamma (\rho ))].  \notag
\end{eqnarray}

\noindent Letting $\rho \rightarrow 0$ in (\ref{2.6.4}), we reach a
contradiction since the left hand side of (\ref{2.6.4}) is bounded from
above while the right hand side goes to +$\infty $ in view of $\log D(\Gamma
(0))$ $=$ $\log D(p)$ $=$ $\log 0$ $=$ $-\infty .$

In case $\bar{a}$ $\in $ $(m+\delta ,M-\delta ),$ there exists a constant $%
C_{2}$ $>$ $0$ such that $DD^{\prime \prime }$ $\leq $ $-C_{2}$ for $\rho $ $%
\in $ $(0,\varepsilon ).$ Multiplying this inequality by $\frac{D^{\prime }}{%
D}$ and integrating as above give%
\begin{eqnarray}
&&\frac{1}{2}(D^{\prime })^{2}(\Gamma (\rho _{0}))-\frac{1}{2}(D^{\prime
})^{2}(\Gamma (\rho ))  \label{2.6.5} \\
&\leq &-C_{2}[\log D(\Gamma (\rho _{0}))-\log D(\Gamma (\rho ))].  \notag
\end{eqnarray}

\noindent Letting $\rho \rightarrow 0$ in (\ref{2.6.5}), we observe that the
left hand side is a finite number while the right hand side goes to $-\infty 
$, a contradiction. Altogether we can conclude that $\bar{a}$ $=$ $m$ or $M.$
We have proved (\ref{1.6}) and the statement (b) in Theorem B follows.

To prove (c) in Theorem B, we observe that $curl\vec{F}$ is continuous,
nonzero, and hence $curl\vec{F}(p)$ and $curl\vec{F}(q)$ have the same sign
if $p$ and $q$ are connected by a characteristic curve $\Gamma $. Now $%
N^{\perp }$ points in an inward (outward, resp.) direction of $\Gamma $ at
both $p$ and $q$ if $curl\vec{F}$ $>$ $0$ ($curl\vec{F}$ $<$ $0,$ resp.) at $%
p$ and $q.$ This contradicts the continuity of $N^{\perp }$ (on $\Gamma ).$
We have shown the nonexistence of $\Gamma $ connecting $p$ and $q,$ and
hence (c).

\bigskip

In the Appendix we generalize equation (\ref{1.2}) and prove a result
analogous to Theorem B (a).

\bigskip

\section{Local configuration of the singular set}

Let $\Omega $ be a domain of $R^{2}.$ Let $u$ $\in $ $C^{1}(\Omega )$ and $%
\vec{F}$ $=$ $(F_{1},$ $F_{2})$ $\in $ $C^{1}(\Omega ).$ Recall that a point 
$p$ $\in $ $\Omega $ is called singular if $\nabla u$ $+$ $\vec{F}$ $=$ $0$
at $p.$ Let $S_{\vec{F}}(u)$ denote the set of all singular points. Define $%
\vec{G}^{\perp }$ $:=$ $(G_{2},$ $-G_{1})$ for $\vec{G}$ $=$ $(G_{1},$ $%
G_{2}).$ Recall that $curl\vec{F}$ $:=$ $(F_{2})_{x}-(F_{1})_{y}.$

\bigskip

\textbf{Lemma 3.1}. \textit{Suppose }$curl\vec{F}$\textit{\ }$\neq $\textit{%
\ }$0$\textit{\ in }$\Omega .$\textit{\ Then }$S_{\vec{F}}(u)$\textit{\ is
nowhere dense in }$\Omega .$

\bigskip

\proof
First note that $S_{\vec{F}}(u)$ is closed. So if $S_{\vec{F}}(u)$\ is not
nowhere dense in\textit{\ }$\Omega $, then there is a point $p_{1}$ $\in $ $%
S_{\vec{F}}(u)$ such that $S_{\vec{F}}(u)$ contains $B_{r_{1}}(p_{1}),$ a
ball of center $p_{1}$ with radius $r_{1}$ $>$ $0.$ Take a sequence of $%
C^{\infty }$ smooth functions $u_{n}$ such that $u_{n}$ converges to $u$ in $%
C^{1}$ norm on the closure of $B_{r_{2}}(p_{1})$ for $0$ $<$ $r_{2}$ $<$ $%
r_{1}.$ Since $\nabla u$ $+$ $\vec{F}$ $=$ $0$ in $B_{r_{1}}(p_{1}),$ we
have ($\nu $ denotes the unit outer normal)%
\begin{eqnarray}
0 &=&\doint\limits_{\partial B_{r_{2}}(p_{1})}(\nabla u+\vec{F})^{\perp
}\cdot \nu  \label{3.1} \\
&=&\lim_{n\rightarrow \infty }\doint\limits_{\partial
B_{r_{2}}(p_{1})}(\nabla u_{n}+\vec{F})^{\perp }\cdot \nu  \notag \\
&=&\int_{B_{r_{2}}(p_{1})}\func{div}(\nabla u_{n}+\vec{F})^{\perp }\text{ \
(by the divergence theorem)}  \notag \\
&=&\int_{B_{r_{2}}(p_{1})}\func{div}\vec{F}^{\perp }\text{ \ (since }\func{%
div}(\nabla u_{n})^{\perp }=0)  \notag \\
&=&\int_{B_{r_{2}}(p_{1})}curl\vec{F}.  \notag
\end{eqnarray}

\noindent Since $curl\vec{F}$ is continuous on $B_{r_{2}}(p_{1})$, a
connected set, we must have either $curl\vec{F}$ $>$ $0$ in $%
B_{r_{2}}(p_{1}) $ or $curl\vec{F}$ $<$ $0$ in $B_{r_{2}}(p_{1}).$ This
contradicts (\ref{3.1}).

\endproof%

\bigskip

When $H$ is the $p($or $H)$-mean curvature (see (\ref{1.1})), $\vec{F}$ $=$ $%
(-y,$ $x),$ so $curl\vec{F}$ $=$ $2$ and hence Lemma 3.1 applies. For
simplicity we will only consider the case of $p($or $H)$-minimal graphs in
the following discussion. Since $H$ $=$ $0,$ the characteristic curves are
straight lines. We often call them characteristic lines (here line may just
mean line segment). Recall that by a domain we mean an open and connected
set. For a subset $A$ $\subset $ $R^{n}$ we define an $\varepsilon $%
-neighborhood $N_{\varepsilon }(A)$ by%
\begin{equation*}
N_{\varepsilon }(A):=\{p\in R^{n}\text{ }|\text{ }d(p,A)<\varepsilon \}
\end{equation*}

\noindent where $d(p,A)$ $:=$ $\inf \{d(p,q)$ $|$ $q\in A\}$ and $d(\cdot
,\cdot )$ is the Euclidean distance. For $A,$ $B$ $\subset $ $R^{n},$ we
define the Hausdorff distance between $A$ and $B$ to be the infimum of $%
\varepsilon $ $>$ $0$ such that $B$ $\subset $ $N_{\varepsilon }(A),$ $A$ $%
\subset $ $N_{\varepsilon }(B).$

\bigskip

\textbf{Lemma 3.2}. \textit{Consider a }$C^{1}$\textit{\ smooth }$p$\textit{%
-minimal graph defined by }$u$\textit{\ over a plane convex domain }$\Omega
. $\textit{\ Let }$\tilde{\Gamma}_{\infty }$\textit{\ be a straight line
such that }$\Gamma _{\infty }:=\tilde{\Gamma}_{\infty }$\textit{\ }$\cap $%
\textit{\ }$\Omega $\textit{\ divides }$\Omega $\textit{\ into two disjoint
nonempty domains }$\Omega ^{+},$\textit{\ }$\Omega ^{-}.$\textit{\ Let }$%
\tilde{\Gamma}_{j}$\textit{\ be a family of straight lines such that all }$%
\Gamma _{j}:=\tilde{\Gamma}_{j}$\textit{\ }$\cap $\textit{\ }$\Omega $%
\textit{\ are characteristic. Suppose }$\{\Gamma _{j}\}$\textit{\ converges
to }$\Gamma _{\infty }$ \textit{in the sense that the Hausdorff distance
between }$\Gamma _{j}$\textit{\ and }$\Gamma _{\infty }$\textit{\ tends to
zero as }$j$\textit{\ }$\rightarrow $\textit{\ }$\infty .$\textit{\ Then }$%
\Gamma _{\infty }$\textit{\ is a characteristic line (segment). In
particular, }$\Gamma _{\infty }$\textit{\ contains no singular point.}

\bigskip

\proof

We first observe that since the singular set $S(u)$ is closed, $\Gamma
_{\infty }\backslash S(u)$ is open in $\Gamma _{\infty }.$ So $\Gamma
_{\infty }\backslash S(u)$ is empty or the union of open line segments $%
\Gamma _{\infty }^{k},$ $k$ $=$ $1,$ $2,...$. Each $\Gamma _{\infty }^{k}$
is characteristic since it is the (Hausdorff) limit of characteristic lines.
If one of the $\Gamma _{\infty }^{k}$'s has two singular end points in $%
\Omega ,$ we reach a contradiction since we cannot have a characteristic
line connecting two singular points by Theorem B (c). Therefore we have only
four possibilities:

Case 1-- $\Gamma _{\infty }\backslash S(u)$ is empty, i.e., $\Gamma _{\infty
}$ $\subset $ $S(u);$

Case 2-- $\Gamma _{\infty }\backslash S(u)$ $=$ $\Gamma _{\infty }^{1};$

Case 3-- $\Gamma _{\infty }\backslash S(u)$ $=$ $\Gamma _{\infty }^{1}$ $%
\cup $ $\Gamma _{\infty }^{2};$

Case 4-- $\Gamma _{\infty }\backslash S(u)$ $=$ $\Gamma _{\infty },$ i.e., $%
\Gamma _{\infty }$ contains no singular points,

\noindent where each of $\Gamma _{\infty }^{1}$ and $\Gamma _{\infty }^{2}$
has only one singular end point in $\Omega .$ In Case 3, if there is only
one singular point $p$ $\in $ $\Gamma _{\infty }$, then we can decompose $%
\Gamma _{\infty }$ as a disjoint union $\Gamma _{\infty }^{+}$ $\cup $ $%
\{p\} $ $\cup $ $\Gamma _{\infty }^{-}$ where $\Gamma _{\infty }^{+}$ $=$ $%
\Gamma _{\infty }^{1}$ and $\Gamma _{\infty }^{-}$ $=$ $\Gamma _{\infty
}^{2} $ are characteristic rays emitted by $p$ in opposite directions
respectively (see Figure 3.1 below). We claim that this is impossible.

\begin{figure}[h]
\begin{center}
\includegraphics[height=4.5cm]{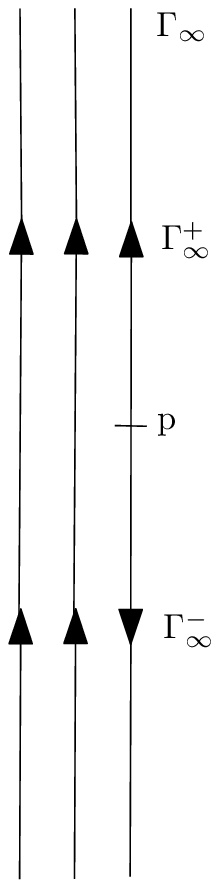}\\[0pt]
Figure 3.1
\end{center}
\par
\end{figure}

Let $B_{r}(p)$ $\subset $ $\Omega $ denote a ball of center $p$ with radius $%
r.$ Take $p_{j}^{+}$ $\in $ $\Gamma _{\infty }^{+}\cap B_{r}(p)$ and $%
p_{j}^{-}$ $\in \ \Gamma _{\infty }^{-}\cap B_{r}(p)$ approaching $p.$ There
exist a large integer $n(j)$ and $q_{j}^{+}$, $q_{j}^{-}$ $\in $ $\Gamma
_{n(j)}$ such that%
\begin{equation}
|N^{\perp }(q_{j}^{\pm })-N^{\perp }(p_{j}^{\pm })|<\frac{1}{j}  \label{3.2}
\end{equation}

\noindent by the continuity of $N^{\perp }.$ On the other hand, we have%
\begin{equation}
\lim_{p_{j}^{+}\rightarrow p}N^{\perp
}(p_{j}^{+})=-\lim_{p_{j}^{-}\rightarrow p}N^{\perp }(p_{j}^{-})  \label{3.3}
\end{equation}

\noindent\ according to Theorem B (b) in Section 1 (for the case of $\Gamma
_{\infty }$ being a straight line, $N^{\perp }(p_{j}^{+})$ $=$ $-N^{\perp
}(p_{j}^{-})$ considered as free vectors$).$ From (\ref{3.2}), we get%
\begin{eqnarray*}
\lim_{p_{j}^{+}\rightarrow p}N^{\perp }(p_{j}^{+}) &=&\lim_{j\rightarrow
\infty }N^{\perp }(q_{j}^{+}) \\
&=&\lim_{j\rightarrow \infty }N^{\perp
}(q_{j}^{-})=\lim_{p_{j}^{-}\rightarrow p}N^{\perp }(p_{j}^{-})
\end{eqnarray*}

\noindent which contradicts (\ref{3.3}) (for the second equality, we have
used $N^{\perp }(q_{j}^{+})$ $=$ $N^{\perp }(q_{j}^{-})$ since $\Gamma
_{n(j)}$ is a straight line). We have proved our claim. So in Case 3, the
remaining situation is that $\Gamma _{\infty }$ $=$ $\Gamma _{\infty }^{1}$ $%
\cup $ $I\cup $ $\Gamma _{\infty }^{2}$ (disjoint union) where $I$ $\subset $
$S(u)$ is a closed line segment. Case 1, Case 2, and this situation of Case
3 have the common feature that $\Gamma _{\infty }$ contains an open (and
hence a shorter closed) line segment which consists of singular points. Let $%
\Upsilon $ $\subset $ $\Gamma _{\infty }$ denote such a closed singular line
segment. We are going to show that this is impossible.

Near $\Upsilon ,$ we can parametrize $\Gamma _{\infty }$ ($\Gamma _{j},$
resp.) by the map $\gamma $ ($\gamma _{j},$ resp.): $(-a,a)$ $\rightarrow $ $%
\Gamma _{\infty }$ ($\Gamma _{j},$ resp.), where $s$ $\in $ $(-a,a)$ is the
unit-speed parameter with $\frac{\partial }{\partial s}$ $=$ $N^{\perp }$ on 
$\Gamma _{j}$, such that $\gamma ([-\varepsilon ,\varepsilon ])$ $\subset $ $%
\Upsilon $ for some $0$ $<$ $\varepsilon $ $<$ $a$ and $\gamma _{j}(s)$ $%
\rightarrow $ $\gamma (s)$ as $j$ $\rightarrow $ $\infty $ for $s$ $\in $ $%
[-\varepsilon ,\varepsilon ]$. We claim that if there is a point $p_{0}$ $%
\in $ $\Gamma _{j}$ with $D^{\prime }(p_{0})$ $>$ $2$ ($D^{\prime }(p_{0})$ $%
=$ $2,$ $1$ $<$ $D^{\prime }(p_{0})$ $<$ $2,$ $D^{\prime }(p_{0})$ $=$ $1,$ $%
D^{\prime }(p_{0})$ $<$ $1,$ resp.), then $D^{\prime }(p)$ $>$ $2$ ($%
D^{\prime }(p)$ $=$ $2,$ $1$ $<$ $D^{\prime }(p)$ $<$ $2,$ $D^{\prime }(p)$ $%
=$ $1,$ $D^{\prime }(p)$ $<$ $1,$ resp.) for all $p$ $\in $ $\Gamma _{j}.$
First observe that $D^{\prime }$ $\equiv $ $2$ ($D^{\prime }$ $\equiv $ $1,$
resp$.)$ on $\Gamma _{j}$ is a solution to (\ref{1.4}). Our claim follows
from the uniqueness of solutions to (\ref{1.4}) with initial data $%
(D(p_{0}), $ $D^{\prime }(p_{0}))$. By the pigeonhole principle, at least
one of the five cases: $D^{\prime }$ $>$ $2,$ $D^{\prime }$ $=$ $2,$ $1$ $<$ 
$D^{\prime }$ $<$ $2,$ $D^{\prime }$ $=$ $1,$ $D^{\prime }$ $<$ $1$ holds on 
$\Gamma _{j}$ for infinitely many $j^{\prime }s.$ Suppose that $D^{\prime }$ 
$>$ $2$ on $\Gamma _{j}$ for infinitely many $j^{\prime }s.$ Then we have%
\begin{equation}
D(\gamma _{j}(\varepsilon ))-D(\gamma _{j}(0))\geq 2\varepsilon  \label{3.4}
\end{equation}

\noindent for infinitely many $j^{\prime }s.$ Letting $j$ $\rightarrow $ $%
\infty $ in (\ref{3.4}) and noting that $D$ $=$ $0$ on $\Upsilon $ give 
\begin{equation*}
0=D(\gamma (\varepsilon ))-D(\gamma (0)\geq 2\varepsilon ,
\end{equation*}

\noindent a contradiction. Similarly we can reach a contradiction for the
cases $D^{\prime }$ $=$ $2,$ $1$ $<$ $D^{\prime }$ $<$ $2,$ $D^{\prime }$ $=$
$1,$ resp.$.$ Now we are left to deal with the remaining case: $D^{\prime }$ 
$<$ $1$ on infinitely many $\Gamma _{j}$'$s$ (still denoted as $\Gamma
_{j}). $ From (\ref{1.4}) we learn that 
\begin{equation*}
D^{\prime \prime }=\frac{2(D^{\prime }-1)(D^{\prime }-2)}{D}>0.
\end{equation*}%
\noindent So $D^{\prime }$ is strictly increasing on each $\Gamma _{j}$.
Hence one of the three cases: $D^{\prime }(\gamma _{j}(-\frac{\varepsilon }{3%
}))$ $\leq $ $-\frac{1}{2},$ $D^{\prime }(\gamma _{j}(\frac{\varepsilon }{3}%
))$ $\geq $ $\frac{1}{2},$ $-\frac{1}{2}$ $\leq $ $D^{\prime }$ $\leq $ $%
\frac{1}{2}$ on $\gamma _{j}([-\frac{\varepsilon }{3},\frac{\varepsilon }{3}%
])$ holds for infinitely many $j^{\prime }s$ by the pigeonhole principle$.$
Suppose that $D^{\prime }(\gamma _{j}(-\frac{\varepsilon }{3}))$ $\leq $ $-%
\frac{1}{2}$ for infinitely many $j^{\prime }s.$ Then $D^{\prime }$ $\leq $ $%
-\frac{1}{2}$ on $\gamma _{j}([-\frac{2\varepsilon }{3},$ $-\frac{%
\varepsilon }{3}])$ since $D^{\prime }$ is strictly increasing on each $%
\Gamma _{j}.$ It follows that%
\begin{equation*}
D(\gamma _{j}(-\frac{2\varepsilon }{3}))-D(\gamma _{j}(-\frac{\varepsilon }{3%
}))\geq \frac{1}{2}\cdot \frac{\varepsilon }{3}.
\end{equation*}

\noindent Taking $j$ $\rightarrow $ $\infty ,$ we obtain that $0$ $=$ $%
D(\gamma (-\frac{2\varepsilon }{3}))$ $-$ $D(\gamma (-\frac{\varepsilon }{3}%
))$ $\geq $ $\frac{\varepsilon }{6},$ a contradiction. Next for the case $%
D^{\prime }(\gamma _{j}(\frac{\varepsilon }{3}))$ $\geq $ $\frac{1}{2},$ we
have $D^{\prime }$ $\geq $ $\frac{1}{2}$ on $\gamma _{j}([\frac{\varepsilon 
}{3},$ $\frac{2\varepsilon }{3}])$ since $D^{\prime }$ is strictly
increasing on each $\Gamma _{j}.$ A similar argument as above will lead to a
contradiction. Now let us assume that $-\frac{1}{2}$ $\leq $ $D^{\prime }$ $%
\leq $ $\frac{1}{2}$ on $\gamma _{j}([-\frac{\varepsilon }{3},\frac{%
\varepsilon }{3}])$ for infinitely many $j^{\prime }s.$ From (\ref{1.4}) we
have%
\begin{eqnarray}
D^{\prime \prime } &=&\frac{2(D^{\prime }-1)(D^{\prime }-2)}{D}  \label{3.5}
\\
&\geq &\frac{3}{2}\cdot \frac{1}{D}  \notag
\end{eqnarray}

\noindent on $\gamma _{j}([-\frac{\varepsilon }{3},\frac{\varepsilon }{3}]).$
On the other hand, we have the Taylor expansion to the second order for $D:$%
\begin{eqnarray}
&&D(\gamma _{j}(\frac{\varepsilon }{3}))-D(\gamma _{j}(-\frac{\varepsilon }{3%
}))  \label{3.6} \\
&=&D^{\prime }(\gamma _{j}(-\frac{\varepsilon }{3}))\frac{2\varepsilon }{3}+%
\frac{D^{\prime \prime }(\gamma _{j}(\xi _{j}))}{2}(\frac{2\varepsilon }{3}%
)^{2}  \notag
\end{eqnarray}

\noindent for some $\xi _{j}$ between $-\frac{\varepsilon }{3}$ and $\frac{%
\varepsilon }{3}.$ From (\ref{3.6}), condition on $D^{\prime },$ and (\ref%
{3.5}), we have%
\begin{eqnarray}
&&D(\gamma _{j}(\frac{\varepsilon }{3}))-D(\gamma _{j}(-\frac{\varepsilon }{3%
}))  \label{3.7} \\
&\geq &-\frac{1}{2}\cdot \frac{2\varepsilon }{3}+\frac{3}{4}\cdot \frac{1}{%
D(\gamma _{j}(\xi _{j}))}\cdot (\frac{2\varepsilon }{3})^{2}.  \notag
\end{eqnarray}

\noindent Letting $j$ $\rightarrow $ $\infty $ in (\ref{3.7}) and noting
that $D(\gamma _{j}(\xi _{j}))$ $\rightarrow $ $0$, we obtain that $0$ $=$ $%
D(\gamma (\frac{\varepsilon }{3}))$ $-$ $D(\gamma (-\frac{\varepsilon }{3}))$
$\geq $ $\infty ,$ a contradiction. Altogether we have excluded Cases 1, 2,
and 3. But Case 4 just means that $\Gamma _{\infty }$ contains no singular
point and hence it is a characteristic line.

\endproof%

\bigskip

We can extend Lemma 3.2 as follows. Although Lemma $3.2^{\prime }$ is more
general than Lemma 3.2, we include the above proof of Lemma 3.2 for the
reader to understand the situation better.

\bigskip

\textbf{Lemma }$\mathbf{3.2}^{\prime }$\textbf{.} \textit{Let }$\Omega $%
\textit{\ be a bounded domain of }$R^{2}.$ \textit{Let }$u$\textit{\ }$\in $%
\textit{\ }$C^{1}(\Omega )$\textit{\ be a weak solution to (}$\ref{1.1})$%
\textit{\ with }$\vec{F}$\textit{\ }$\in $\textit{\ }$C^{1}(\Omega )$\textit{%
\ and }$H$\textit{\ }$\in $\textit{\ }$C^{0}(\Omega ).$\textit{\ Assume
further }$N^{\perp }(curl\vec{F})$\textit{\ and }$N(H)$\textit{\ exist and
are continuous (extended over singular points) in }$\Omega .$\textit{\
Suppose }$\func{curl}\vec{F}$ $\neq $ $0.$\textit{\ Take }$p_{0}$ $\in $ $%
\Omega ,$ $B_{r_{1}}(p_{0})$ $\subset \subset $ $\Omega $ \textit{such that }%
$0$ $<$ $r_{1}$\textit{\ }$\leq $ $(\sup_{B_{r_{1}}(p_{0})}|H|)^{-1}.$%
\textit{\ Take a sequence of nonsingular points }$p_{j}$\textit{\ }$\in $%
\textit{\ }$B_{r_{1}}(p_{0})$\textit{\ converging to }$p_{0}.$\textit{\
Suppose for each }$j,$\textit{\ the characteristic curve passing through }$%
p_{j}$\textit{\ does not hit }$p_{0}$\textit{\ or any singular points in }$%
B_{r_{1}}(p_{0})$\textit{\ before it meets }$\partial B_{r_{1}}(p_{0})$%
\textit{\ at two points }$q_{j}^{1},$\textit{\ }$q_{j}^{2}.$\textit{\ Let }$%
\tilde{\Gamma}_{j}$\textit{\ denote the characteristic curve passing through 
}$p_{j}$\textit{\ with }$q_{j}^{1},$\textit{\ }$q_{j}^{2}$\textit{\ as two
end points}$.$\textit{\ Then}

\textit{(a) There exists }$0$\textit{\ }$<$\textit{\ }$r_{2}$\textit{\ }$<$%
\textit{\ }$r_{1}$\textit{\ such that a subsequence of closed arcs }$\Gamma
_{j}$\textit{\ }$\subset $\textit{\ }$\tilde{\Gamma}_{j}$\textit{\ }$\cap $%
\textit{\ }$B_{r_{2}}(p_{0}),$\textit{\ containing }$p_{j}$ \textit{and an
open arc}$,$ \textit{converges to a closed arc }$\Gamma _{\infty }$\textit{\ 
}$\subset $\textit{\ }$B_{r_{2}}(p_{0})$\textit{\ in }$C^{2}$ \textit{with
respect to a certain parametrization.}

\textit{(b) }$\Gamma _{\infty }$ \textit{contains }$p_{0}$\textit{\ and an
open arc, but\ contains no singular points. Moreover, }$\Gamma _{\infty }$ 
\textit{is a characteristic curve passing through} $p_{0}.$

\bigskip

\proof
Take $0$ $<$ $r_{2}$\textit{\ }$<<$\textit{\ }$r_{1}$ ($\leq $ $%
(\sup_{B_{r_{1}}(p_{0})}|H|)^{-1}$ by assumption) such that $\Gamma _{j}^{0}$%
\textit{\ }$:=$\textit{\ }$\tilde{\Gamma}_{j}$\textit{\ }$\cap $\textit{\ }$%
B_{r_{2}}(p_{0})$\textit{\ }is a connected arc passing through $p_{j}$ (note
that -$H$ is the curvature of the curve $\tilde{\Gamma}_{j}$ \cite{CHY2}$)$
(see Figure $3.1^{\prime}$ below).

\begin{figure}[h]
\begin{center}
\includegraphics[height=4.5cm]{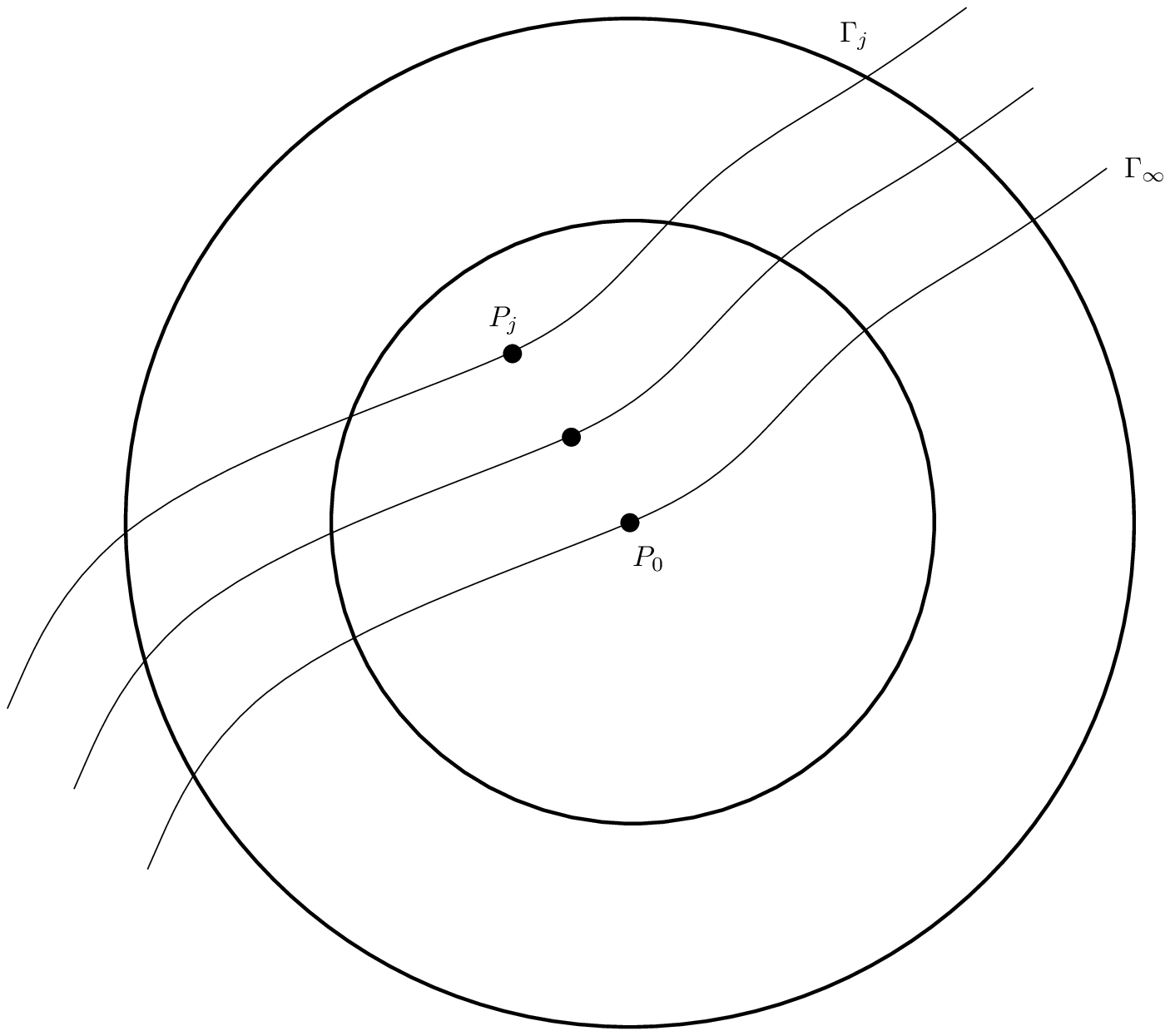}\\[0pt]
Figure $3.1^{\prime}$%
\end{center}
\par
\end{figure}

For $j$ large, there exists $s_{0}$ (independent of $j)$ $>$ $0$ such that
near $p_{j}$ we can parametrize $\Gamma _{j}^{0}$ by the ($C^{2}$ smooth
according to \cite{CHY2}) map $\gamma _{j}$ $:$ $[-s_{0},s_{0}]$ $%
\rightarrow $ $\Gamma _{j}^{0}$ $\subset $ $B_{r_{2}}(p_{0})$, $\gamma
_{j}(0)$ $=$ $p_{j},$ where $s$ $\in $ $[-s_{0},s_{0}]$ is the unit-speed
parameter with $\frac{\partial }{\partial s}$ $=$ $N^{\perp }$ on $\Gamma
_{j}.$ We choose an angular function $\theta $ ranging in $[0,2\pi ),$ which
works for all $\Gamma _{j}^{0\prime }s.$ Let $\theta _{j}(s)$ $:=$ $\theta
(\gamma _{j}(s)).$ It follows that for $j$ large and $s$ $\in $ $%
[-s_{0},s_{0}],$ we have%
\begin{eqnarray}
|\theta _{j}(s)| &\leq &C_{1}  \label{3.8} \\
|\theta _{j}^{\prime }(s)| &=&|-H(\gamma _{j}(s)|\leq C_{2}  \notag
\end{eqnarray}

\noindent (Theorem A in \cite{CHY2}) where the constants $C_{1}$ and $C_{2}$
are independent of $j$ and $s.$ By the Arzela-Ascoli theorem in view of (\ref%
{3.8}), we can find a subsequence of $\theta _{j},$ still denoted $\theta
_{j},$ such that $\theta _{j}$ converges to $\theta _{\infty }$ in $%
C^{0}([-s_{0},s_{0}]).$ So \{$\theta _{j}\}$ is Cauchy. We claim that \{$%
\gamma _{j}\}$ is Cauchy (with respect to $C^{0}$-topology). Write 
\begin{eqnarray}
\gamma _{j}(s)-\gamma _{k}(s) &=&(\gamma _{j}(s)-\gamma _{j}(0))-(\gamma
_{k}(s)-\gamma _{k}(0))  \label{3.9} \\
&&+(\gamma _{j}(0)-\gamma _{k}(0)).  \notag
\end{eqnarray}

\noindent From the definition of characteristic curves, we have%
\begin{equation*}
\gamma _{j}(s)-\gamma _{j}(0)=\int_{0}^{s}(\sin \theta _{j}(\tau ),-\cos
\theta _{j}(\tau ))d\tau
\end{equation*}%
\noindent and hence 
\begin{eqnarray}
&&(\gamma _{j}(s)-\gamma _{j}(0))-(\gamma _{k}(s)-\gamma _{k}(0))
\label{3.10} \\
&=&\int_{0}^{s}[(\sin \theta _{j}(\tau )-\sin \theta _{k}(\tau ),-\cos
\theta _{j}(\tau )+\cos \theta _{k}(\tau ))]d\tau .  \notag
\end{eqnarray}

\noindent Since \{$\theta _{j}\}$ is Cauchy, we have $(\gamma _{j}(s)-\gamma
_{j}(0))-(\gamma _{k}(s)-\gamma _{k}(0))$ is small uniformly in $s$ $\in $ $%
[-s_{0},s_{0}]$ for $j,$ $k$ large enough by (\ref{3.10}). On the other
hand, $\gamma _{j}(0)$ $=$ $p_{j}$ $(\in $ $R^{2})$ converges to $p_{0},$ so 
$\{p_{j}\}$ is Cauchy and hence $\gamma _{j}(0)-\gamma _{k}(0)$ is small for 
$j,$ $k$ large enough. Altogether we have shown that \{$\gamma _{j}\}$ is
Cauchy in view of (\ref{3.9}). Thus $\gamma _{j}$ converges in $%
C^{0}([-s_{0},s_{0}],R^{2})$ and we denote the limit by $\gamma _{\infty }.$
So if we take $\Gamma _{j}$ in (a) to be $\gamma _{j}([-s_{0},s_{0}]),$ then 
$\Gamma _{\infty }$ $=$ $\gamma _{\infty }([-s_{0},s_{0}]).$ Since $\theta
_{j}$ converges to $\theta _{\infty }$ in $C^{0}([-s_{0},s_{0}])$, $\gamma
_{j}$ converges to $\gamma _{\infty }$ in $C^{1}([-s_{0},s_{0}],R^{2})$ with 
$\frac{d\gamma _{\infty }(s)}{ds}$ $=$ $(\sin \theta _{\infty }(s),$ $-\cos
\theta _{\infty }(s))$ as the following argument shows$.$ Write $\gamma
_{j}(s)$ $=$ $(x_{j}(s),$ $y_{j}(s))$ and $\gamma _{\infty }(s)$ $=$ $%
(x_{\infty }(s),$ $y_{\infty }(s)).$ From the mean value theorem we have%
\begin{equation*}
\frac{x_{j}(s_{2})-x_{j}(s_{1})}{s_{2}-s_{1}}=\sin \theta _{j}(\hat{s}),%
\frac{y_{j}(s_{2})-y_{j}(s_{1})}{s_{2}-s_{1}}=-\cos \theta _{j}(\bar{s})
\end{equation*}

\noindent for $s_{1}$ $<$ $\hat{s},$ $\bar{s}$ $<$ $s_{2}.$ Taking $%
j\rightarrow \infty $ in the above formulas, we obtain 
\begin{equation*}
\frac{x_{\infty }(s_{2})-x_{\infty }(s_{1})}{s_{2}-s_{1}}=\sin \theta
_{\infty }(\hat{s}),\frac{y_{\infty }(s_{2})-y_{\infty }(s_{1})}{s_{2}-s_{1}}%
=-\cos \theta _{\infty }(\bar{s}).
\end{equation*}

\noindent Then letting $s_{2}\rightarrow s_{1}$ and hence $\hat{s},$ $\bar{s}
$ $\rightarrow $ $s_{1},$ we get%
\begin{equation}
\frac{dx_{\infty }(s_{1})}{ds}=\sin \theta _{\infty }(s_{1}),\frac{%
dy_{\infty }(s_{1})}{ds}=-\cos \theta _{\infty }(s_{1}).  \label{3.11}
\end{equation}

\noindent Therefore $\gamma _{j}^{\prime }$ $=$ $(\sin \theta _{j},$ $-\cos
\theta _{j})$ converges in $C^{0}([-s_{0},s_{0}],R^{2})$ to $(\sin \theta
_{\infty },$ $-\cos \theta _{\infty })$ which equals $\gamma _{\infty
}^{\prime }$ by (\ref{3.11}). Let $H_{j}(s)$ $:=$ $H(\gamma _{j}(s))$ and $%
H_{\infty }(s)$ $:=$ $H(\gamma _{\infty }(s)).$ A similar argument with $%
x_{j},$ $\sin \theta _{j}$ replaced by $\theta _{j},$ $-H_{j},$ respectively
in the above argument shows that 
\begin{equation}
\frac{d\theta _{\infty }(s_{1})}{ds}=-H_{\infty }(s_{1})  \label{3.12}
\end{equation}

\noindent (noting that $\theta _{j}$ is $C^{1}$ smooth and $\theta
_{j}^{\prime }(s)$ $=$ $-H_{j}(s)$ from Theorem A in \cite{CHY2}$).$ Now we
compute%
\begin{eqnarray*}
\gamma _{j}^{^{\prime \prime }} &=&(\cos \theta _{j},\sin \theta _{j})\theta
_{j}^{\prime } \\
&=&(\cos \theta _{j},\sin \theta _{j})(-H_{j})\rightarrow (\cos \theta
_{\infty },\sin \theta _{\infty })(-H_{\infty })
\end{eqnarray*}%
\noindent in $C^{0}([-s_{0},s_{0}]),$ which equals $\gamma _{\infty
}^{\prime \prime }$ according to (\ref{3.11}), (\ref{3.12}). We have shown
that $\Gamma _{j}$ converges to $\Gamma _{\infty }$ in $C^{2}$ with respect
to the parametrization given by $\gamma _{j},$ $\gamma _{\infty },$
respectively. We have completed the proof of (a).

For the proof of (b), it is clear that $p_{0}$ $\in $ $\Gamma _{\infty }$
and $\Gamma _{\infty }$ contains the open arc $\gamma _{\infty
}((-s_{0},s_{0})).$ To show that $\Gamma _{\infty }$ contains no singular
points, we mimic the reasoning as in the proof of Lemma 3.2. In view of
Theorem B (b), (c) we only have to deal with and exclude (by contradiction)
the situation that $\Gamma _{\infty }$ contains an open (and hence a shorter
closed) arc which consists of singular points. Without loss of generality we
may just assume that $\Gamma _{\infty }$ is such a closed singular arc. That
is to say, $D$ $\equiv $ $0$ on $\Gamma _{\infty }.$ Let%
\begin{equation*}
\lambda :=\inf_{p\in B_{r_{2}}(p_{0})}\func{curl}\vec{F}(p).
\end{equation*}

\noindent From the assumption $\func{curl}\vec{F}$ $\neq $ $0$, we may
assume $\lambda $ $>$ $0$ without loss of generality. Let $\bar{\Gamma}_{j}$ 
$:=$ $\gamma _{j}([-\bar{s}_{0},\bar{s}_{0}])$ for $0$ $<$ $\bar{s}_{0}$ $<$ 
$s_{0}$, independent of $j.$ Let%
\begin{equation*}
m_{j}:=\inf_{p\in \bar{\Gamma}_{j}}D^{\prime }(p),\text{ }M_{j}:=\sup_{p\in 
\bar{\Gamma}_{j}}D^{\prime }(p).
\end{equation*}

\noindent We claim that for $j$ large enough, there holds $M_{j}$ $<$ $\frac{%
\lambda }{4}.$ Suppose the converse holds. Then there exists a subsequence $%
j_{k}$ such that $M_{j_{k}}$ $\geq $ $\frac{\lambda }{4}.$ Hence we can find
a sequence of points $q_{k}$ $\in $ $\bar{\Gamma}_{j_{k}}$ such that $%
D^{\prime }(q_{k})$ $\geq $ $\frac{\lambda }{8}.$ We then extract a
convergent subsequence of $q_{k},$ still denoted as $q_{k}.$ Let $q_{\infty
} $ $=$ $\lim q_{k}.$ It follows that $q_{\infty }$ $\in $ $\bar{\Gamma}%
_{\infty }$ $:=$ $\gamma _{\infty }([-\bar{s}_{0},\bar{s}_{0}])$ since $%
\gamma _{j}$ converges to $\gamma _{\infty }$ (in $C^{2}).$ Let $s_{k}$ $%
\leq $ $\bar{s}_{0}$ such that $\gamma _{j_{k}}(s_{k})$ $=$ $q_{k}.$ We are
going to show that for $k$ large, there holds

\begin{equation}
D^{\prime }(\gamma _{j_{k}}(s))\geq \frac{\lambda }{8}  \label{3.13}
\end{equation}%
\noindent for all $s$ $\in $ $[s_{k},s_{0}]$ ($\supset $ $[\bar{s}%
_{0},s_{0}])$. If not, there are a subsequence of $k$, still denoted as $k,$
and a sequence $t_{k}$ $\in $ $[s_{k},s_{0}]$ such that $D^{\prime }(\gamma
_{j_{k}}(t_{k}))$ $<$ $\frac{\lambda }{8}.$ May assume that $D^{\prime
}(\gamma _{j_{k}}(t_{k}))$ achieves its minimum over $[s_{k},s_{0}]$ at $%
t_{k}.$ From (\ref{1.2}) we evaluate%
\begin{equation}
D^{\prime \prime }=\frac{2(D^{\prime }-\frac{curl\vec{F}}{2})(D^{\prime
}-curl\vec{F})}{D}+(N^{\perp }(curl\vec{F}))+(H^{2}+N(H))D  \label{3.13'}
\end{equation}

\noindent at $\gamma _{j_{k}}(t_{k}).$ It is easy to see that $D^{\prime
\prime }(\gamma _{j_{k}}(t_{k}))$ $>$ $0$ for $k$ large enough since $%
D(\gamma _{j_{k}}(t_{k}))$ $\rightarrow $ $0$ by the assumption: $D$ $\equiv 
$ $0$ on $\Gamma _{\infty }.$ So $D^{\prime }$ is strictly increasing at $%
\gamma _{j_{k}}(t_{k}).$ This contradicts that it achieves a minimum at $%
\gamma _{j_{k}}(t_{k})$ unless $t_{k}$ $=$ $s_{k}.$ But at $s_{k}$ $%
D^{\prime }(\gamma _{j_{k}}(s_{k}))$ $=$ $D^{\prime }(q_{k})$ $\geq $ $\frac{%
\lambda }{8},$ a contradiction. Now by the mean-value theorem and (\ref{3.13}%
), we have%
\begin{eqnarray}
D(\gamma _{j_{k}}(s_{0}))-D(\gamma _{j_{k}}(\bar{s}_{0})) &=&D^{\prime
}(\gamma _{j_{k}}(\tilde{s}_{k}))(s_{0}-\bar{s}_{0})  \label{3.14} \\
&\geq &\frac{\lambda }{8}(s_{0}-\bar{s}_{0}).  \notag
\end{eqnarray}

\noindent where $\tilde{s}_{k}$ $\in $ $[\bar{s}_{0},s_{0}].$ Letting $k$ $%
\rightarrow $ $\infty $ in the left-hand side of (\ref{3.14}), we get%
\begin{eqnarray*}
0 &=&0-0=D(\gamma _{\infty }(s_{0}))-D(\gamma _{\infty }(\bar{s}_{0})) \\
&\geq &\frac{\lambda }{8}(s_{0}-\bar{s}_{0})>0,
\end{eqnarray*}

\noindent a contradiction. We have proved that for $j$ large enough, there
holds%
\begin{equation}
M_{j}<\frac{\lambda }{4}.  \label{3.15}
\end{equation}

Similarly, we can show that for $j$ large enough, there holds%
\begin{equation}
m_{j}>-\frac{\lambda }{4}.  \label{3.16}
\end{equation}

\noindent We then write%
\begin{eqnarray}
D(\gamma _{j_{k}}(s_{0}))-D(\gamma _{j_{k}}(\bar{s}_{0})) &=&D^{\prime
}(\gamma _{j_{k}}(\bar{s}_{0}))(s_{0}-\bar{s}_{0})  \label{3.17} \\
&&+D^{\prime \prime }(\gamma _{j_{k}}(\hat{s}_{k}))\frac{(s_{0}-\bar{s}%
_{0})^{2}}{2}  \notag
\end{eqnarray}

\noindent for $\hat{s}_{k}$ $\in $ $[\bar{s}_{0},s_{0}].$ Evaluate (\ref%
{3.13'}) at $\gamma _{j_{k}}(\hat{s}_{k}).$ Observe that in the right-hand
side of (\ref{3.13'}), the numerator of the first term is bounded away from
zero by (\ref{3.15}), (\ref{3.16}) while the denominator $D$ goes to zero,
the second term is bounded, and the third term goes to zero. It follows that 
$D^{\prime \prime }(\gamma _{j_{k}}(\hat{s}_{k}))$ $\rightarrow $ $\infty $
as $k$ $\rightarrow $ $\infty .$ On the other hand, the two terms in the
left-hand side of (\ref{3.17}) tends to zero as $k$ $\rightarrow $ $\infty $
while the first term in the right-hand side of (\ref{3.17}) is bounded due
to (\ref{3.15}), (\ref{3.16}). Altogether we get $0$ $=$ $\infty $ as $k$ $%
\rightarrow $ $\infty $ in (\ref{3.17}), a contradiction. We can then
conclude that $\Gamma _{\infty }$ contains no singular points. Since $\Gamma
_{\infty }$ is a $C^{2}$ limit of $\Gamma _{j},$ it must be characteristic.
We have proved (b).

\endproof%

\bigskip

\textbf{Proposition 3.3}. \textit{Consider a }$C^{1}$\textit{\ smooth }$p$%
\textit{-minimal graph over a planar domain }$\Omega .$ \textit{Let }$p$%
\textit{\ be a singular point of }$\Omega $\textit{. Then }$p$\textit{\
emits at least one characteristic ray}$.$\textit{\ That is to say, there
exists at least one characteristic line with }$p$\textit{\ as an end point
in }$\Omega .$

\bigskip

\proof
By Lemma 3.1 we can find a sequence of nonsingular points $q_{j}$ converging
to $p.$ Let $\Gamma _{j}$ denote the characteristic line passing through $%
q_{j}.$ Observe that $\Gamma _{j}$ can only hit at most one singular point
in $\Omega $ by Theorem B (c) in Section 1. Suppose there exists a
subsequence of $\Gamma _{j}$ which do not hit singular points in a ball $%
B_{\delta }(p)$ $\subset $ $\Omega ,$ $\delta $ $>$ $0.$ Then $\Gamma
_{\infty },$ the limit of $\Gamma _{j},$ must pass through $p.$ On the other
hand, $\Gamma _{\infty }$ is a characteristic line by Lemma 3.2. This
contradicts the fact that $p$ $\in $ $\Gamma _{\infty }$ is a singular
point. So we may assume that each $\Gamma _{j}$ hits a singular point $s_{j}$
$\in $ $B_{\varepsilon }(p)$ for some $\varepsilon $ $>$ $0$ ($s_{j}$ may
coincide with $p).$ We still let $\Gamma _{\infty }$ be the limit of $\Gamma
_{j}.$ Let $s_{\infty }$ be a limit of (any subsequence) of $s_{j}.$ Since $%
q_{j}\rightarrow p$ and a characteristic line can only hit at most one
singular point (see Theorem B (c)), we must have $s_{\infty }$ $=$ $p$ and $%
\Gamma _{\infty }$ hits $p.$ So $\Gamma _{\infty }$ is the characteristic
ray that we want.

\endproof%

\bigskip

Next we will discuss the situation in which $p$ emits two different
characteristic rays $\Gamma _{1}$, $\Gamma _{2}.$ Denote the fan-shaped
region surrounded by $\Gamma _{1}$, $p,$ and $\Gamma _{2}$ by \{$\Gamma
_{1}p\Gamma _{2}\}$ (see Figure 3.2 below).

\begin{figure}[th]
\begin{center}
\includegraphics[width=9cm]{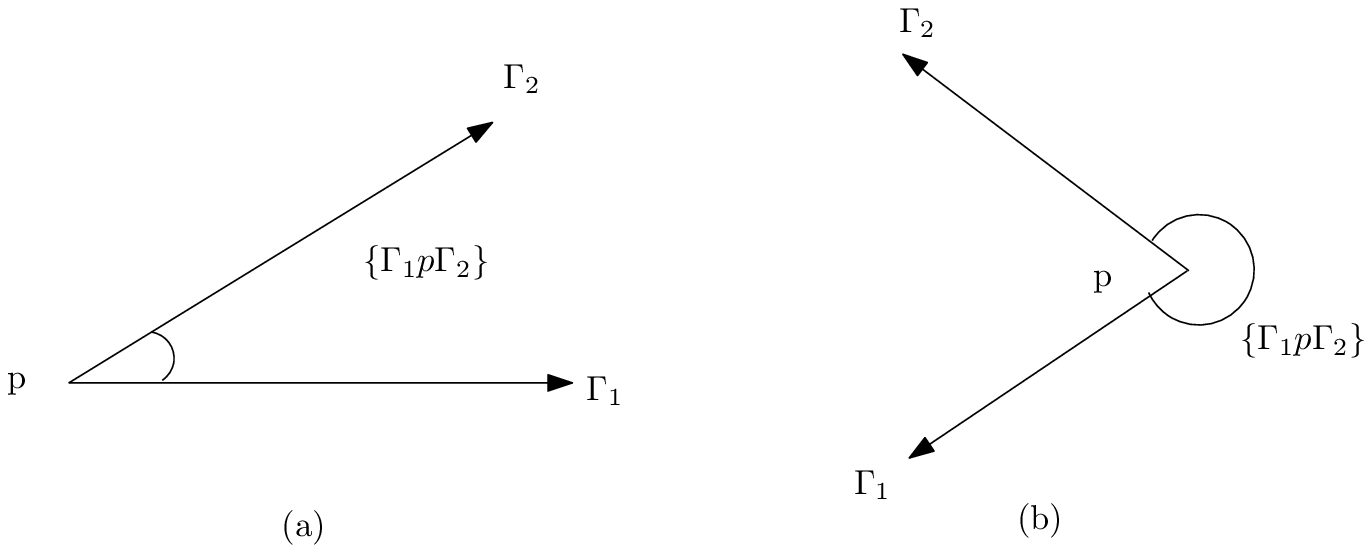}\\[0pt]
Figure 3.2
\end{center}
\par
\end{figure}

\textbf{Proposition 3.4}. \textit{Consider a }$C^{1}$\textit{\ smooth }$p$%
\textit{-minimal graph over a planar domain }$\Omega .$ \textit{Let }$p$%
\textit{\ be a singular point of }$\Omega $\textit{. Suppose that }$p$%
\textit{\ emits two different characteristic rays }$\Gamma _{1}$\textit{, }$%
\Gamma _{2}.$\textit{\ Suppose \{}$\Gamma _{1}p\Gamma _{2}\}$\textit{\ }$%
\cap $\textit{\ }$\Omega $\textit{\ contains no singular points. Then }$p$%
\textit{\ emits characteristic rays (which stop when hitting }$\partial
\Omega )$\textit{\ at all the directions pointing inside \{}$\Gamma
_{1}p\Gamma _{2}\}$\textit{\ (see Figure 3.3)}$.$

\bigskip

\proof
\textbf{Case 1}. Suppose that the angle between $\Gamma _{1}$ and $\Gamma
_{2}$ in the region \{$\Gamma _{1}p\Gamma _{2}\}$ is less than $\pi $ (see
Figure 3.3 (a))$.$ Then there exists $\varepsilon $ $>$ $0$ such that the
tangent line $L_{q}$ at each point $q$ $\in $ $\partial B_{\varepsilon }(p)$ 
$\cap $ \{$\Gamma _{1}p\Gamma _{2}\}$ hits either $\Gamma _{1}$ or $\Gamma
_{2}$ and lies in $\Omega $ before hitting one of them. Now the
characteristic line $\Gamma _{q}$ passing through $q$ cannot be $L_{q}$
since two characteristic lines do not intersect at a nonsingular points (see
Theorem $B^{\prime }$ in \cite{CHY2}). So $\Gamma _{q}$ must hit ($\Gamma
_{1}$ $\cup $ $\{p\}$ $\cup $ $\Gamma _{2})$ $\cap $ $B_{\varepsilon }(p).$
Since $\Gamma _{q}$ cannot intersect with either $\Gamma _{1}$ or $\Gamma
_{2},$ $\Gamma _{q}$ has to hit $p.$ We are done.

\begin{figure}[ht]
\begin{center}
\includegraphics[width=9cm]{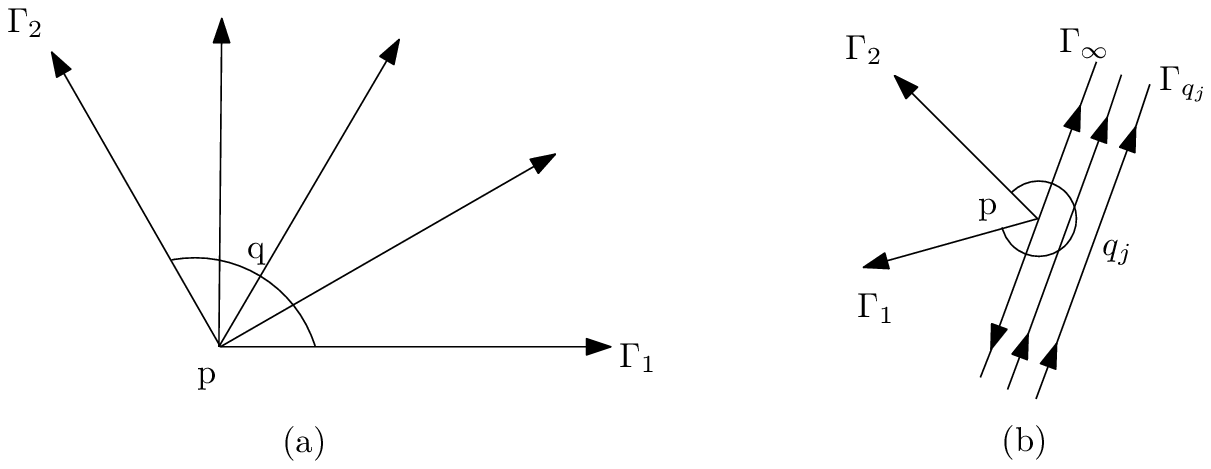}\\[0pt]
Figure 3.3
\end{center}
\par
\end{figure}

\textbf{Case 2.} Suppose that the angle between $\Gamma _{1}$ and $\Gamma
_{2}$ in the region \{$\Gamma _{1}p\Gamma _{2}\}$ is larger than $\pi .$
Assume that the conclusion fails. Then there exists a sequence of
characteristic lines $\Gamma _{q_{j}}$ through $q_{j}$ as $q_{j}\rightarrow
p,$ which do not hit $p$ (see Figure 3.3 (b))$.$ Observe that \{$\Gamma
_{1}p\Gamma _{2}\}$ $\cap $ $\Omega $ contains no singular points and $%
\Gamma _{q_{j}}$'$s$ do not intersect, say, in a ball $B_{r}(p)$ $\subset $ $%
\Omega .$ It follows that a subsequence of $\Gamma _{q_{j}}$ converges to a
straight line $\Gamma _{\infty }$ passing through $p.$ According to Lemma
3.2, $\Gamma _{\infty }$ is a characteristic line (so containing no singular
point). We have reached a contradiction.

\bigskip 
\endproof%

\textbf{Lemma 3.5}. \textit{Consider a }$C^{1}$\textit{\ smooth }$p$\textit{%
-minimal graph over a planar domain }$\Omega .$ \textit{Let }$p$\textit{\
and }$q$\textit{\ be two singular points in }$\Omega $\textit{\ such that }$%
\bar{B}_{d(p,q)}(p)$\textit{\ }$\subset $\textit{\ }$\Omega $\textit{\ where 
}$d(p,q)$\textit{\ denotes the Euclidean distance between }$p$\textit{\ and }%
$q.$\textit{\ Then there exists a }$C^{0}$\textit{\ singular curve }$\gamma $%
\textit{\ }$:$\textit{\ }$[0,1]$\textit{\ }$\rightarrow $\textit{\ }$\bar{B}%
_{d(p,q)}(p)$\textit{\ such that }$\gamma (0)$\textit{\ }$=$\textit{\ }$p$%
\textit{\ and }$\gamma (1)$\textit{\ }$\in $\textit{\ }$\partial
B_{d(p,q)}(p)$ (\textit{note that }$\gamma (1)$\textit{\ may not be }$q).$

\bigskip

\proof
First by Proposition 3.3 we have a characteristic ray $\Gamma _{p}$ emitted
from $p.$ We orient the circle $\mathfrak{S}$ $:=$ $\partial B_{d(p,q)}(p)$
counterclockwise and view the point $q_{0}$ $=$ $\Gamma _{p}$ $\cap $ $%
\mathfrak{S}$ as the starting point. Then there exists a point $\bar{q}$ $%
\in $ $\mathfrak{S}$ (may be $q_{0})$ between $q_{0}$ and $q$ such that for
any point $q^{\prime }$ on the closed arc $q_{0}\bar{q},$ the characteristic
line passing through $q^{\prime }$ meets $p$ and for $\zeta $ $\in $ $%
\mathfrak{S}$ beyond $\bar{q},$ but near $\bar{q},$ the characteristic line $%
\Gamma _{\zeta }$ passing through $\zeta $ does not meet $p$ (see Figure 3.4
below)$.$ Note that since $\bar{q}$ is nonsingular (otherwise, we contradict
Lemma 3.2), we may assume that all nearby $\zeta $ are nonsingular.

\begin{figure}[ht]
\begin{center}
\includegraphics[width=5.5cm]{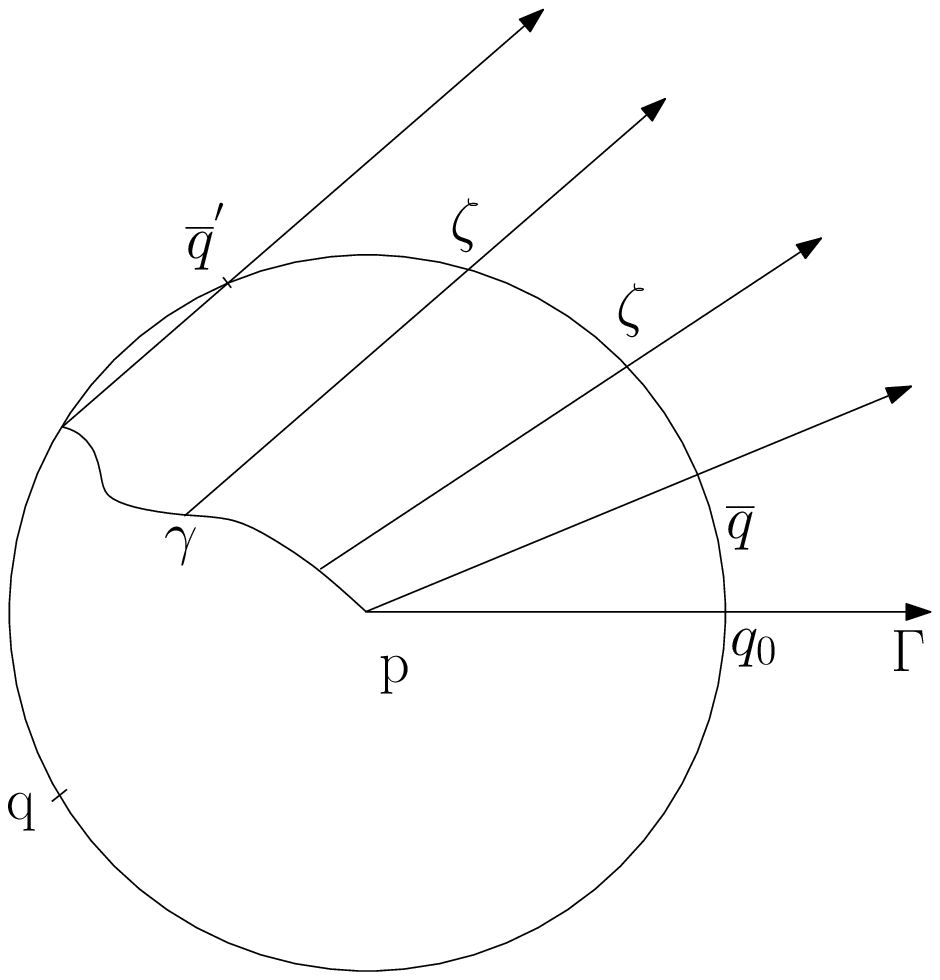}\\[0pt]
Figure 3.4
\end{center}
\par
\end{figure}

We claim that $\Gamma _{\zeta }$ has to meet a singular point other than $p$
in $B_{d(p,q)}(p)$ for any $\zeta $ in the open arc $\bar{q}\bar{q}^{\prime
} $ where $\bar{q}^{\prime }$ $\in $ $\mathfrak{S}\backslash \{$closed arc $%
q_{0}\bar{q}\}$ is close enough to $\bar{q}.$ Suppose the converse holds.
Then there exists a sequence \{$\zeta _{j}\}$ such that $\zeta _{j}$ $%
\rightarrow $ $\bar{q}$ and $\Gamma _{\zeta _{j}}$ does not meet any
singular point in $B_{d(p,q)}(p)$ for any $\zeta _{j}.$ On the other hand,
the sequence of lines $\Gamma _{\zeta _{j}}$ $\cap $ $B_{d(p,q)}(p)$
converges to a straight line containing $\Gamma _{\bar{q}}$ as $\zeta _{j}$ $%
\rightarrow $ $\bar{q}.$ Since $\Gamma _{\bar{q}}$ contains the singular
point $p,$ we have reached a contradiction to Lemma 3.2. Now let $\mathfrak{s%
}(\zeta )$ denote the first singular point that $\Gamma _{\zeta }$ meets in $%
B_{d(p,q)}(p)$. We claim that $\mathfrak{s}$ is continuous on the open arc $%
\bar{q}\bar{q}^{\prime }$ and can be continuously extended to $\bar{q}$ and $%
\bar{q}^{\prime }$ so that $\mathfrak{s}(\bar{q})$ $=$ $p.$ Suppose the
converse holds. Then there exists a sequence \{$\zeta _{j}\}$ $\rightarrow $ 
$\bar{\zeta}$ in the closed arc $\bar{q}\bar{q}^{\prime }$ such that $\lim 
\mathfrak{s}(\zeta _{j})$ $\in $ $\bar{B}_{d(p,q)}(p)$ exists and is
different from $\mathfrak{s}(\bar{\zeta}).$ Let $\mathfrak{\bar{s}}$ $=$ $%
\lim \mathfrak{s}(\zeta _{j}).$ Since the set of singular points is closed, $%
\mathfrak{\bar{s}}$ is a singular point. Observe that \{$\Gamma _{\zeta
_{j}}\}$ converges to a line segment $\bar{\Gamma}$ ending at $\mathfrak{%
\bar{s}}.$ By the continuity of $N^{\perp }$ at $\bar{\zeta},$ $\bar{\Gamma}$
contains $\Gamma _{\bar{\zeta}}$ since $\mathfrak{s}(\bar{\zeta})$ is the
first singular point that $\Gamma _{\bar{\zeta}}$ meets. We are in a
situation in which Lemma 3.2 applies. But $\mathfrak{\bar{s}}$ $\neq $ $%
\mathfrak{s}(\bar{\zeta})$ implies that the limit $\bar{\Gamma}$ contains $%
\mathfrak{s}(\bar{\zeta}),$ a singular point. This contradicts the
conclusion of Lemma 3.2. Thus the map $\mathfrak{s}$ with $\mathfrak{s}(\bar{%
q})$ $=$ $p$ is continuously defined and the domain closed arc $\bar{q}\bar{q%
}^{\prime }$ can be extended so that either $\mathfrak{s}(\bar{q}^{\prime })$
$\in $ $\mathfrak{S}$ or $\bar{q}^{\prime }$ becomes a singular point
(before hitting $q$ or $=$ $q)$. In the latter case, we claim that $%
\lim_{\zeta \rightarrow \bar{q}^{\prime }}\mathfrak{s}(\zeta )$ $=$ $\bar{q}%
^{\prime }$ ($\in $ $\mathfrak{S}).$ Otherwise, there is a sequence \{$\zeta
_{j}^{\prime }\}$ $\rightarrow $ $\bar{q}^{\prime }$ such that $\lim_{\zeta
_{j}^{\prime }\rightarrow \bar{q}^{\prime }}\mathfrak{s}(\zeta _{j}^{\prime
})$ exists in $\bar{B}_{d(p,q)}(p),$ but is different from $\bar{q}^{\prime
}.$ Now the limit line seqment $\Gamma _{\infty }$ of $\Gamma _{\zeta
_{j}^{\prime }}$ is characteristic by Lemma 3.2 while $\Gamma _{\infty }$
meets two distinct singular points, which contradicts Theorem B (c). We have
shown that $\mathfrak{s}$ $:$ closed arc $\bar{q}\bar{q}^{\prime }$ $%
\rightarrow $ $\bar{B}_{d(p,q)}(p)$ is continuous with $\mathfrak{s}(\bar{q}%
) $ $=$ $p$ and $\mathfrak{s}(\bar{q}^{\prime })$ $\in $ $\mathfrak{S}$ ($:=$
$\partial B_{d(p,q)}(p)).$ Parametrize the closed arc $\bar{q}\bar{q}%
^{\prime }$ continuously by the map $l$ $:$ $[0,1]$ $\rightarrow $ closed
arc $\bar{q}\bar{q}^{\prime }$ so that $l(0)$ $=$ $\bar{q},$ $l(1)$ $=$ $%
\bar{q}^{\prime }.$ Let $\gamma $ $=$ $\mathfrak{s}\circ l.$ It is now clear
that $\gamma $ satisfies the required property.

\endproof%

\bigskip

By a similar argument replacing characteristic lines by characteristic
curves, we can extend Proposition 3.3 and Lemma 3.5 to a general situation,
based on Lemma 3.2$^{\prime }.$

\bigskip

\textbf{Proposition 3.3}$^{\prime }.$ \textit{Let }$\Omega $\textit{\ be a
bounded domain of }$R^{2}.$ \textit{Let }$u$\textit{\ }$\in $\textit{\ }$%
C^{1}(\Omega )$\textit{\ be a weak solution to (}$\ref{1.1})$\textit{\ with }%
$\vec{F}$\textit{\ }$\in $\textit{\ }$C^{1}(\Omega )$\textit{\ and }$H$%
\textit{\ }$\in $\textit{\ }$C^{0}(\Omega ).$\textit{\ Assume further }$%
N^{\perp }(curl\vec{F})$\textit{\ and }$N(H)$\textit{\ exist and are
continuous (extended over singular points) in }$\Omega .$\textit{\ Suppose }$%
\func{curl}\vec{F}$ $\neq $ $0.$ \textit{Let }$p$\textit{\ be a singular
point of }$\Omega $\textit{. Then there exists at least one characteristic
curve with }$p$\textit{\ as an end point in }$\Omega .$

\bigskip

\textbf{Lemma 3.5}$^{\prime }.$ \textit{Let }$\Omega $\textit{\ be a bounded
domain of }$R^{2}.$ \textit{Let }$u$\textit{\ }$\in $\textit{\ }$%
C^{1}(\Omega )$\textit{\ be a weak solution to (}$\ref{1.1})$\textit{\ with }%
$\vec{F}$\textit{\ }$\in $\textit{\ }$C^{1}(\Omega )$\textit{\ and }$H$%
\textit{\ }$\in $\textit{\ }$C^{0}(\Omega ).$\textit{\ Assume further }$%
N^{\perp }(curl\vec{F})$\textit{\ and }$N(H)$\textit{\ exist and are
continuous (extended over singular points) in }$\Omega .$\textit{\ Suppose }$%
\func{curl}\vec{F}$ $\neq $ $0.$ \textit{Let }$p$\textit{\ and }$q$\textit{\
be two singular points in }$\Omega $\textit{\ such that }$\bar{B}%
_{d(p,q)}(p) $\textit{\ }$\subset $\textit{\ }$\Omega $\textit{\ where }$%
d(p,q)$\textit{\ denotes the Euclidean distance between }$p$\textit{\ and }$%
q.$\textit{\ Assume }$|H|$ $\leq $ $\frac{1}{2d(p,q)}$ in $B_{d(p,q)}(p).$ 
\textit{Then there exists a }$C^{0}$\textit{\ singular curve }$\gamma $%
\textit{\ }$:$\textit{\ }$[0,1]$\textit{\ }$\rightarrow $\textit{\ }$\bar{B}%
_{d(p,q)}(p)$\textit{\ such that }$\gamma (0)$\textit{\ }$=$\textit{\ }$p$%
\textit{\ and }$\gamma (1)$\textit{\ }$\in $\textit{\ }$\partial
B_{d(p,q)}(p)$ (\textit{note that }$\gamma (1)$\textit{\ may not be }$q).$

\bigskip

Note that the condition on the bound of $H$ is to guarantee that
characteristic curves behave like straight lines in a ball.

\bigskip

\proof
\textbf{(of Theorem C)} Suppose that $p$ is not an isolated singular point.
Then near $p$ we can find another singular point $q$ such that $\bar{B}%
_{d(p,q)}(p)$\textit{\ }$\subset $\textit{\ }$\Omega $ and $|H|$ $\leq $ $%
\frac{1}{2d(p,q)}$ in $B_{d(p,q)}(p).$\textit{\ }Now it is clear that the $%
C^{0}$ singular curve obtained in Lemma 3.5$^{\prime }$ satisfies the
desired property. We have proved $(a)$.

Let $A$ $\subset $ $S(u)$ denote the path-connected component containing $p.$
To prove $(b),$ suppose the converse holds. Then there exist a sequence $%
\varepsilon _{j}$ $\rightarrow $ $0$ and a sequence $p_{j}$ $\in $ [$S(u)$ $%
\cap $ $B_{\varepsilon _{j}}(p)]\backslash A.$ By Proposition 3.3$^{\prime }$
we can find a characteristic curve $\Gamma _{p_{j}}$ emitted from $p_{j}.$
Take a ball $B_{\delta }(p)$ of radius $\delta $ such that $|H|$ $\leq $ $%
\frac{1}{2\delta }$ in $B_{\delta }(p).$ For $j$ large, $\Gamma _{p_{j}}$
must go out of the ball $B_{\delta }(p)$ and hit a nonsingular point $\tilde{%
p}_{j}$ on $\partial B_{\delta }(p)$ by Theorem B (c)$.$ There is a
subsequence, still denoted $\tilde{p}_{j},$ converging to $\tilde{p}_{\infty
}$ $\in $ $\partial B_{\delta }(p).$ By a similar argument as in the proof
of Lemma 3.2$^{\prime },$ we can show that $\Gamma _{p_{j}}$ converges to a
characteristic curve $\Gamma _{\infty }$ $\subset $ $B_{\delta }(p),$
emitted from $\tilde{p}_{\infty },$ in $C^{1}$ on a compact parameter
interval. Since $p_{j}$ $\rightarrow $ $p,$ $\Gamma _{\infty }$ cannot hit
another singular point before hitting $p$ by Lemma 3.2$^{\prime }$ and the
uniqueness of characteristic curves through a point (see Theorem B$^{\prime
} $ of \cite{CHY2}). So $\Gamma _{\infty }$ connects $p$ with $\tilde{p}%
_{\infty }.$ From Theorem B (c) we conclude that $\tilde{p}_{\infty }$ is a
nonsingular point since $p$ is singular. Let $\mathcal{V}$ $\subset $ $%
\partial B_{\delta }(p)\backslash S(u)$ be a neighborhood of $\tilde{p}%
_{\infty }$. For $\tilde{p}$ $\in $ $\mathcal{V}$, we define $\mathfrak{s(}%
\tilde{p})$ to be the (first) singular point which the characteristic curve
through $\tilde{p}$ hits in $B_{\delta }(p).$ Now in view of Lemma 3.2$%
^{\prime }$ the map $\mathfrak{s:}$ $\tilde{p}$ $\in $ $\mathcal{V}^{\prime
} $ $\rightarrow $ $\mathfrak{s(}\tilde{p})$ $\in S(u)$ $\cap $ $B_{\delta
}(p) $ is defined and continuous for a perhaps smaller neighborhood $%
\mathcal{V}^{\prime }$ $\subset $ $\mathcal{V}$. Thus $p$ $=$ $\mathfrak{s(}%
\tilde{p}_{\infty })$ and $p_{j}$ $=$ $\mathfrak{s(}\tilde{p}_{j})$ are
path-connected in $\mathfrak{s}(\mathcal{V}^{\prime }).$ So $p_{j}$ $\in $ $%
A,$ a contradiction.

\endproof%

\bigskip

\textbf{Proposition 3.6. }\textit{Let }$\Omega $\textit{\ be a bounded
domain of }$R^{2}.$ \textit{Let }$u$\textit{\ }$\in $\textit{\ }$%
C^{1}(\Omega )$\textit{\ be a weak solution to (}$\ref{1.1})$\textit{\ with }%
$\vec{F}$\textit{\ }$\in $\textit{\ }$C^{1}(\Omega )$\textit{\ and }$H$%
\textit{\ }$\in $\textit{\ }$C^{0}(\Omega ).$\textit{\ Assume further }$%
N^{\perp }(curl\vec{F})$\textit{\ and }$N(H)$\textit{\ exist and are
continuous (extended over singular points) in }$\Omega .$\textit{\ Suppose }$%
\func{curl}\vec{F}$ $\neq $ $0.$ \textit{Let }$p_{0}$\textit{\ be a singular
point in }$\Omega .$\textit{\ Take }$B_{r_{1}}(p_{0})$ $\subset \subset $ $%
\Omega $ \textit{such that }$0$ $<$ $r_{1}$\textit{\ }$\leq $ $%
(\sup_{B_{r_{1}}(p_{0})}|H|)^{-1}.$\textit{\ Then there exists }$r_{2},$%
\textit{\ }$0$\textit{\ }$<$\textit{\ }$r_{2}$\textit{\ }$<$\textit{\ }$%
r_{1},$\textit{\ such that for any nonsingular point }$p$\textit{\ }$\in $%
\textit{\ }$B_{r_{2}}(p_{0})\backslash \{p_{0}\},$\textit{\ the
characteristic curve }$\Gamma _{p}$\textit{\ passing through }$p$\textit{\
must hit a singular point in }$B_{r_{1}}(p_{0}).$

\bigskip

\proof
Suppose the converse is true. Then there exists a sequence of points $p_{j},$
converging to $p_{0},$ such that each characteristic curve $\Gamma _{p_{j}}$
does not hit any singular point in $B_{r_{1}}(p_{0}).$ By Lemma 3.2$^{\prime
}$ (a) we conclude that a subsequence of closed arcs $\subset $ $\Gamma
_{p_{j}}$ converges to a closed arc $\Gamma _{\infty }$ while $\Gamma
_{\infty }$ is a characteristic curve containing $p_{0}$ by Lemma 3.2$%
^{\prime }$ (b)$.$ We have reached a contradiction.

\endproof%

\bigskip

\proof
\textbf{(of Theorem D) }Take\textbf{\ }$r_{1}$ $>$ $0$ such that on $\bar{B}%
_{r_{1}}(p),$ $|H|$ $\leq $ $C_{1}$ $<<$ $\frac{1}{r_{1}}$ for some positive
constant $C_{1}.$ In $B_{r_{1}}(p)$ any characteristic curve has to hit two
points on the boundary of $B_{r_{1}}(p)$ if it does not meet $p.$ We claim
that there exists $0$ $<$ $r_{0}$ $<$ $r_{1}$ such that for any $q$\ $\in $\ 
$\bar{B}_{r_{0}}(p)\backslash \{p\}$, the characteristic curve $\tilde{\Gamma%
}_{q}$ passing through $q$\ has to meet $p.$ Suppose the converse holds.
Then we can find a sequence of points $q_{j}$ approaching $p$ such that all
the $\tilde{\Gamma}_{q_{j}}$'s do not meet $p$ and satisfy the condition in
Lemma $3.2^{\prime }.$ By Lemma $3.2^{\prime }$ (a) we can find a
subsequence of closed arcs $\Gamma _{q_{j}}$ $\subset $ $\tilde{\Gamma}%
_{q_{j}}$ $\cap $ $B_{r_{2}}(p)$ for $0$ $<$ $r_{2}$ $<$ $r_{1}$, which
converges to $\Gamma _{\infty }.$ According to Lemma $3.2^{\prime }$ (b)$,$ $%
\Gamma _{\infty }$ contains $p,$ a singular point$.$ This contradicts
another statement of Lemma $3.2^{\prime }$ (b) that $\Gamma _{\infty }$
contains no singular points. We have proved the first part of the theorem,
that is, $\Gamma _{q}$ meets $p.$

Next we observe that $|\theta ^{\prime }(x)|$ $=$ $|-H|$ $\leq $ $C_{1}$ for 
$x$ $\in $ $\Gamma _{q}$ $\cap $ $\bar{B}_{r_{1}}(p),$ and hence $N^{\perp
}(x)$ $=$ $(\sin \theta (x),$ $-\cos \theta (x))$ is Cauchy in $x$ near $p.$
Therefore the unit tangent vector $N^{\perp }$\ of $\Gamma _{q}$\ has a
limit at $p,$ denoted $v(q).$\textit{\ }Define the map $\psi $\ $:$\ $q$\ $%
\in $\ $\partial B_{r_{0}}(p)$\ $\rightarrow $\ $v(q)$\ $\in $\ $T_{p}\Omega
.$\ We can extend by a similar proof Theorem $B^{\prime }$ in \cite{CHY2}
for the uniqueness of characteristic curves to include the case that $p$ is
a singular point by a similar proof. We can then conclude that $\psi $ is
injective. Note that $q$ $\in $ $\partial B_{r_{0}}(p)$ has two "sides"
locally in $\partial B_{r_{0}}(p)$. When $q^{\prime }$ $\rightarrow $ $q$
clockwise (counterclockwise, resp.), we write $q^{\prime }$ $\rightarrow $ $%
q^{+}$ ($q^{\prime }$ $\rightarrow $ $q^{-},$ resp.). Observe that $%
\lim_{q^{\prime }\rightarrow q^{+}}$ $v(q^{\prime })$ exists since
characteristic curves do not intersect in $B_{r_{0}}(p)\backslash \{p\}$
(and hence \{$v(q^{\prime })\}$ is "ordered"). Let $w$ :$=$ $\lim_{q^{\prime
}\rightarrow q^{+}}$ $v(q^{\prime }).$ Suppose $w$ $\neq $ $v(q).$ From the
standard O.D.E. theory (continuous dependence of the solution on initial
data), the solution curve of (\ref{1.7}) with initial tangent $w$ has to
meet $q.$ We now have two distinct characteristic curves passing through $q,$
a contradiction, so $w$ $=$ $v(q).$ Similarly we also have $\lim_{q^{\prime
}\rightarrow q^{-}}$ $v(q^{\prime })$ $=$ $v(q).$ We have shown that $\psi $
is $C^{0}.$

Take two different points $q_{1},$ $q_{2}$ $\in $ $\partial B_{r_{0}}(p).$
Consider the image $\psi ([q_{1}q_{2}])$ of $\psi $ from the small (large,
resp.) arc $[q_{1}q_{2}]$ formed by $q_{1},$ $q_{2}$ into the small or large
arc [$v(q_{1})v(q_{2})]$ formed by $v(q_{1}),$ $v(q_{2}).$ Since $\psi $ is
continuous, the image of a path-connected set is path-connected, so $\psi
([q_{1}q_{2}])$ is path-connected. On the other hand, $\psi ([q_{1}q_{2}])$
contains $v(q_{1})$ and $v(q_{2}),$ and hence contains [$v(q_{1})v(q_{2})].$
It follows that $\psi ([q_{1}q_{2}])$ $=$ $[v(q_{1})v(q_{2})]$ (similar
formula holds for another arc). We have shown that $\psi $ is surjective
onto the space of unit tangent vectors at $p$. The continuity of $\psi ^{-1}$
follows from the standard O.D.E. theory (continuous dependence of the
solution on initial data). We have completed the proof.

\endproof%

\bigskip

According to Theorem D and Theorem B (c), it is impossible to have a $p$%
-minimal graph over a convex planar domain with two isolated singular
points. However, if the domain is not convex, this is possible as shown by
the following configuration of characteristic lines and singular points (see
Figure 3.5 below: $S_{1},S_{2}$ are two isolated singular points and the
straight lines denote the characteristic lines).

\begin{figure}[ht]
\begin{center}
\includegraphics[width=4cm]{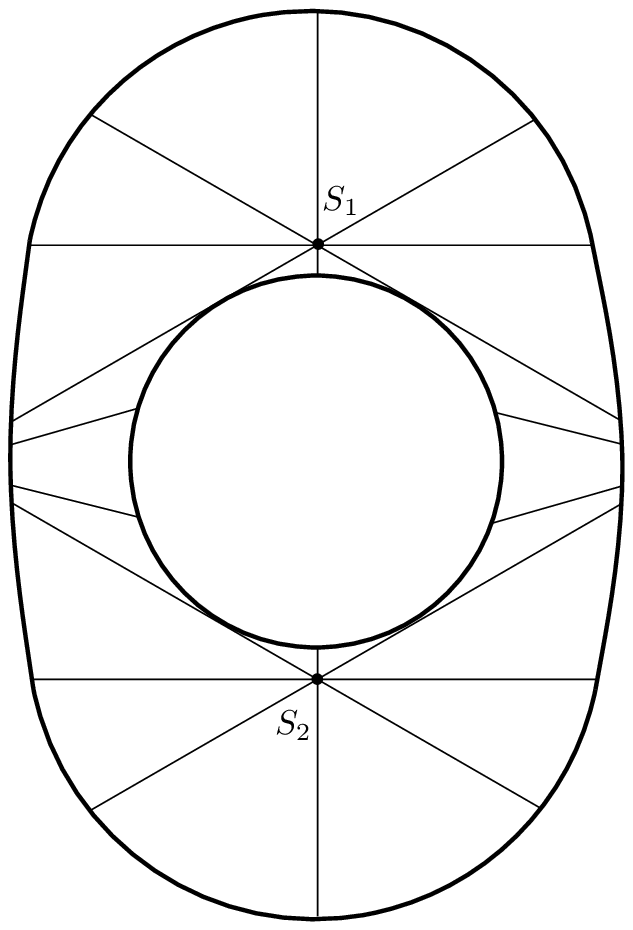}\\[0pt]
Figure 3.5
\end{center}
\par
\end{figure}

\bigskip

\section{Examples}

For a $C^{2}$ smooth $p$-minimal surface, we showed (\cite{CHMY04}) that,
among others, a singular curve (which must be $C^{1}$) never has an end
(boundary) point. But this situation can occur for a $C^{1}$ (hence not $%
C^{2})$ smooth $p$-minimal surface as shown in the following examples.

\bigskip

\textbf{Example 4.1. }According to Proposition 3.3, any singular point emits
at least one characteristic ray. Can we have an example in which a singular
point emits exactly one characteristic ray? We\textbf{\ }are going to
construct such an example. First we want the union of the negative $x$-axis
and the origin to be the singular set. Each singular point on the negative $%
x $-axis emits two characteristic rays having equal angles with the positive 
$x $-axis and the origin emits only one characteristic ray, the positive $x$%
-axis (see Figure 4.1).

\begin{figure}[ht]
\begin{center}
\includegraphics[width=7cm]{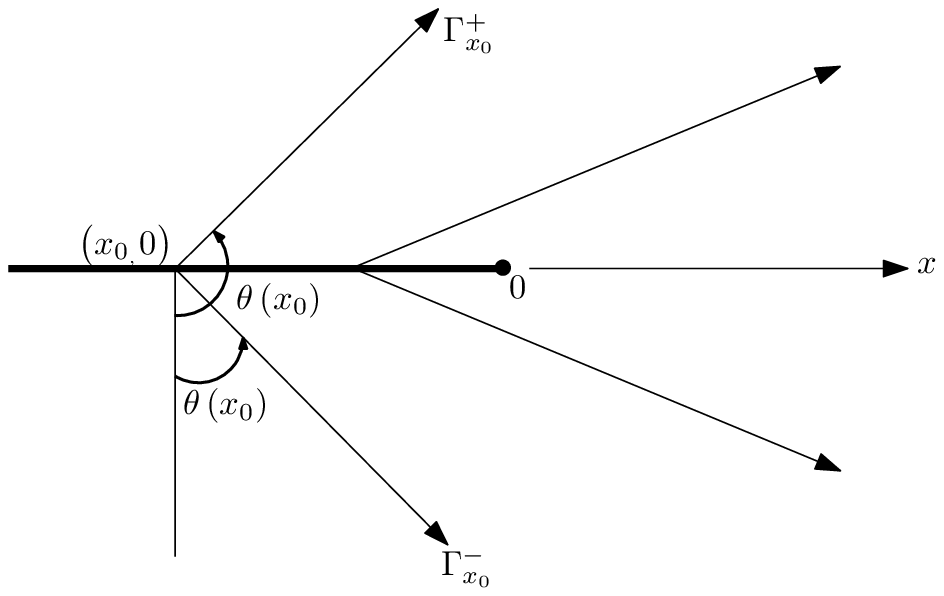}\\[0pt]
Figure 4.1
\end{center}
\par
\end{figure}

\bigskip Let $\Gamma _{x_{0}}^{+}$ ($\Gamma _{x_{0}}^{-},$ resp.) denote the
characteristic ray emitted from $(x_{0},0),$ $x_{0}$ $<$ $0,$ on the upper
(lower, resp.) half plane. Describe a point $(x,y)$ $\in $ $\Gamma
_{x_{0}}^{+}$ $\cup $ $\{(x_{0},0)\}$ as follows:%
\begin{eqnarray}
x &=&x_{0}+s\sin \theta (x_{0})  \label{3.6.1} \\
y &=&-s\cos \theta (x_{0})  \notag
\end{eqnarray}

\noindent where $s$ $\geq $ $0$ and $\theta (x_{0})$ $:=$ $\frac{\pi }{2}$ $%
+ $ ($-,$ resp.) "\textit{the angle between }$\Gamma _{x_{0}}^{+}$\textit{\ (%
}$\Gamma _{x_{0}}^{-},$ resp.) \textit{and} \textit{the positive }$x$\textit{%
-axis" }(we choose $\theta $ such that the characteristic direction is $%
N^{\perp }$ $=$ $(\sin \theta ,-\cos \theta )$ accordingly to the previous
choice). Integrating $du$ $+$ $xdy$ $-$ $ydx$ $=$ $0$ along the ray
described by (\ref{3.6.1}), we determine the $u$-value uniquely once we know 
$u(x_{0},0).$ Taking $u(x_{0},0)$ $=$ $0,$ we then obtain%
\begin{equation}
u(x,y)=-x_{0}y.  \label{3.6.2}
\end{equation}

\noindent (For characteristic rays on the lower half plane, we use a similar
argument to write $u$ as in (\ref{3.6.2})). We require $\theta $ to be $%
C^{1} $ in $x_{0}$ $\in $ $(-\infty ,0]$ such that $\frac{\pi }{2}$ $<$ $%
\theta (x_{0})$ $<$ $\frac{3\pi }{2}$ for $x_{0}$ $\in $ $(-\infty ,0),$ $%
\theta ^{\prime }(x_{0})$ $\leq $ $0,$ $\theta ^{\prime }(0)$ $=$ $a$ $<$ $%
0, $ and $\lim_{x_{0}\rightarrow 0^{-}}\theta (x_{0})$ $=$ $\theta (0)$ $:=$ 
$\frac{\pi }{2}.$ It follows from (\ref{3.6.1}) that%
\begin{equation}
\det \frac{\partial (x,y)}{\partial (x_{0},s)}=-s\theta ^{\prime
}(x_{0})-\cos \theta (x_{0})>0  \label{3.6.3}
\end{equation}

\noindent unless $x_{0}$ $=$ 0, where $\theta (0)$ $=$ $\frac{\pi }{2}.$
Observe that the $C^{1}$ smooth map $\Psi _{+}$ from $(-\infty ,0)$ $\times $
$R^{+}$ to the upper half plane, defined by $(x_{0},s)$ $\rightarrow $ $%
(x,y) $ according to (\ref{3.6.1}), is globally one to one and onto. So $%
\Psi _{+}^{-1}$ exists and is $C^{1}$ smooth by the inverse function theorem
due to (\ref{3.6.3}). Similarly for the case of \textit{the lower half
plane, }we take $-\frac{\pi }{2}$ $<$ $\theta (x_{0})$ $<$ $\frac{\pi }{2}$
for $x_{0}$ $\in $ $(-\infty ,0)$, $\theta ^{\prime }(x_{0})$ $\geq $ $0,$ $%
\theta ^{\prime }(0)$ ($=$ $-a)$ $>$ $0,$ and $\lim_{x_{0}\rightarrow
0^{-}}\theta (x_{0})$ $=$ $\theta (0)$ $:=$ $\frac{\pi }{2}.$ The $C^{1}$
smooth map $\Psi _{-}$ $:$ $(x_{0},s)$ $\in $ $(-\infty ,0)$ $\times $ $%
R^{+} $ $\rightarrow $ $(x,y)$ $\in $ \textit{the lower half plane} defined
still by (\ref{3.6.1}) (but with different range of $\theta )$ is globally
one to one and onto. Moreover, the inverse $\Psi _{-}^{-1}$ is $C^{1}$
smooth by the inverse function theorem again since we have%
\begin{equation*}
\det \frac{\partial (x,y)}{\partial (x_{0},s)}=-s\theta ^{\prime
}(x_{0})-\cos \theta (x_{0})<0
\end{equation*}

\noindent for the new range of $\theta $ and the nonnegativity of $\theta
^{\prime }.$ Thus we can view $x_{0}$ as a $C^{1}$ smooth function of $x$
and $y$ on the whole plane except the $x$-axis. Now we claim that the
function $u$ defined by%
\begin{eqnarray}
u(x,y) &=&-x_{0}y\text{ \ for }y\neq 0  \label{3.6.4} \\
u(x,y) &=&0\text{ \ for }y=0  \notag
\end{eqnarray}

\noindent is $C^{1}$ smooth on the whole plane. Moreover, $\Gamma
_{x_{0}}^{+},$ $\Gamma _{x_{0}}^{-},$ and the positive $x$-axis are
characteristic rays of $u$ while $(-\infty ,0]$ $\times $ $\{0\}$ is the
singular set of $u.$ For the $C^{1}$ smoothness of $u,$ we only need to
check whether $u_{x}$ and $u_{y}$ are continuous on the $x$-axis. It is
clear that $u_{x}(x,0)$ $=$ $0.$ Compute $u_{y}(x,0)$ from the definition of 
$u$ as follows:%
\begin{eqnarray}
u_{y}(x,0) &=&\lim_{y\rightarrow 0}\frac{u(x,y)-u(x,0)}{y}  \label{3.6.5} \\
&=&\lim_{y\rightarrow 0}\frac{-x_{0}y-0}{y}=\lim_{y\rightarrow
0}(-x_{0}(x,y))  \notag \\
&=&0\text{ \ if }x\geq 0;\text{ }=-x\text{ \ if }x<0.  \notag
\end{eqnarray}

\noindent On the other hand, we learn from (\ref{3.6.4}) and (\ref{3.6.1})
that for $y$ $\neq $ $0$%
\begin{eqnarray}
u_{x}(x,y) &=&-\frac{\partial x_{0}}{\partial x}y  \label{3.6.6} \\
&=&\frac{\cos \theta (x_{0})}{s\theta ^{\prime }(x_{0})+\cos \theta (x_{0})}%
s\cos \theta (x_{0}).  \notag
\end{eqnarray}

\noindent Observe that 
\begin{equation}
0<\frac{\cos \theta (x_{0})}{s\theta ^{\prime }(x_{0})+\cos \theta (x_{0})}%
\leq 1  \label{3.6.7}
\end{equation}

\noindent for $s$ $\geq $ $0,$ $x_{0}$ $<$ $0.$ So from (\ref{3.6.6}) and (%
\ref{3.6.7}) we can estimate%
\begin{equation*}
|u_{x}(x,y)|\leq s|\cos \theta (x_{0})|\rightarrow 0
\end{equation*}

\noindent as $(x,y)$ $\rightarrow $ $(\bar{x},0)$ since $s$ $\rightarrow $ $%
0 $ for $\bar{x}$ $\leq $ $0$ and $\cos \theta (x_{0})$ $\rightarrow $ $\cos
\theta (0)$ $=$ $\cos \frac{\pi }{2}$ $=$ $0$ while $s$ $\rightarrow $ $\bar{%
x}$ for $\bar{x}$ $>$ $0.$ It follows that $u_{x}(x,y)$ $\rightarrow $ $%
u_{x}(\bar{x},0)$ $=$ $0$ as $(x,y)$ $\rightarrow $ $(\bar{x},0).$ That is
to say, $u_{x}$ is continuous at the points of the $x$-axis. Next for $y$ $%
\neq $ $0,$ we compute%
\begin{eqnarray}
u_{y}(x,y) &=&-\frac{\partial x_{0}}{\partial y}y-x_{0}  \label{3.6.8} \\
&=&(\frac{\sin \theta (x_{0})}{s\theta ^{\prime }(x_{0})+\cos \theta (x_{0})}%
)s\cos \theta (x_{0})-x_{0}  \notag
\end{eqnarray}

\noindent by (\ref{3.6.4}) and (\ref{3.6.1}). For $(x,y)$ $\rightarrow $ $(%
\bar{x},0)$ with $\bar{x}$ $\leq $ $0,$ we have $s$ $\rightarrow $ $0$ and $%
x_{0}$ $\rightarrow $ $\bar{x}.$ Therefore $u_{y}(x,y)$ $\rightarrow $ $0$ $-%
\bar{x}$ $=$ $-\bar{x}$ by (\ref{3.6.8}) and (\ref{3.6.7}). For $(x,y)$ $%
\rightarrow $ $(\bar{x},0)$ with $\bar{x}$ $>$ $0,$ we have $x_{0}$ $%
\rightarrow $ $0,$ $\theta (x_{0})$ $\rightarrow $ $\theta (0)$ $=$ $\frac{%
\pi }{2},$ $\cos \theta (x_{0})$ $\rightarrow $ $0,$ $s$ $\rightarrow $ $%
\bar{x},$ and $\lim_{(x,y)\rightarrow (\bar{x},0^{+})}\theta ^{\prime
}(x_{0})$ $=$ $\theta ^{\prime }(0)$ $=$ $a$ $<$ $0$, $\lim_{(x,y)%
\rightarrow (\bar{x},0^{-})}\theta ^{\prime }(x_{0})$ $=$ $\theta ^{\prime
}(0)$ $=$ $-a$ $>$ $0$ by assumption. Observe that $(\frac{\sin \theta
(x_{0})}{s\theta ^{\prime }(x_{0})+\cos \theta (x_{0})})s$ in (\ref{3.6.8})
is bounded in the limit. So $\lim_{(x,y)\rightarrow (\bar{x},0)}$ $%
u_{y}(x,y) $ $=$ $0$ for $\bar{x}$ $>$ $0.$ Thus we have shown that $u_{y}$
is continuous at points of the $x$-axis in view of (\ref{3.6.5}). Altogether
on the whole plane $u_{x}$ and $u_{y}$ are $C^{0}$ and hence $u$ is $C^{1}$
smooth. Next we want to compute $D$ $:=$ $\sqrt{(u_{x}-y)^{2}+(u_{y}+x)^{2}}%
. $ From (\ref{3.6.6}), (\ref{3.6.1}), and $u$ $=$ $0$ on the $x$-axis we
learn that%
\begin{eqnarray}
u_{x}-y &=&(\frac{\cos \theta (x_{0})}{s\theta ^{\prime }(x_{0})+\cos \theta
(x_{0})}+1)s\cos \theta (x_{0})\text{ \ for }y\neq 0  \label{3.6.9} \\
u_{x}-y &=&0\text{ \ for }y=0.  \notag
\end{eqnarray}

\noindent From (\ref{3.6.8}), (\ref{3.6.1}), and (\ref{3.6.5}) we obtain%
\begin{eqnarray}
u_{y}+x &=&(\frac{\cos \theta (x_{0})}{s\theta ^{\prime }(x_{0})+\cos \theta
(x_{0})}+1)s\sin \theta (x_{0})\text{ \ for }y\neq 0  \label{3.6.10} \\
u_{y}+x &=&x\text{ \ for }y=0\text{ and }x>0;\text{ }=0\text{ for }y=0\text{
and }x\leq 0.  \notag
\end{eqnarray}

\noindent Therefore by (\ref{3.6.9}) and (\ref{3.6.10}) we have%
\begin{eqnarray}
D &=&(\frac{\cos \theta (x_{0})}{s\theta ^{\prime }(x_{0})+\cos \theta
(x_{0})}+1)s\text{ \ for }y\neq 0  \label{3.6.11} \\
D &=&x\text{ \ for }y=0\text{ and }x>0;\text{ }=0\text{ for }y=0\text{ and }%
x\leq 0.  \notag
\end{eqnarray}

\noindent It follows that $N^{\perp }$ $=$ $(u_{y}+x,$ $-(u_{x}-y))D^{-1}$ $%
= $ $(\sin \theta (x_{0}),-\cos \theta (x_{0}))$ for $y\neq 0$ and $N^{\perp
}$ $=$ $(1,$ $0)$ for $y=0$ and $x>0$ while $D$ $=$ $0$ for $y=0$ and $x\leq
0.$ We have shown that $\Gamma _{x_{0}}^{+},$ $\Gamma _{x_{0}}^{-},$ and the
positive $x$-axis are characteristic rays of $u$ while $(-\infty ,0]$ is the
singular set of $u.$ To verify Theorem B (a), we compute%
\begin{eqnarray}
D^{\prime } &=&\frac{\partial D}{\partial s}  \label{3.6.12} \\
&=&\frac{\cos \theta (x_{0})}{s\theta ^{\prime }(x_{0})+\cos \theta (x_{0})}%
+1-\frac{s\theta ^{\prime }(x_{0})\cos \theta (x_{0})}{(s\theta ^{\prime
}(x_{0})+\cos \theta (x_{0}))^{2}}  \notag
\end{eqnarray}

\noindent for $y$ $\neq $ $0$ and 
\begin{equation*}
D^{\prime }=\frac{\partial D}{\partial x}=1
\end{equation*}

\noindent for $y$ $=$ $0$ and $x$ $>$ $0$ by (\ref{3.6.11}). It follows from
(\ref{3.6.12}) that $D^{\prime }$ $\rightarrow $ $2$ as $s$ $\rightarrow $ $%
0 $ ($D$ $\rightarrow $ $0$ along $\Gamma _{x_{0}}^{+}$ or $\Gamma
_{x_{0}}^{-}).$ Since each singular point on the negative $x$-axis emits two
characteristic rays having equal angles with the positive $x$-axis, the
graph defined by $u$ and restricted to any bounded plane domain is a ($%
C^{1}) $ weak solution to (\ref{1.1}) (with $\vec{F}$ $=$ $(-y,x)$ and $H$ $%
= $ $0)$ and hence a $p$-minimizer in view of Theorem 6.3 and Theorem 3.3 in 
\cite{CHY1}.

\bigskip

The following two examples are inspired by \cite{Ri} (in which Ritor\'{e}
constructed locally Lipschitz continuous $p$-minimal graphs $(x,$ $y,$ $%
u(x,y))$ with finitely many singular half-lines emitted from a singular
point).

\bigskip

\textbf{Example 4.2. }We are going to construct a $C^{1}$ smooth $p$-minimal
graph $(x,$ $y,$ $u(x,y))$ with two singular half-lines emitted from a
singular point (see Figure 4.2 below).

\begin{figure}[ht]
\begin{center}
\includegraphics[width=6cm]{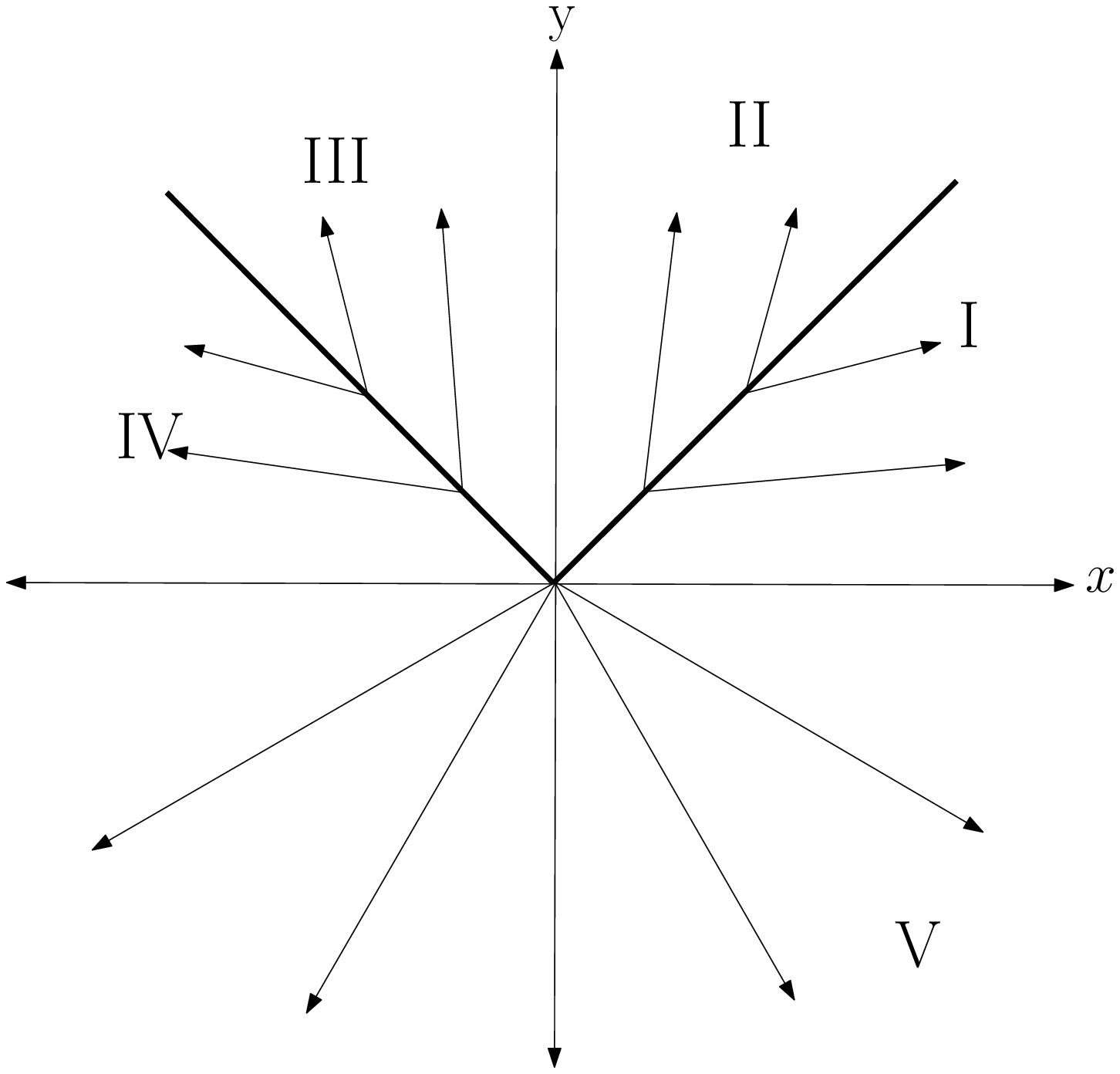}\\[0pt]
Figure 4.2
\end{center}
\par
\end{figure}

Let $\alpha $ $:$ $[0,\infty )$ $\rightarrow $ $[0,\frac{\pi }{4})$ be a $%
C^{1}$ smooth function such that $\alpha (0)$ $=$ $0,$ $\alpha ^{\prime
}(0^{+})$ $=$ $1,$ and $\alpha ^{\prime }(t)$ $>$ $0$ for all $t$ $\in $ $%
(0,\infty ).$ Let $\beta $ $:$ $[0,\infty )$ $\rightarrow $ $[0,\infty )$ be
a $C^{1}$ smooth function such that $\beta (0)$ $=$ $0,$ $\beta ^{\prime
}(0^{+})$ $=$ $0,$ $\beta ^{\prime }(t)$ $>$ $0$ for all $t$ $\in $ $%
(0,\infty ),$ and $\beta \rightarrow \infty $ as $t\rightarrow \infty .$ We
divide the plane into five regions (see Figure 4.3): region I :=\{ $x$ $>$ $%
y $ $>$ $0\},$ region II :=\{ $y$ $>$ $x$ $>$ $0\},$ region III :=\{ $y$ $>$ 
$-x$ $>$ $0\},$ region IV :=\{ $-x$ $>$ $y$ $>$ $0\},$ and region V :=\{ $y$ 
$<$ $0\}.$ In region I, we require the characteristic ray emitted from $%
(\beta (t),$ $\beta (t))$ to be parametrized by%
\begin{eqnarray}
x &=&s\cos \alpha (t)+\beta (t)  \label{4.3.1} \\
y &=&s\sin \alpha (t)+\beta (t)  \notag
\end{eqnarray}

\noindent for $s$ $>$ $0.$ The $C^{1}$ smooth map $\Psi $ $:$ $(s,t)$ $\in $ 
$(0,\infty )$ $\times $ $(0,\infty )\rightarrow $ $(x,y)$ $\in $ \textit{%
region I} is globally one to one and onto, hence $\Psi ^{-1}$ exists. We
compute the Jacobian matrix of $\Psi $ from (\ref{4.3.1}) as follows:%
\begin{eqnarray}
\frac{\partial (x,y)}{\partial (s,t)} &=&\left( 
\begin{array}{cc}
\frac{\partial x}{\partial s} & \frac{\partial x}{\partial t} \\ 
\frac{\partial y}{\partial s} & \frac{\partial y}{\partial t}%
\end{array}%
\right)  \label{4.3.2} \\
&=&\left( 
\begin{array}{cc}
\cos \alpha (t) & -s\alpha ^{\prime }(t)\sin \alpha (t)+\beta ^{\prime }(t)
\\ 
\sin \alpha (t) & s\alpha ^{\prime }(t)\cos \alpha (t)+\beta ^{\prime }(t)%
\end{array}%
\right) .  \notag
\end{eqnarray}

\noindent So we have the Jacobian%
\begin{equation}
\det \frac{\partial (x,y)}{\partial (s,t)}=s\alpha ^{\prime }(t)+\beta
^{\prime }(t)(\cos \alpha (t)-\sin \alpha (t))  \label{4.3.3}
\end{equation}

\noindent which is positive for $(s,t)$ $\in $ $(0,\infty )$ $\times $ $%
(0,\infty ).$ By the inverse function theorem, the map $\Psi ^{-1}$ is $%
C^{1} $ smooth and 
\begin{eqnarray}
\frac{\partial (s,t)}{\partial (x,y)} &=&\left( 
\begin{array}{cc}
\frac{\partial s}{\partial x} & \frac{\partial s}{\partial y} \\ 
\frac{\partial t}{\partial x} & \frac{\partial t}{\partial y}%
\end{array}%
\right)  \label{4.3.4} \\
&=&\frac{1}{s\alpha ^{\prime }+\beta ^{\prime }(\cos \alpha -\sin \alpha )}%
\left( 
\begin{array}{cc}
s\alpha ^{\prime }\cos \alpha +\beta ^{\prime } & s\alpha ^{\prime }\sin
\alpha -\beta ^{\prime } \\ 
-\sin \alpha & \cos \alpha%
\end{array}%
\right) .  \notag
\end{eqnarray}

Along a characteristic ray $u(s,t)$ is determined by $u(0,t)$ by integrating 
$du$ $+$ $xdy$ $-$ $ydx$ $=$ $0.$ We have%
\begin{eqnarray*}
u(s,t) &=&u(0,t)+\int_{0}^{s}(ydx-xdy) \\
&=&u(0,t)+s\beta (t)(\cos \alpha (t)-\sin \alpha (t)).
\end{eqnarray*}

\noindent by (\ref{4.3.1}). Set $u(0,t)$ $:=$ $0.$ That is, $u(\bar{x},\bar{x%
})$ $:=$ $0$ for $\bar{x}$ $\geq $ $0$ since $\bar{x}$ $=$ $\beta (t)$ for
some $t$ $\geq $ $0$ (note that (\ref{4.3.1}) also describes the boundary of
region I by enlarging the domain of $(s,t)$ to $[0,\infty )$ $\times $ $%
[0,\infty )).$ Now it is reasonable to define $u$ in the region I by%
\begin{equation}
u(s,t)=s\beta (t)(\cos \alpha (t)-\sin \alpha (t)).  \label{4.3.5}
\end{equation}

\noindent Similarly in region II (region III, resp.) we have%
\begin{eqnarray}
x &=&s\sin \alpha (t)+\beta (t)  \label{4.3.6} \\
&&\text{(}-s\sin \alpha (t)-\beta (t),\text{ resp.})  \notag \\
y &=&s\cos \alpha (t)+\beta (t)  \notag \\
u &=&s\beta (t)(\sin \alpha (t)-\cos \alpha (t))  \notag \\
&&(s\beta (t)(\cos \alpha (t)-\sin \alpha (t),\text{ resp.}).  \notag
\end{eqnarray}

\noindent and in region IV we have%
\begin{eqnarray}
x &=&-s\cos \alpha (t)-\beta (t)  \label{4.3.7} \\
y &=&s\sin \alpha (t)+\beta (t)  \notag \\
u &=&s\beta (t)(\sin \alpha (t)-\cos \alpha (t)).  \notag
\end{eqnarray}

\noindent We set%
\begin{equation}
u=0  \label{4.3.8'}
\end{equation}

\noindent on region V, the $x$-axis \{$y$ $=$ $0\},$ the half-lines \{$x$ $=$
$y$ $>$ $0\},$ $\{x$ $=$ $-y$ $<$ $0\},$ and the positive $y$-axis $\{x$ $=$ 
$0,$ $y$ $>$ $0\}.$ We can verify that $u$ defined by (\ref{4.3.5}), (\ref%
{4.3.6}), (\ref{4.3.7}), and (\ref{4.3.8'}) is $C^{1}$ smooth by showing the
continuity of $u_{x}$ and $u_{y}$. We leave the details to the reader as an
exercise. Moreover, it is a direct verification that the graph defined by $u$
is $p$-minimal with the expected singular set and characteristic lines.

Since each singular point on the half-lines $\{x$ $=$ $y$ $>$ $0\}$ and $\{x$
$=$ $-y$ $<$ $0\}$ emits two characteristic rays having equal angles with
the half-line, the graph defined by $u$ and restricted to any bounded planar
domain is a ($C^{1})$ weak solution to (\ref{1.1}) (with $\vec{F}$ $=$ $%
(-y,x)$ and $H$ $=$ $0)$ and hence a $p$-minimizer in view of Theorem 6.3
and Theorem 3.3 in \cite{CHY1}. On the lower half plane, we can easily see
from $u$ $=$ $0$ that each ray emitted from the origin is characteristic.
Compute $D$ $:=$ $\sqrt{(u_{x}-y)^{2}+(u_{y}+x)^{2}}$ $=$ $\sqrt{y^{2}+x^{2}}%
.$ It follows that along a characteristic ray, we have $N^{\perp }$ $=$ $%
\frac{\partial }{\partial r}$ and hence $D^{\prime }$ $=$ $\frac{\partial D}{%
\partial r}$ $=$ $1$ where $r$ :$=$ $\sqrt{x^{2}+y^{2}}.$ On the positive $y$%
-axis (which is a special characteristic ray also emitted from the origin),
we have $D$ $=$ $\sqrt{(0-y)^{2}+(0+0)^{2}}$ $=$ $y,$ $N^{\perp }$ $=$ $%
\frac{\partial }{\partial y},$ and hence $D^{\prime }$ $=$ $\frac{\partial D%
}{\partial y}$ $=$ $1.$ On the other hand, we claim that $D^{\prime }$ $%
\rightarrow $ $2$ as its argument tends to a singular point $(\bar{x},\bar{x}%
)$ or $(-\bar{x},\bar{x}),$ $\bar{x}$ $>$ $0,$ along a characteristic ray.
We check this for the characteristic rays in region I and leave the
remaining cases to the reader as an exercise. First observe that from (\ref%
{2.1}), (\ref{2.4}), $curl\vec{F}$ $=$ $2,$ and $\theta $ $=$ $\alpha $ $+$ $%
\frac{\pi }{2}$ in region I$,$ we have%
\begin{equation}
N\alpha =\frac{1}{D}(2-D^{\prime }).  \label{4.3.28}
\end{equation}

\noindent Note that in region I, we have $N$ $:=$ $(\cos \theta ,$ $\sin
\theta )$ $=$ $(-\sin \alpha ,$ $\cos \alpha ).$ Computing 
\begin{eqnarray*}
N\alpha &=&-\alpha ^{\prime }t_{x}\sin \alpha +\alpha ^{\prime }t_{y}\cos
\alpha \\
&=&\frac{\alpha ^{\prime }}{s\alpha ^{\prime }+\beta ^{\prime }(\cos \alpha
-\sin \alpha )}
\end{eqnarray*}%
\noindent by (\ref{4.3.4}),.we conclude that $N\alpha $ is bounded as $(x,y)$
$\rightarrow $ $(\bar{x},\bar{x}),$ $\bar{x}$ $>$ $0,$ while $s$ $%
\rightarrow $ $0,$ $t$ $\rightarrow $ $t_{0}$ (recall $\beta (t_{0})$ $=$ $%
\bar{x})$ along a characteristic ray in region I. So $D^{\prime }$ $%
\rightarrow $ $2$ from (\ref{4.3.28}) since $D$ $\rightarrow $ $0$ as its
argument tends to a singular point.

\bigskip

\textbf{Example 4.3}. Continuing the construction in Example 4.2, we are
going to build a $C^{1}$ smooth $p$-minimal graph $(x,$ $y,$ $u(x,y))$ with
three singular half-lines emitted from a singular point (see Figure 4.3
below).

\begin{figure}[ht]
\begin{center}
\includegraphics[width=6cm]{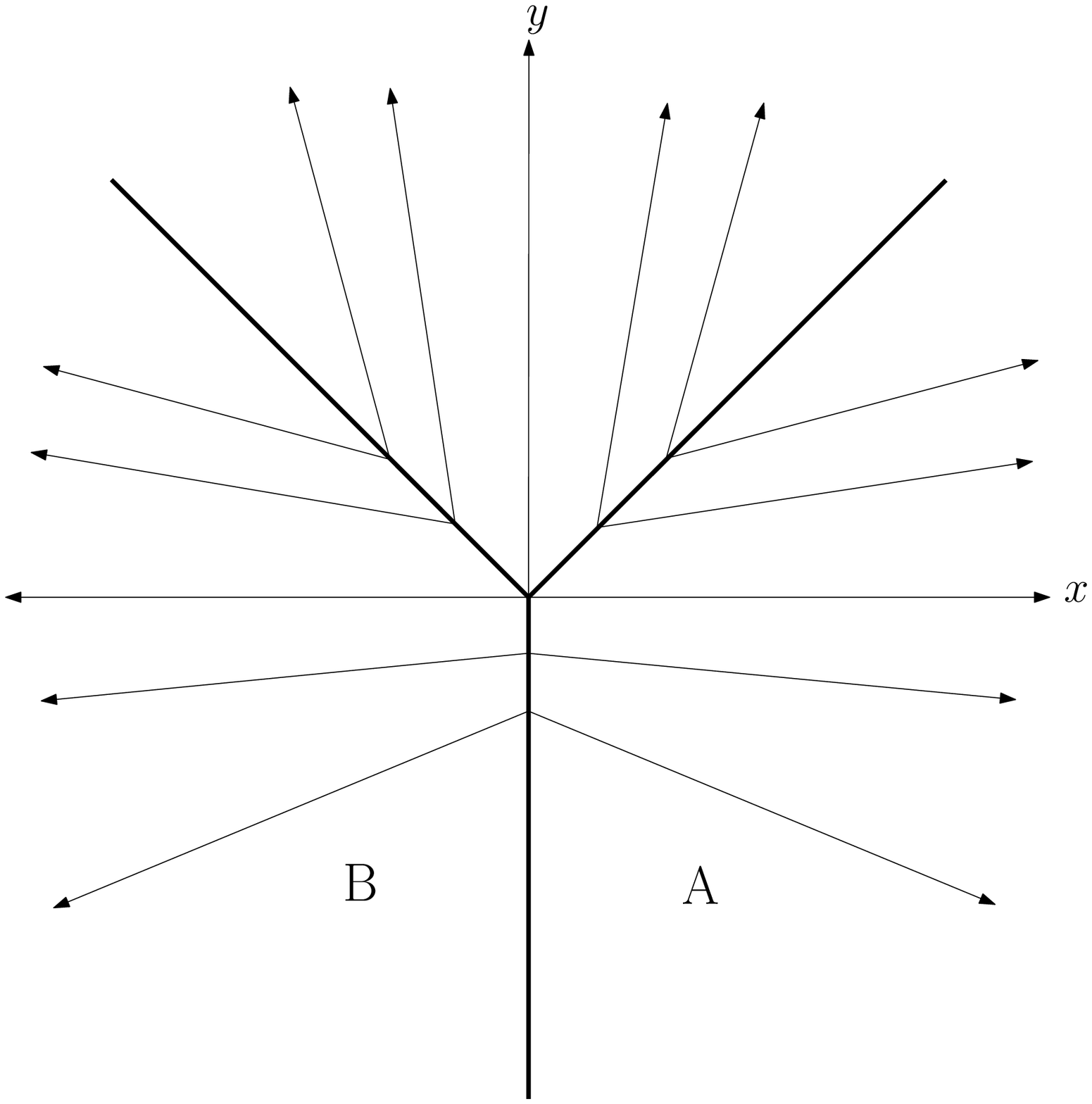}\\[0pt]
Figure 4.3
\end{center}
\par
\end{figure}

On the upper half plane and the $x$-axis, we define $u$ as in Example 4.2.
For the lower half plane, we divide it into two regions $A:=$ $\{x$ $>$ $0,$ 
$y$ $<$ $0\}$ and $B$ $:=$ $\{x$ $<$ $0,$ $y$ $<$ $0\}.$ In region A (region
B, resp.) we require the characteristic ray emitted from $(0,$ $-\beta (t))$
to be parametrized by%
\begin{eqnarray*}
x &=&s\cos \alpha (t) \\
\text{(}x &=&-s\cos \alpha (t),\text{ resp.)} \\
y &=&-s\sin \alpha (t)-\beta (t).
\end{eqnarray*}

\noindent for $(s,t)$ $\in $ $(0,\infty )$ $\times $ $(0,\infty ).$ As
before we can show that the $C^{1}$ smooth map $:$ $(s,t)$ $\rightarrow $ $%
(x,y)$ is a diffeomorphism from $(0,\infty )$ $\times $ $(0,\infty )$ onto
region A (region B, resp.). We then integrate $du$ $+$ $xdy$ $-$ $ydx$ $=$ $%
0 $ along the characteristic rays to get the expected $u$-value as follows:%
\begin{eqnarray*}
u(s,t) &=&u(0,t)+\int_{0}^{s}(ydx-xdy) \\
&=&u(0,t)-s\beta (t)\cos \alpha (t) \\
( &=&u(0,t)+s\beta (t)\cos \alpha (t),\text{ resp.).}
\end{eqnarray*}

\noindent Set $u$ $:=$ $0$ on the negative $y$-axis. It follows that $u(0,t)$
$=$ $0$ and hence we define $u$ on region A (region B, resp.) as

\begin{eqnarray*}
u(s,t) &=&-s\beta (t)\cos \alpha (t) \\
( &=&s\beta (t)\cos \alpha (t),\text{ resp.).}
\end{eqnarray*}

\noindent It is a direct verification that $u$ $\in $ $C^{1}$ and the graph
defined by $u$ is $p$-minimal with the expected singular set and
characteristic rays. Moreover, restricted to any bounded plane domain, $u$
is a ($C^{1})$ weak solution to (\ref{1.1}) (with $\vec{F}$ $=$ $(-y,x)$ and 
$H$ $=$ $0)$ and hence a $p$-minimizer in view of Theorem 6.3 and Theorem
3.3 in \cite{CHY1}. Along a characteristic ray in the upper half plane, we
have shown the behavior of $D^{\prime }$ in Example 4.2. By a similar
argument we can show that $D^{\prime }$ $\rightarrow $ $2$ as its argument
tends to a singular point $(0,\bar{y}),$ $\bar{y}$ $<$ $0,$ along a
characteristic ray in the lower half plane.

\bigskip

\textbf{Example 4.4. }We are going to construct a $C^{1}$ smooth $p$%
-minimizer $v$ defined on a bounded plane domain $\Omega $ having the line
segment $L\equiv \overline{(0,0),(0,1)}$ as the singular set (see Figure
4.4). The $p$-minimizer $v$ is $C^{\infty }$ smooth on $\Omega \backslash L$
and $C^{1}\backslash C^{2}$ on $L.$

Let $\alpha ,$ $\beta :[0,1]\rightarrow R$ be $C^{\infty }$ smooth functions
with the properties:

\begin{eqnarray}
\alpha (0) &=&\alpha (\frac{1}{2})=\alpha (1)=0,\alpha ^{\prime }(0)=\alpha
^{\prime }(1)=1,  \label{eqn6.7} \\
\alpha (t) &>&0\text{ for }0<t<\frac{1}{2},\text{ }\alpha (t)<0\text{ for }%
\frac{1}{2}<t<1,\text{ }|\alpha (t)|<\frac{\pi }{2},  \notag \\
\beta (0) &=&0,\beta (\frac{1}{2})=\frac{1}{2},\beta (1)=1,\beta ^{\prime
}(t)>0,\text{ }0<t<1,  \notag \\
\beta ^{(n)}(0) &=&\beta ^{(n)}(1)=0,\text{ }n=1,2,....  \notag
\end{eqnarray}

\noindent We then define a $p$-minimal surface in the parameters $s$ and $t$
as follows (note that $\theta (\tau )=\frac{\pi }{2}+\alpha (\tau )$ in
(4.9) of \cite{CHMY04}):

\begin{eqnarray}
x &=&s\cos \alpha (t)  \label{eqn6.8} \\
y &=&s\sin \alpha (t)+\beta (t)  \notag \\
z &=&s\beta (t)\cos \alpha (t)  \notag
\end{eqnarray}

\noindent for $(s,t)\in \lbrack 0,\infty )\times \lbrack 0,1].$ Compute the
Jacobian $J_{\varphi }$ of the map $\varphi :$ $(s,t)\rightarrow (x,y)$
given by (\ref{eqn6.8}):%
\begin{eqnarray}
J_{\varphi } &\equiv &\left| 
\begin{array}{c}
\frac{\partial x}{\partial s}\text{ \ }\frac{\partial x}{\partial t} \\ 
\frac{\partial y}{\partial s}\text{ \ }\frac{\partial y}{\partial t}%
\end{array}%
\right| =\left| 
\begin{array}{c}
\cos \alpha (t)\text{ \ \ \ \ \ \ \ \ }-s\alpha ^{\prime }(t)\sin \alpha (t)
\\ 
\sin \alpha (t)\text{ \ }s\alpha ^{\prime }(t)\cos \alpha (t)+\beta ^{\prime
}(t)%
\end{array}%
\right|  \label{eqn6.9} \\
&=&s\alpha ^{\prime }(t)+\beta ^{\prime }(t)\cos \alpha (t).  \notag
\end{eqnarray}

\noindent Observe that $J_{\varphi }>0$ in $\{[0,s_{0})\times \lbrack
0,1]\}\backslash $ $\{(0,0),(0,1)\}$ for a small positive number $s_{0}.$
Furthermore, for a possibly smaller positive number $s_{+},$ we can show
that $\varphi $ is a $C^{\infty }$ diffeomorphism (homeomorphism,
respectively) from $\{[0,s_{+})\times \lbrack 0,1]\}\backslash $ $%
\{(0,0),(0,1)\}$ ($[0,s_{+})\times \lbrack 0,1],$ respectively) into $R^{2}$
in view of the behavior of $\alpha $ near $t=0$, $1$ and the compactness of
a closed subinterval away from the boundary points $0,1$ of $(0,1).$ So $z$
in (\ref{eqn6.8}) can be viewed as a function of $x,y$ defining a $p$%
-minimal graph $z$ $=$ $u_{+}(x,y)$ over $\varphi ([0,s_{+})\times \lbrack
0,1])\subset R^{2}$. Similarly we describe another piece of $p$-minimal
surface by

\begin{eqnarray}
x &=&s\cos \alpha (t)  \label{eqn6.10} \\
y &=&-s\sin \alpha (t)+\beta (t)  \notag \\
z &=&s\beta (t)\cos \alpha (t)  \notag
\end{eqnarray}

\noindent for $(s,t)\in (-\infty ,0]\times \lbrack 0,1].$ Computing the
Jacobian $J_{\psi }$ of the map $\psi :$ $(s,t)\rightarrow (x,y)$ given by (%
\ref{eqn6.10}) in a similar way as in (\ref{eqn6.9}), we obtain 
\begin{equation*}
J_{\psi }=-s\alpha ^{\prime }(t)+\beta ^{\prime }(t)\cos \alpha (t).
\end{equation*}

\noindent Observe that $J_{\psi }>0$ in $\{(s_{1},0]\times \lbrack
0,1]\}\backslash $ $\{(0,0),(0,1)\}$ for some negative number $s_{1}.$ For a
negative number $s_{-}$ possibly closer to $0,$ we can show that $\psi $ is
a $C^{\infty }$ diffeomorphism (homeomorphism, respectively) from $%
\{(s_{-},0]\times \lbrack 0,1]\}\backslash $ $\{(0,0),(0,1)\}$ ($%
(s_{-},0]\times \lbrack 0,1],$ respectively) into $R^{2}.$ Thus we have
another piece of $p$-minimal graph $z$ $=$ $u_{-}(x,y)$ defined by (\ref%
{eqn6.10}) over $\psi ((s_{-},0]\times \lbrack 0,1])\subset R^{2}$.

Now we define a $p$-minimal graph $z$ $=$ $u(x,y)$ for $y\leq 0$ and $y\geq
1 $ by%
\begin{eqnarray}
u &=&0\text{ for }y\leq 0  \label{eqn6.11} \\
u &=&x\text{ for }y\geq 1,  \notag
\end{eqnarray}

\noindent $z$ $=$ $u_{+}(x,y)$ in $\varphi ([0,s_{+})\times \lbrack 0,1])$,
and $z$ $=$ $u_{-}(x,y)$ in $\psi ((s_{-},0]\times \lbrack 0,1]).$ Observe
that $u$ coincides with $u_{+}$ ($u_{-}$, respectively) on $[0,s_{+})\times
\{0,1\}$ ($(s_{-},0]\times \{0,1\}$, respectively) while $u_{+}$ $=$ $u_{-}$
on $L\equiv \{0\}\times \lbrack 0,1].$ Write $u_{+}$ ($u_{-}$, respectively)$%
=$ $x\beta (t)$ by (\ref{eqn6.8})((\ref{eqn6.10}), respectively). It follows
from $\beta ^{(k)}(0)$ $=$ $\beta ^{(k)}(1)$ $=$ $0$ in (\ref{eqn6.7}) that
for all positive integers $k,$%
\begin{eqnarray}
\frac{\partial ^{k}u_{+}}{\partial y^{k}} &=&0\text{ \ on }[0,s_{+})\times
\{0,1\}  \label{eqn6.12} \\
(\frac{\partial ^{k}u_{-}}{\partial y^{k}} &=&0\text{ \ on }(s_{-},0]\times
\{0,1\},\text{ respectively).}  \notag
\end{eqnarray}

\noindent So $u$ matches with $u_{+}$ ($u_{-}$, respectively) on $%
[0,s_{+})\times \{0,1\}$ ($(s_{-},0]\times \{0,1\}$, respectively) $%
C^{\infty }$ smoothly by (\ref{eqn6.12}). On $L$ (where $x$ $=$ $0$), we have%
\begin{equation}
\frac{\partial u_{+}}{\partial x}=\frac{\partial u_{-}}{\partial x}=\beta
(t).  \label{eqn6.13}
\end{equation}

\noindent Differentiating the first equation of (\ref{eqn6.8}) or (\ref%
{eqn6.10}) with respect to $x$ at a point on $L$ (where $s=0$) gives 
\begin{equation}
\frac{\partial s}{\partial x}=\frac{1}{\cos \alpha (t)}.  \label{eqn6.14}
\end{equation}

\noindent Differentiating the second equation of (\ref{eqn6.8}) and (\ref%
{eqn6.10}) with respect to $x$ at a point on $L$ (where $s=0$) $\backslash
\{(0,0),$ $(0,1)\}$ and substituting (\ref{eqn6.14}) into the resulting
formulas, we obtain

\begin{equation}
\frac{\partial t}{\partial x}=\mp \frac{\tan \alpha (t)}{\beta ^{\prime }(t)}
\label{eqn6.15}
\end{equation}

\noindent for the \textquotedblright $\pm $ sides\textquotedblright ,
respectively. Now computing the second derivative of $u_{\pm }$ in the $x$
direction on $L$ (where $x$ $=$ $0$ or $s$ $=$ $0$)$,$ we get%
\begin{equation}
\frac{\partial ^{2}u_{\pm }}{\partial x^{2}}=2\beta ^{\prime }(t)\frac{%
\partial t}{\partial x}=\mp 2\tan \alpha (t)  \label{eqn6.16}
\end{equation}

\noindent by (\ref{eqn6.15}) (for $0<t<1,$ but the final result also holds
for $t=0$ and $1).$ It follows from (\ref{eqn6.16}) and (\ref{eqn6.7}) that

\begin{equation}
\frac{\partial ^{2}u_{+}}{\partial x^{2}}\neq \frac{\partial ^{2}u_{-}}{%
\partial x^{2}}  \label{eqn6.17}
\end{equation}

\noindent on $L\backslash (\{0\}\times \{0,\frac{1}{2},1\})$ while at $%
(0,0), $ $(0,\frac{1}{2}),$ and $(0,1),$ we have%
\begin{equation}
\frac{\partial ^{2}u_{+}}{\partial x^{2}}=\frac{\partial ^{2}u_{-}}{\partial
x^{2}}.  \label{eqn6.18}
\end{equation}

We define $v$ to be $u,$ $u_{+},$ $u_{-}$ on the corresponding domains. Glue
a patch of suitable domain around $(0,0)$ from $\{y<0\}$ and a patch of
suitable domain from $\{y>1\}$ to $\varphi ([0,s_{+})\times \lbrack 0,1])$ $%
\cup $ $\psi ((s_{-},0]\times \lbrack 0,1])$ to form a $C^{\infty }$ smooth
bounded domain $\Omega $ (see Figure 4.4 below)$.$

\begin{figure}[ht]
\begin{center}
\includegraphics[width=6cm]{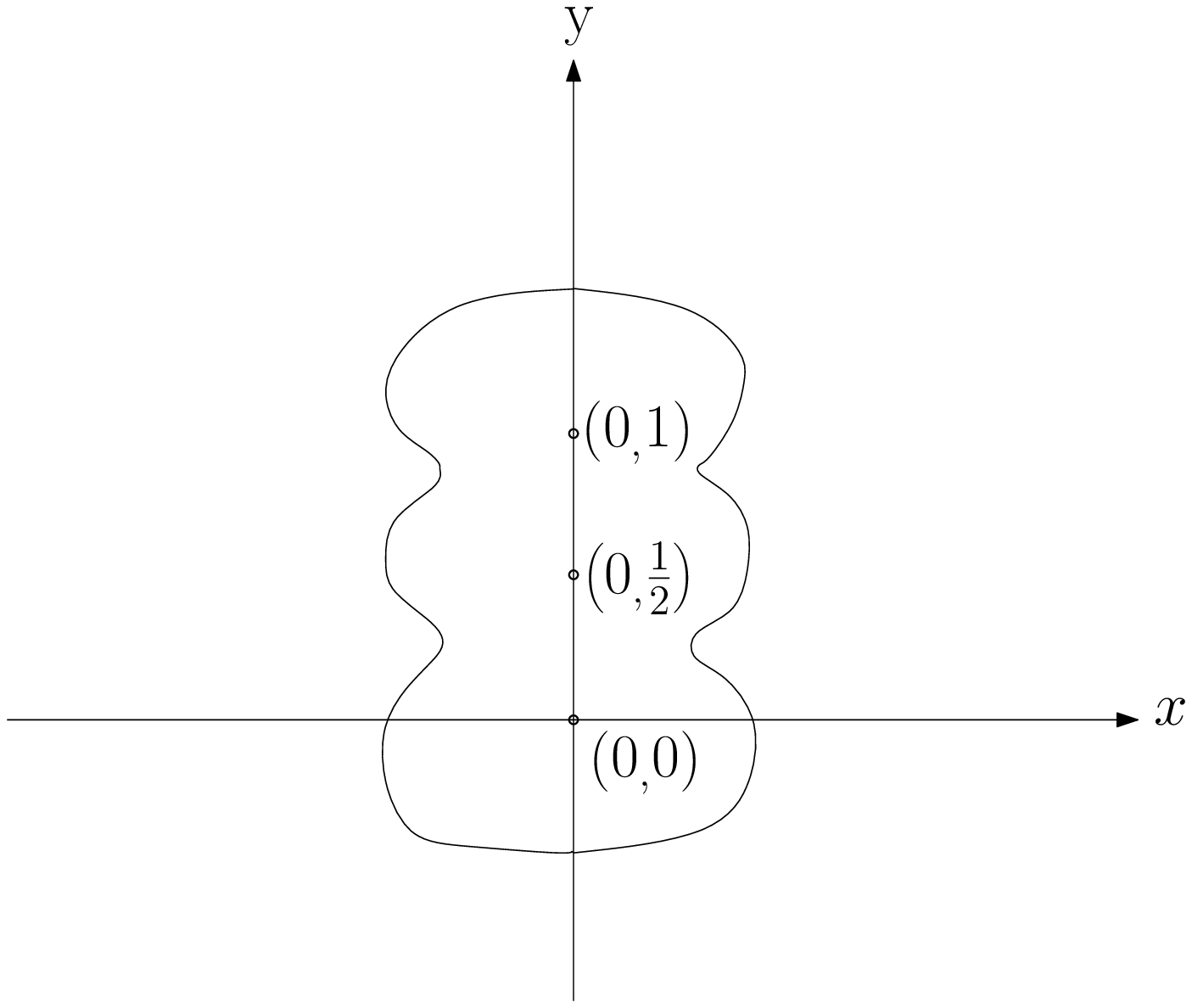}\\[0pt]
Figure 4.4
\end{center}
\par
\end{figure}


Consider $v$ restricted to $\Omega .$ From (\ref{eqn6.13}), (\ref{eqn6.17}),
and from similar arguments as in Example 4.2, we learn that $v$ is $%
C^{1}\backslash C^{2}$ on $L.$ Moreover, $v$ is $C^{\infty }$ smooth on $%
\bar{\Omega}\backslash L$ by (\ref{eqn6.12}). In view of the extension
theorem (Proposition 3.5 in \cite{CHMY04}) of characteristic lines for a $%
C^{2}$ smooth solution and $J_{\varphi }>0$ (see (\ref{eqn6.9}))$,$ $J_{\psi
}>0$, we can easily show that $L$ is the only singular set in $\Omega .$
Observe that $0$ and $x$ are solutions to the $p$-minimal surface equation
(see (\ref{1.1}) with $\vec{F}$ $=$ $(-y,x)$ and $H$ $=$ $0$). It then
follows that $v$ is a $p$-minimizer in view of Theorem 6.3 and Theorem 3.3
in \cite{CHY1}. We have proved our claim in the beginning of this example.

\bigskip

\textbf{Example 4.5. } According to Theorem C, if a singular point is not
isolated, then it emits at least one $C^{0}$ singular curve. One may ask if
there are only finitely many such singular curves. The configuration of
singular lines in Figure 4.5 (b) shows that it is possible to have infinitly
many singular curves emitted from a singular point and shrinking to this
singular point. On the other hand, it is not possible to have a
configuration of singular lines as shown in Figure 4.5 (a), which converges
to another singular line of length $>0$. The reason is that any
characteristic line emitted from a point in this limit singular line must
span an angle, and hence hit an approaching singular line, contradicting
Theorem B (c).

\begin{figure}[ht]
\begin{center}
\includegraphics[width=6cm]{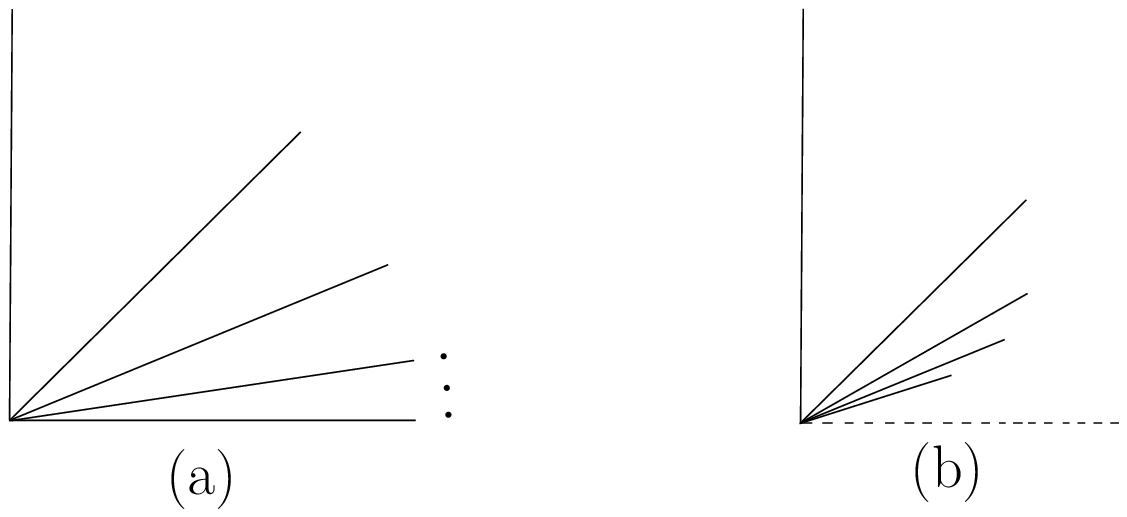}\\[0pt]
Figure 4.5
\end{center}
\par
\end{figure}

\bigskip

\section{Size and regularity of the singular set}

In this section we first study the size of the singular set.

\proof
\textbf{(of Theorem F) }It is enough to show that for any $p$ $\in $ $S_{%
\vec{F}}(u)$ $\cap $ $\Omega $, there exists $r$ $>$ $0$ such that $B_{r}(p)$
$\subset $ $\Omega $ and $\mathcal{H}^{2}(S_{\vec{F}}(u)\cap B_{r}(p))$ $=$ $%
0.$ First take $r_{0}$ $>$ $0$ such that $B_{r_{0}}(p)$ $\subset \subset $ $%
\Omega .$ Take%
\begin{equation*}
r=\frac{1}{2}\min \{\frac{1}{\max_{\bar{B}_{r_{0}}(p)}|H|},r_{0}\}.
\end{equation*}

\noindent It is clear that $B_{r}(p)$ $\subset $ $B_{2r}(p)$ $\subset $ $%
B_{r_{0}}(p)$ $\subset \subset $ $\Omega .$ Next we are going to show $%
\mathcal{H}^{2}(S_{\vec{F}}(u)\cap B_{r}(p))$ $=$ $0.$

By Proposition $3.3^{\prime }$ each $q$ $\in $ $S_{\vec{F}}(u)$ $\cap $ $%
B_{r}(p)$ emits a characteristic curve $Lq$ which hits the boundary $%
\partial B_{2r}(p)$ of a bigger ball $B_{2r}(p)$ at $q\prime $ transversely
(noting that $\pm H$ is the line curvature of $Lq)$. Since $q\prime $ must
be nonsingular by Theorem B(c), we can find a small open interval $%
I_{q^{\prime }},$ contained in $\partial B_{2r}(p)$ and consisting of
nonsingular points, such that the characteristic curves emitted from $%
I_{q^{\prime }}$ hit $S_{\vec{F}}(u)$ in a set $S_{q}$. Note that each
characteristic curve emitted from $I_{q^{\prime }}$ must hit a singular
point in $B_{2r}(p)$ by Lemma $3.2^{\prime }.$ Since the number of disjoint
open intervals (take the union if two $I_{q^{\prime }}$ overlap) in $%
\partial B_{2r}(p)$ is countable, we need only to show $\mathcal{H}^{2}$($%
S_{q}$) $=$ $0$ (two $Lq^{\prime }s$ do not hit the same $q^{\prime };$
otherwise, $q^{\prime }$ becomes singular). Since $I_{q^{\prime }}$ is
transverse to characteristic curves, we may use a parameter $\tau $ such
that $I_{q^{\prime }}$ $=$ $\beta (\tau )$ $(\in $ $C^{1}$ for $\tau $ $\in $
$(\tau _{0},$ $\tau _{1}),$ say) to describe $S_{q}$ as a (folded) graph $%
G_{q}$ $=$ $\{(\tau ,$ $\sigma _{1}(\tau )\}.$ Here we parametrize
characteristic curves by the arc length $\sigma $ such that $\sigma $ $=$ $0$
for $I_{q^{\prime }}.$ We have known that $\sigma _{1}(\tau )$ is $C^{0}$ in 
$\tau .$ Define a map $\varphi $ : $(\tau ,\sigma )$ $\rightarrow $ $(x,y)$
such that ($x(\tau ,\sigma )$, $y(\tau ,\sigma ))$ describes a
characteristic curve for each $\tau $ with the initial data ($x(\tau ,0)$, $%
y(\tau ,0))$ $=$ $\beta (\tau ).$ Namely we have%
\begin{eqnarray}
\frac{\partial x(\tau ,\sigma )}{\partial \sigma } &=&\sin \theta (x(\tau
,\sigma ),y(\tau ,\sigma ))  \label{5.6'} \\
\frac{\partial y(\tau ,\sigma )}{\partial \sigma } &=&-\cos \theta (x(\tau
,\sigma ),y(\tau ,\sigma ))  \notag
\end{eqnarray}

\noindent with $x(\tau ,0)$ $=$ $x_{0}(\tau ),$ $y(\tau ,0)$ $=$ $y_{0}(\tau
)$ where $\beta (\tau )$ $=$ $(x_{0}(\tau ),$ $y_{0}(\tau )).$ Since $\theta 
$ is $C^{1}$ smooth by Theorem D in \cite{CHY2}, we conclude that the
solution to (\ref{5.6'}) is $C^{1}$ smooth in the parameter $\tau $ for $%
C^{1}$ smooth initial data $x_{0}(\tau ),$ $y_{0}(\tau ).$ Hence the map $%
\varphi $ is $C^{1}$ smooth (for $\tau $ $\in $ $(\tau _{0},$ $\tau _{1})$
and $\sigma $ $\in $ $(0,$ $\sigma _{1}(\tau ))).$ To extend $\varphi $
beyond $G_{q}$ so that $\varphi (G_{q})$ $=$ $S_{q}$, we consider $\theta $
as an independent variable and ($x(\tau ,\sigma )$, $y(\tau ,\sigma ))$ to
be the first two components of the unique solution to the following
equations:%
\begin{eqnarray}
\frac{dx}{d\sigma } &=&\sin \theta ,\text{ \ }\frac{dy}{d\sigma }=-\cos
\theta  \label{5.7'} \\
\frac{d\theta }{d\sigma } &=&-H(x,y)  \notag
\end{eqnarray}

\noindent with the initial condition $x(0)$ $=$ $x_{0}(\tau ),$ $y(0)$ $=$ $%
y_{0}(\tau ),$ and $\theta (0)$ $=$ $\theta _{0}(\tau )$ $\in $ $C^{1}$.
Note that the third equation of (\ref{5.7'}) can be deduced if $\theta $ is
the angular function associated to the horizontal normal $N$ (see Theorem A
in \cite{CHY2}). By Theorem 3.1 on page 95 in \cite{Ha}, $x(\tau ,\sigma )$
and $y(\tau ,\sigma )$ are $C^{1}$ smooth in $\tau $ (and $\sigma $ of
course), and apparently they are defined on an open neighborhood of $G_{q}.$
That is, the map $\varphi $ extends $C^{1}$ smoothly over $G_{q}.$ Let $%
J(\varphi )$ denote the Jacobian of $\varphi .$ Clearly we have%
\begin{equation*}
|J(\varphi )|\leq C_{j}
\end{equation*}%
\noindent on a compact set $K_{j}$ where $\cup _{j=1}^{\infty }K_{j}$
exhausts the whole domain. It follows (e.g., 2.10.11 in \cite{Fe}) that%
\begin{equation}
\mathcal{H}^{2}(\varphi (G_{q}\cap K_{j}))\leq C_{j}\mathcal{H}%
^{2}(G_{q}\cap K_{j}).  \label{5.8'}
\end{equation}

\noindent On the other hand, we have 
\begin{equation*}
\mathcal{H}^{2}(G_{q}\cap K_{j})=\int \mathcal{H}^{1}((G_{q}\cap K_{j})\cap
\{\tau =c\})dc
\end{equation*}%
\noindent by Fubini's theorem. But $(G_{q}\cap K_{j})\cap \{\tau =c\}$
consists of only one point ($c$, $\sigma _{1}(c)$), and hence $\mathcal{H}%
^{1}((G_{q}\cap K_{j})\cap \{\tau =c\})$ $=$ $0$. So $\mathcal{H}%
^{2}(G_{q}\cap K_{j})$ $=$ $0$. It follows from (\ref{5.8'}) that \noindent $%
\mathcal{H}^{2}(\varphi (G_{q}\cap K_{j}))$ $=$ $0.$ We then have%
\begin{eqnarray*}
\mathcal{H}^{2}(S_{q}) &=&\mathcal{H}^{2}(\varphi (G_{q})) \\
&=&\mathcal{H}^{2}(\varphi (\cup _{j=1}^{\infty }(G_{q}\cap K_{j}))) \\
&=&\mathcal{H}^{2}(\cup _{j=1}^{\infty }\varphi (G_{q}\cap K_{j})) \\
&\leq &\sum_{j=1}^{\infty }\mathcal{H}^{2}(\varphi (G_{q}\cap K_{j}))=0.
\end{eqnarray*}

\noindent Therefore $\mathcal{H}^{2}(S_{q})$ $=$ $0.$

\endproof%

\bigskip

From now on, we are going to study the regularity of nondegenerate crack
points (see Definition 5.1 and Definition 5.2). Let $N_{\varepsilon }$
denote a mollification of $N.$

\bigskip

\textbf{Lemma 5.1}. \textit{Suppose }$\mathcal{H}^{2}(K)$\textit{\ }$=$%
\textit{\ }$0$\textit{\ for a subset }$K$\textit{\ of }$\Omega .$ \textit{%
Let }$\theta $\textit{\ }$\in $\textit{\ }$C^{0}(\Omega \backslash K)$%
\textit{\ be a weak solution to (\ref{1.8.1}) with }$H$ $\in $ $C^{0}(\Omega
)$\textit{.} \textit{Let }$\Omega ^{\prime }$\textit{\ }$\subset \subset $%
\textit{\ }$\Omega $\textit{\ be a Lipschitzian domain. Then we have}%
\begin{equation}
\lim_{\varepsilon \rightarrow 0}\doint\limits_{\mathcal{\partial }\Omega
^{\prime }}N_{\varepsilon }\cdot \nu =\int_{\Omega ^{\prime }}H  \label{5.6}
\end{equation}%
\textit{\noindent where }$\nu $\textit{\ is the outward unit normal to }$%
\partial \Omega ^{\prime }$\textit{\ (defined a.e. with respect to the
boundary measure).}

\bigskip

\proof
From the definition of weak solution, we have%
\begin{equation}
\int_{\Omega }N\cdot \nabla \varphi +\int_{\Omega }H\varphi =0  \label{5.7}
\end{equation}

\noindent for all $\varphi $ $\in $ $C_{0}^{\infty }(\Omega ).$ Let $\varphi
_{\varepsilon }$ ($V_{\varepsilon },$ resp.) denote a mollification of $%
\varphi $ (of a vector field $V,$ resp.)$.$ Take a bounded domain $\Omega
^{\prime \prime }$ such that $\Omega ^{\prime }$\textit{\ }$\subset $ $%
\Omega ^{\prime \prime }$ $\subset \subset $\textit{\ }$\Omega .$ Let $%
\varphi $ $\in $ $C_{0}^{\infty }(\Omega ^{\prime \prime }).$ Then from (\ref%
{5.7}) we have%
\begin{eqnarray}
0 &=&\int_{\Omega ^{\prime \prime }}N\cdot \nabla \varphi _{\varepsilon
}+\int_{\Omega ^{\prime \prime }}H\varphi _{\varepsilon }\text{ (}\varphi
_{\varepsilon }\in C_{0}^{\infty }(\Omega ^{\prime \prime })\text{ for }%
\varepsilon \text{ small)}  \label{5.8} \\
&=&\int_{\Omega ^{\prime \prime }}N\cdot (\nabla \varphi )_{\varepsilon
}+\int_{\Omega ^{\prime \prime }}H\varphi _{\varepsilon }  \notag \\
&=&\int_{\Omega ^{\prime \prime }}N_{\varepsilon }\cdot \nabla \varphi
+\int_{\Omega ^{\prime \prime }}H_{\varepsilon }\varphi  \notag
\end{eqnarray}

\noindent (note that $|N|$ $=$ $1$, $H$ $\in $ $L_{loc}^{1}(\Omega ),$ and
hence both $N$ and $H$ are in $L^{1}(\Omega ^{\prime \prime })).$ Here we
have used the fact that $\int gf_{\varepsilon }=\int g_{\varepsilon }f$ for $%
g$ $\in $ $L^{1}$ and $f$ $\in $ $C_{0}^{\infty }$ (which can be easily
proved by Fubini's theorem). It follows from (\ref{5.8}) that%
\begin{equation}
\func{div}N_{\varepsilon }=H_{\varepsilon }\text{ in }\Omega ^{\prime \prime
}.  \label{5.9}
\end{equation}

We can now integrate (\ref{5.9}) over $\Omega ^{\prime }$ and apply the
divergence theorem to obtain%
\begin{equation}
\doint\limits_{\partial \Omega ^{\prime }}N_{\varepsilon }\cdot \nu
=\int_{\Omega ^{\prime }}H_{\varepsilon }.  \label{5.10}
\end{equation}

\noindent Since $H_{\varepsilon }$ is bounded on $\Omega ^{\prime }$ for $%
\varepsilon $ small, we have $\lim_{\varepsilon \rightarrow 0}\int_{\Omega
^{\prime }}H_{\varepsilon }$ $=$ $\int_{\Omega ^{\prime }}H$ by Lebesgue's
dominated convergence theorem. From this and (\ref{5.10}), we get (\ref{5.6}%
).

\endproof%

\bigskip

We remark that Lemma 5.1 holds even if $\mathcal{H}^{1}(\partial \Omega
^{\prime }\cap K)$ $\neq $ $0$ where $\mathcal{H}^{1}$ denotes the
1-dimensional Hausdorff measure$.$ Since $N$ is not defined on $\partial
\Omega ^{\prime }$ $\cap $ $K,$ $N_{\varepsilon }$ may not converge to $N$
a.e. on $\partial \Omega ^{\prime }.$ Now suppose that we are in the
situation of Lemma 5.1 with $H$ $\in $ $C^{1}(\Omega )$. Take $p$ $\in $ $K$%
. Suppose

(a) there is an open neighborhood $U$ $\subset $ $\Omega $ of $p$ such that $%
U\cap K$ is a $C^{0}$ curve $\gamma $ (passing through $p)$ dividing $U$
into two connected regions $U^{+}$ and $U^{-}$ and

(b) for any $q$ $\in $ $\gamma $ there are exactly two characteristic curves 
$\Gamma _{q}^{+},$ $\Gamma _{q}^{-}$ issuing from $q,$ such that $\Gamma
_{q}^{+}$ $\subset $ $U^{+}$ and $\Gamma _{q}^{-}$ $\subset $ $U^{-}.$

Take $q^{+}$ $\in $ $\Gamma _{q}^{+}\backslash \{q\}$ ($q^{-}$ $\in $ $%
\Gamma _{q}^{-}\backslash \{q\},$ resp.). There passes a seed curve (i.e.,
integral curve of $N$) $\gamma _{+}$ ($\gamma _{-},$ resp.) each point of
which emits a characteristic curve in $U^{+}$ ($U^{-},$ resp.) hitting a
point in a neighborhood of $p$ in $\gamma .$ We parametrize $\gamma _{+}$ ($%
\gamma _{-},$ resp.) by the arc-length parameter $\tau ^{+}$ ($\tau ^{-},$
resp.) in some open interval including $\tau _{0}^{+}$ ($\tau _{0}^{-},$
resp.) so that $\frac{\partial }{\partial \tau ^{+}}$ ($\frac{\partial }{%
\partial \tau ^{-}},$ resp.) coincides with $\pm N$ and the characteristic
curve issuing from $\gamma _{+}(\tau ^{+})$ ($\gamma _{-}(\tau ^{-}),$
resp.) hits $\gamma $ at $X_{+}(\tau ^{+})$ ($X_{\_}(\tau ^{-}),$ resp.). We
assume that $X_{+}(\tau _{0}^{+})$ $=$ $X_{\_}(\tau _{0}^{-}),$ denoted as $%
p_{0}.$ Consider the following system of ordinary differential equations:%
\begin{eqnarray}
\frac{dx}{d\sigma } &=&\sin \theta ,\text{ }\frac{dy}{d\sigma }=-\cos \theta
\label{5.10.1} \\
\frac{d\theta }{d\sigma } &=&-H.  \notag
\end{eqnarray}%
\noindent Note that the first two equations of (\ref{5.10.1}) describe the
characteristic curves and $\sigma $ is the unit-speed parameter (on the $xy$%
-plane). Let $P_{\pm }(\sigma ,\tau ^{\pm })$ :$=$ $(x_{\pm }(\sigma ,\tau
^{\pm }),$ $y_{\pm }(\sigma ,\tau ^{\pm }))$ where $x_{\pm }(\sigma ,\tau
^{\pm })$, $y_{\pm }(\sigma ,\tau ^{\pm }),$ and $\theta _{\pm }(\sigma
,\tau ^{\pm })$ are the solutions to (\ref{5.10.1}) with $(x_{\pm }(0,\tau
^{\pm }),$ $y_{\pm }(0,\tau ^{\pm }))$ $=$ $\gamma _{\pm }(\tau ^{\pm })$
and $\theta _{\pm }(0,\tau ^{\pm })$ $=$ $\theta (\gamma _{\pm }(\tau ^{\pm
})).$ Write%
\begin{eqnarray}
X_{+}(\tau ^{+}) &=&P_{+}(\sigma _{+}(\tau ^{+}),\tau ^{+})  \label{5.11} \\
X_{\_}(\tau ^{-}) &=&P_{-}(\sigma _{-}(\tau ^{-}),\tau ^{-})  \notag
\end{eqnarray}%
\noindent for some continuous functions $\sigma _{\pm }.$

\bigskip

\textbf{Definition 5.1. }Suppose $H^{2}(K)$\ $=$\ $0$\ for a subset $K$\ of $%
\Omega .$ Let $\theta $\ $\in $\ $C^{0}(\Omega \backslash K)$\ be a weak
solution to (\ref{1.8.1}) with $H$ $\in $ $C^{1}(\Omega )$. Take $p$ $\in $ $%
K$. Suppose conditions (a) and (b) above hold. Take $q^{\pm }$ $\in $ $%
\Gamma _{q}^{\pm }\backslash \{q\}$. We have seed curves $\gamma _{\pm },$
parametrized by $\tau ^{\pm },$ $X_{\pm }(\tau ^{\pm }),$ and $X_{+}(\tau
_{0}^{+})$ $=$ $X_{\_}(\tau _{0}^{-})$ $=$ $p_{0}$ as above. We call $p_{0}$
nondegenerate if the Jacobian of $P_{\pm }$ (see (\ref{5.11})) satisfies%
\begin{equation}
\det \frac{\partial (x_{\pm }(\sigma ,\tau ^{\pm }),y_{\pm }(\sigma ,\tau
^{\pm }))}{\partial (\sigma ,\tau ^{\pm })}\neq 0  \label{5.11.1}
\end{equation}%
\noindent at $(\sigma _{\pm }(\tau _{0}^{\pm }),$ $\tau _{0}^{\pm })$ (note
that $P_{\pm }$ is defined on a neighborhood of $(\sigma _{\pm }(\tau
_{0}^{\pm }),$ $\tau _{0}^{\pm })$ by the basic ODE theory).

\bigskip

Note that Definition 5.1 is independent of the choice of seed curves $\gamma
_{\pm }$ by observing that ($\vec{u}\wedge \vec{v}$ $:=$ $u_{1}v_{2}$ $-$ $%
u_{2}v_{1}$ for $\vec{u}$ $=$ $(u_{1},u_{2})$ and $\vec{v}$ $=$ $%
(v_{1},v_{2}))$%
\begin{eqnarray}
&&\frac{\partial P_{\pm }}{\partial \sigma }\wedge \frac{\partial P_{\pm }}{%
\partial \tau ^{\pm }}  \label{5.11.1a} \\
&=&\frac{\partial P_{\pm }}{\partial \tilde{\sigma}}\wedge (\frac{\partial
P_{\pm }}{\partial \tilde{\tau}^{\pm }}\frac{\partial \tilde{\tau}^{\pm }}{%
\partial \tau ^{\pm }}+\frac{\partial P_{\pm }}{\partial \tilde{\sigma}}%
\frac{\partial \tilde{\sigma}}{\partial \tau ^{\pm }})  \notag \\
&=&(\frac{\partial P_{\pm }}{\partial \tilde{\sigma}}\wedge \frac{\partial
P_{\pm }}{\partial \tilde{\tau}^{\pm }})\frac{\partial \tilde{\tau}^{\pm }}{%
\partial \tau ^{\pm }}  \notag
\end{eqnarray}%
\noindent where $\tilde{\tau}^{\pm }$ are arc-length parameters for another
choice of seed curves $\tilde{\gamma}_{\pm }$ and $\tilde{\sigma}$ is the
corresponding unit-speed parameter for characteristic curves. It follows
that (\ref{5.11.1}) holds with respect to $(\tilde{\sigma},$ $\tilde{\tau}%
^{\pm })$ since $\frac{\partial \tilde{\tau}^{\pm }}{\partial \tau ^{\pm }}$ 
$\neq $ $0$ in the above equality. Now assume that $p_{0}$ is nondegenerate.
Then by the inverse function theorem, $P_{\pm }$ is a local $C^{1}$
diffeomorphism from $(\sigma ,\tau ^{\pm })$ to $(x_{\pm }(\sigma ,\tau
^{\pm }),y_{\pm }(\sigma ,\tau ^{\pm }))$ near $(\sigma _{\pm }(\tau
_{0}^{\pm }),$ $\tau _{0}^{\pm }).$ So $N_{\pm }$ :$=$ $(\cos \theta _{\pm
}, $ $\sin \theta _{\pm })$ is well defined near $p_{0}$ on the $xy$-plane.
Since $\gamma _{\pm }$ is $C^{1}$ smooth in $\tau ^{\pm },$ $\theta _{\pm }$
is $C^{1}$ smooth in $\sigma $ and $\tau ^{\pm }$ by the ODE theory, and
hence $C^{1}$ smooth near $p_{0}$ on the $xy$-plane. Observe that $N_{\pm }$
satisfies the equation 
\begin{equation}
\func{div}N_{\pm }=H  \label{5.11.2}
\end{equation}%
\noindent near $p_{0},$ say$,$ in a neighborhood $U^{\prime }$ of $p_{0}$
since 
\begin{eqnarray*}
\func{div}N_{\pm } &=&(\cos \theta _{\pm })_{x}+(\sin \theta _{\pm })_{y} \\
&=&-(\sin \theta _{\pm })(\theta _{\pm })_{x}+(\cos \theta _{\pm })(\theta
_{\pm })_{y} \\
&=&-\frac{d\theta _{\pm }}{d\sigma }=H
\end{eqnarray*}

\noindent by (\ref{5.10.1}).

\bigskip

\textbf{Lemma 5.2}. \textit{Suppose we are in the situation of Definition
5.1. In particular, assume }$p_{0}$\textit{\ is nondegenerate}. \textit{Let }%
$N_{\pm }$\textit{\ be defined in a neighborhood }$U^{\prime }$ \textit{of }$%
p_{0},$\textit{\ satisfying (\ref{5.11.2}) as above. Let }$\Upsilon $\textit{%
\ }$\subset $\textit{\ }$U^{\prime }$\textit{\ be a }$C^{1}$\textit{\ smooth
arc joining }$p_{0}$\textit{\ to }$\bar{p}$\textit{\ }$\in $\textit{\ }$%
\gamma $\textit{\ }$\cap $\textit{\ }$U^{\prime }$\textit{. Then }%
\begin{equation}
\int_{\Upsilon }(N_{+}-N_{-})\cdot \nu =0.  \label{5.11.3}
\end{equation}

\noindent \textit{where }$\nu $\textit{\ is the unit normal to }$\Upsilon $ 
\textit{such that }$\nu $\textit{\ and }$\Upsilon ^{\prime }$\textit{\ are
positively oriented.}

\bigskip

\proof
Take $C^{1}$ smooth arcs $\Upsilon ^{+}$ $\subset $ $U^{+}$ $\cap $ $%
U^{\prime }$ and $\Upsilon ^{-}$ $\subset $ $U^{-}$ $\cap $ $U^{\prime }$
joining $p_{0}$ to $\bar{p}$ (oriented from $p_{0}$ towards $\bar{p}).$ Let $%
\Omega _{+}$ ($\Omega _{-},$ resp.) denote the region surrounded by $%
\Upsilon $ and $\Upsilon ^{+}$ ($\Upsilon ^{-},$ resp.)$.$ From (\ref{5.6})
we have%
\begin{eqnarray}
\int_{\Omega _{+}}H &=&\lim_{\varepsilon \rightarrow 0}\doint\limits_{%
\mathcal{\partial }\Omega _{+}}N_{\varepsilon }\cdot \nu  \label{5.11.4} \\
&=&\lim_{\varepsilon \rightarrow 0}\int_{\Upsilon }N_{\varepsilon }\cdot \nu
-\int_{\Upsilon ^{+}}N_{+}\cdot \nu  \notag
\end{eqnarray}

\noindent by the Lebesgue Dominated Convergence Theorem since $%
N_{\varepsilon }$ $\rightarrow $ $N$ $=$ $N_{+}$ on $\Upsilon ^{+}$ a.e.
(pointwise convergence except 2 end points) and $|N_{\varepsilon }|$ $\leq $ 
$2,$ say, for $\varepsilon $ small. On the other hand, we deduce from (\ref%
{5.11.2}) that%
\begin{eqnarray}
\int_{\Omega _{+}}H &=&\doint\limits_{\mathcal{\partial }\Omega
_{+}}N_{+}\cdot \nu  \label{5.11.5} \\
&=&\int_{\Upsilon }N_{+}\cdot \nu -\int_{\Upsilon ^{+}}N_{+}\cdot \nu . 
\notag
\end{eqnarray}

\noindent Comparing (\ref{5.11.4}) with (\ref{5.11.5}) gives%
\begin{equation}
\int_{\Upsilon }N_{+}\cdot \nu =\lim_{\varepsilon \rightarrow
0}\int_{\Upsilon }N_{\varepsilon }\cdot \nu .  \label{5.11.6}
\end{equation}

\noindent Similarly we also have%
\begin{equation}
\int_{\Upsilon }N_{-}\cdot \nu =\lim_{\varepsilon \rightarrow
0}\int_{\Upsilon }N_{\varepsilon }\cdot \nu .  \label{5.11.7}
\end{equation}

\noindent Now (\ref{5.11.3}) follows from (\ref{5.11.6}) and (\ref{5.11.7}).

\endproof%

\bigskip

\textbf{Definition 5.2}. Let $\theta $ $\in $ $C^{0}(\Omega \backslash K)$
be a weak solution to (\ref{1.8.1}) with $H$ $\in $ $C^{0}(\Omega )$. We
call $p$ $\in $ $K$ (where $\theta $ is not defined) a crack point if

(i) there is an open neighborhood $U$ $\subset $ $\Omega $ of $p$ such that $%
U\cap K$ is a $C^{0}$ curve dividing $U$ into two connected regions $U^{+}$, 
$U^{-}$;

(ii) $N_{\pm }(p)$ exists as the limit of $N(q)$ for $q$ $\in $ $U^{\pm }$ $%
\rightarrow $ $p,$ resp. and $N_{+}(p)$ $\neq $ $N_{-}(p).$

\bigskip

\textbf{Corollary 5.3}\textit{. Suppose we are in the situation of Lemma
5.2. In particular, assume }$p_{0}$\textit{\ is nondegenerate. Moreover, we
assume that }$p_{0}$\textit{\ is a crack point}. \textit{Then we have}

\begin{equation}
\lim_{\Delta \tau ^{+}\rightarrow 0^{+}(0^{-},\text{ resp.)}}\frac{%
X_{+}(\tau _{0}^{+}+\Delta \tau ^{+})-X_{+}(\tau _{0}^{+})}{|X_{+}(\tau
_{0}^{+}+\Delta \tau ^{+})-X_{+}(\tau _{0}^{+})|}=\pm \frac{N_{+}-N_{-}}{%
|N_{+}-N_{-}|}(p_{0}).  \label{5.11.8}
\end{equation}

\begin{equation*}
(\lim_{\Delta \tau ^{-}\rightarrow 0^{+}(0^{-},\text{ resp.)}}\frac{%
X_{-}(\tau _{0}^{-}+\Delta \tau ^{-})-X_{-}(\tau _{0}^{-})}{|X_{-}(\tau
_{0}^{-}+\Delta \tau ^{-})-X_{-}(\tau _{0}^{-})|}=\pm \frac{N_{+}-N_{-}}{%
|N_{+}-N_{-}|}(p_{0}),\text{ resp.)}
\end{equation*}

\bigskip

\proof
Take $\Upsilon $ to be the line segment joining $X_{\pm }(\tau _{0}^{\pm })$
to $X_{\pm }(\tau _{0}^{\pm }+\Delta \tau ^{\pm })$ in (\ref{5.11.3}). By
the mean-value theorem we then have%
\begin{equation}
(N_{+}-N_{-})(\tilde{p})\cdot (X_{\pm }(\tau _{0}^{\pm }+\Delta \tau ^{\pm
})-X_{\pm }(\tau _{0}^{\pm }))^{\perp }=0  \label{5.11.8a}
\end{equation}

\noindent for $\tilde{p}$ $\in $ $\Upsilon .$ Taking $\Delta \tau ^{\pm }$ $%
\rightarrow $ $0$ (hence $\tilde{p}$ $\rightarrow $ $p_{0})$ we get 
\begin{equation}
\lim_{\Delta \tau ^{\pm }\rightarrow 0}\frac{X_{\pm }(\tau _{0}^{\pm
}+\Delta \tau ^{\pm })-X_{\pm }(\tau _{0}^{\pm })}{|X_{\pm }(\tau _{0}^{\pm
}+\Delta \tau ^{\pm })-X_{\pm }(\tau _{0}^{\pm })|}\cdot (N_{+}^{\perp
}-N_{-}^{\perp })(p_{0})=0.  \label{5.11.9}
\end{equation}

\noindent if the limit exists. In case $N_{+}$ $\neq $ $N_{-}$ at $p_{0},$
the limit exists in view of (\ref{5.11.8a}), and (\ref{5.11.8}) follows from
(\ref{5.11.9}).

\endproof%

\bigskip

\proof
\textbf{(of Theorem G}$^{\prime })$ From (\ref{5.11}) we have%
\begin{eqnarray}
&&\frac{X_{\pm }(\tau _{0}^{\pm }+\Delta \tau ^{\pm })-X_{\pm }(\tau
_{0}^{\pm })}{\Delta \tau ^{\pm }}  \label{5.12} \\
&=&\frac{\partial P_{\pm }}{\partial \sigma _{\pm }}(\sigma _{\pm }(\tau
_{0}^{\pm }+\Delta \tau _{1}^{\pm }),\tau _{0}^{\pm })\frac{\sigma _{\pm
}(\tau _{0}^{\pm }+\Delta \tau ^{\pm })-\sigma _{\pm }(\tau _{0}^{\pm })}{%
\Delta \tau ^{\pm }}  \notag \\
&&+\frac{\partial P_{\pm }}{\partial \tau ^{\pm }}(\sigma _{\pm }(\tau
_{0}^{\pm }+\Delta \tau ^{\pm }),\tau _{0}^{\pm }+\Delta \tau _{2}^{\pm }) 
\notag
\end{eqnarray}%
\noindent where $\Delta \tau _{1}^{\pm }$ and $\Delta \tau _{2}^{\pm }$ are
numbers between $0$ and $\Delta \tau ^{\pm }.$ Observe that $\frac{\partial
P_{\pm }}{\partial \sigma _{\pm }}$ $=$ $N_{\pm }^{\perp }$ from ($\ref%
{5.10.1}).$ Let $\Delta \sigma _{\pm }$ $=$ $\sigma _{\pm }(\tau _{0}^{\pm
}+\Delta \tau ^{\pm })-\sigma _{\pm }(\tau _{0}^{\pm }).$ Denote%
\begin{eqnarray}
&&N_{\pm }^{\perp }(P_{\pm }(\sigma _{\pm }(\tau _{0}^{\pm }+\Delta \tau
_{1}^{\pm }),\tau _{0}^{\pm }))\cdot \frac{N_{+}-N_{-}}{|N_{+}-N_{-}|}%
(p_{0}),  \label{5.12.2} \\
&&N_{\pm }^{\perp }(P_{\pm }(\sigma _{\pm }(\tau _{0}^{\pm }+\Delta \tau
_{1}^{\pm }),\tau _{0}^{\pm }))\cdot \frac{N_{+}^{\perp }-N_{-}^{\perp }}{%
|N_{+}^{\perp }-N_{-}^{\perp }|}(p_{0}),  \notag \\
&&\frac{\partial P_{\pm }}{\partial \tau ^{\pm }}(\sigma _{\pm }(\tau
_{0}^{\pm }+\Delta \tau ^{\pm }),\tau _{0}^{\pm }+\Delta \tau _{2}^{\pm
})\cdot \frac{N_{+}-N_{-}}{|N_{+}-N_{-}|}(p_{0}),  \notag \\
&&\frac{\partial P_{\pm }}{\partial \tau ^{\pm }}(\sigma _{\pm }(\tau
_{0}^{\pm }+\Delta \tau ^{\pm }),\tau _{0}^{\pm }+\Delta \tau _{2}^{\pm
})\cdot \frac{N_{+}^{\perp }-N_{-}^{\perp }}{|N_{+}^{\perp }-N_{-}^{\perp }|}%
(p_{0})  \notag
\end{eqnarray}

\noindent by $A_{1},$ $A_{2},$ $B_{1},$ $B_{2},$ resp.. So we can write (\ref%
{5.12}) as follows:%
\begin{eqnarray}
&&\frac{X_{\pm }(\tau _{0}^{\pm }+\Delta \tau ^{\pm })-X_{\pm }(\tau
_{0}^{\pm })}{\Delta \tau ^{\pm }}  \label{5.12.3} \\
&=&(A_{1}\frac{\Delta \sigma _{\pm }}{\Delta \tau ^{\pm }}+B_{1})\frac{%
N_{+}-N_{-}}{|N_{+}-N_{-}|}(p_{0})+(A_{2}\frac{\Delta \sigma _{\pm }}{\Delta
\tau ^{\pm }}+B_{2})\frac{N_{+}^{\perp }-N_{-}^{\perp }}{|N_{+}^{\perp
}-N_{-}^{\perp }|}(p_{0}).  \notag
\end{eqnarray}

Observe from (\ref{5.12.2}) that%
\begin{eqnarray}
\lim_{\Delta \tau ^{\pm }\rightarrow 0}A_{1} &=&(N_{\pm }^{\perp }\cdot 
\frac{N_{+}-N_{-}}{|N_{+}-N_{-}|})(p_{0})=\frac{N_{-}^{\perp }\cdot N_{+}}{%
|N_{+}-N_{-}|}(p_{0}),  \label{5.12.4} \\
\lim_{\Delta \tau ^{\pm }\rightarrow 0}A_{2} &=&(N_{\pm }^{\perp }\cdot 
\frac{N_{+}^{\perp }-N_{-}^{\perp }}{|N_{+}^{\perp }-N_{-}^{\perp }|}%
)(p_{0})=\pm \frac{1-N_{+}\cdot N_{-}}{|N_{+}-N_{-}|}(p_{0}),  \notag \\
\lim_{\Delta \tau ^{\pm }\rightarrow 0}B_{1} &=&\frac{\partial P_{\pm }}{%
\partial \tau ^{\pm }}(p_{0})\cdot \frac{N_{+}-N_{-}}{|N_{+}-N_{-}|}(p_{0}),
\notag \\
\lim_{\Delta \tau ^{\pm }\rightarrow 0}B_{2} &=&\frac{\partial P_{\pm }}{%
\partial \tau ^{\pm }}(p_{0})\cdot \frac{N_{+}^{\perp }-N_{-}^{\perp }}{%
|N_{+}^{\perp }-N_{-}^{\perp }|}(p_{0}).  \notag
\end{eqnarray}

\noindent Comparing (\ref{5.12.3}) with (\ref{5.11.8}), we obtain%
\begin{equation}
\lim_{\Delta \tau ^{\pm }\rightarrow 0}\frac{A_{2}\frac{\Delta \sigma _{\pm }%
}{\Delta \tau ^{\pm }}+B_{2}}{A_{1}\frac{\Delta \sigma _{\pm }}{\Delta \tau
^{\pm }}+B_{1}}=0.  \label{5.12.5}
\end{equation}

\noindent Note that $N_{+}$ $\neq $ $N_{-}$ at $p_{0}$ by the assumption
that $p_{0}$ is a crack point. It follows from (\ref{5.12.4}) that $%
\lim_{\Delta \tau ^{\pm }\rightarrow 0}A_{2}$ $\neq $ $0$. Hence from (\ref%
{5.12.5}) we conclude that the limit of $\frac{\Delta \sigma _{\pm }}{\Delta
\tau ^{\pm }}$ exists and%
\begin{equation}
\lim_{\Delta \tau ^{\pm }\rightarrow 0}\frac{\Delta \sigma _{\pm }}{\Delta
\tau ^{\pm }}=-\frac{\lim_{\Delta \tau ^{\pm }\rightarrow 0}B_{2}}{%
\lim_{\Delta \tau ^{\pm }\rightarrow 0}A_{2}}  \label{5.12.6}
\end{equation}%
\noindent since all the limits of $A_{1},$ $A_{2},$ $B_{1},$ $B_{2},$as $%
\Delta \tau ^{\pm }$ $\rightarrow $ $0$ exist by (\ref{5.12.4}). From (\ref%
{5.12.3}) we can then compute 
\begin{eqnarray}
&&\lim_{\Delta \tau ^{\pm }\rightarrow 0}\frac{X_{\pm }(\tau _{0}^{\pm
}+\Delta \tau ^{\pm })-X_{\pm }(\tau _{0}^{\pm })}{\Delta \tau ^{\pm }}
\label{5.12.7} \\
&=&[\lim_{\Delta \tau ^{\pm }\rightarrow 0}A_{1}(-\frac{\lim_{\Delta \tau
^{\pm }\rightarrow 0}B_{2}}{\lim_{\Delta \tau ^{\pm }\rightarrow 0}A_{2}}%
)+\lim_{\Delta \tau ^{\pm }\rightarrow 0}B_{1}]\frac{N_{+}-N_{-}}{%
|N_{+}-N_{-}|}(p_{0})  \notag \\
&&+[\lim_{\Delta \tau ^{\pm }\rightarrow 0}A_{2}(-\frac{\lim_{\Delta \tau
^{\pm }\rightarrow 0}B_{2}}{\lim_{\Delta \tau ^{\pm }\rightarrow 0}A_{2}}%
)+\lim_{\Delta \tau ^{\pm }\rightarrow 0}B_{2}]\frac{N_{+}^{\perp
}-N_{-}^{\perp }}{|N_{+}^{\perp }-N_{-}^{\perp }|}(p_{0})  \notag \\
&=&[\lim_{\Delta \tau ^{\pm }\rightarrow 0}A_{1}(-\frac{\lim_{\Delta \tau
^{\pm }\rightarrow 0}B_{2}}{\lim_{\Delta \tau ^{\pm }\rightarrow 0}A_{2}}%
)+\lim_{\Delta \tau ^{\pm }\rightarrow 0}B_{1}]\frac{N_{+}-N_{-}}{%
|N_{+}-N_{-}|}(p_{0})  \notag
\end{eqnarray}

\noindent by (\ref{5.12.6}). $(a)$ follows from (\ref{5.12.7}) since all the
quantities in the formula are continuous at $p_{0}.$ From (\ref{5.11.1}) and 
$\frac{\partial P_{\pm }}{\partial \sigma _{\pm }}$ $=$ $N_{\pm }^{\perp },$
we write%
\begin{equation}
\frac{\partial P_{\pm }}{\partial \tau ^{\pm }}(p_{0})=\lambda _{\pm
}(p_{0})N_{\pm }(p_{0})+\mu _{\pm }(p_{0})N_{\pm }^{\perp }(p_{0})
\label{5.12.8}
\end{equation}

\noindent for $\lambda _{\pm }(p_{0}),$ $\mu _{\pm }(p_{0})$ $\in $ $R$ and $%
\lambda _{\pm }(p_{0})$ $\neq $ $0.$ Substituting (\ref{5.12.8}) into (\ref%
{5.12.4}), we compute%
\begin{eqnarray}
&&\lim_{\Delta \tau ^{\pm }\rightarrow 0}A_{1}(-\frac{\lim_{\Delta \tau
^{\pm }\rightarrow 0}B_{2}}{\lim_{\Delta \tau ^{\pm }\rightarrow 0}A_{2}}%
)+\lim_{\Delta \tau ^{\pm }\rightarrow 0}B_{1}  \label{5.12.9} \\
&=&\pm \lambda _{\pm }(p_{0})[\frac{|N_{+}-N_{-}|}{1-N_{+}\cdot N_{-}}%
](p_{0})\neq 0  \notag
\end{eqnarray}

\noindent where we have used the identity ($1-N_{+}\cdot N_{-})^{2}$ $+$ $%
(N_{+}\cdot N_{-}^{\perp })^{2}$ $=$ $|N_{+}-N_{-}|^{2}$ (the term involving 
$\mu _{\pm }(p_{0})$ vanishes). Now $(b)$ follows from (\ref{5.12.7}) in
view of (\ref{5.12.9}).

\endproof%

\bigskip

In the above proof, we observe from (\ref{5.12.8}), (\ref{5.11.1}), and $%
\frac{\partial P_{\pm }}{\partial \sigma _{\pm }}$ $=$ $N_{\pm }^{\perp }$
that nondegeneracy of $p_{0}$ is equivalent to the condition 
\begin{equation}
\lambda _{\pm }(p_{0})=\frac{\partial P_{\pm }}{\partial \tau ^{\pm }}%
(p_{0})\cdot N_{\pm }(p_{0})\neq 0  \label{5.12.10}
\end{equation}

\noindent for both "$+"$ and "$-$" (the expanding rate of characteristic
curves$).$

\bigskip

\proof
\textbf{(of Theorem G) }Let $p_{0}$ be a nondegenerate singular point. By
Theorem B (b) the directions of $N_{+}^{\perp }$ and $N_{-}^{\perp }$ at $%
p_{0}$ must point inwards (outwards, resp.) of $U^{+}$ and $U^{-},$ resp. if 
$\func{curl}\vec{F}(p_{0})$ $>$ $0$ ($\func{curl}\vec{F}(p_{0})$ $<$ $0,$
resp.), where $U^{+}$ and $U^{-}$ are the regions in which a singular curve $%
\gamma $ passing through $p_{0}$ divides a small neighborhood $U$ of $p_{0}.$
Let $\Gamma _{+}$ $\subset $ $U^{+}$ ($\Gamma _{-}$ $\subset $ $U^{-},$
resp.) denote the characteristic curve meeting $p_{0}$ with tangent vector $%
N_{+}^{\perp }(p_{0})$ ($N_{-}^{\perp }(p_{0}),$ resp.) at $p_{0}.$ Suppose
that $N_{+}(p_{0})$ $=$ $N_{-}(p_{0})$ (and hence $N_{+}^{\perp }(p_{0})$ $=$
$N_{-}^{\perp }(p_{0})).$ By the uniqueness of the characteristic curves
(extending Theorem B (b) or Theorem B$^{\prime }$ in \cite{CHY2} to the case
that $p$ is singular by a similar argument), $\Gamma _{+}$ must coincide
with $\Gamma _{-}$ near $p_{0},$ a contradiction due to $U^{+}$ $\cap $ $%
U^{-}$ being empty. We have shown $N_{+}(p_{0})$ $\neq $ $N_{-}(p_{0}).$
I.e. $p_{0}$ is a crack point. Now $(a),$ $(b)$ follow from $(a),$ $(b)$ of
Theorem G$^{\prime },$ resp..

\endproof%

\bigskip

In the above proof of Theorem G we showed that a singular point is a crack
point in a certain situation.\ We now want to prove the converse. Let $%
\gamma $ be a $C^{1}$ smooth curve dividing a planar domain $U$ into two
connected regions $U^{+}$ and $U^{-}.$ Let $\tilde{u}$ $\in $ $%
C^{1}(U\backslash \gamma )$ $\cap $ $C^{0}(U)$ be such that $U\backslash
\gamma $ is a nonsingular domain. Moreover, suppose $\tilde{u}$ is a weak
solution to ($\ref{1.1})$\ with $\vec{F}$ $\in $ $C^{1}(U),$ $\func{curl}%
\vec{F}$ $\neq $ $0$ and $H$\ $\in $\ $C^{1}(U)$ in the sense that for any $%
\varphi \in C_{0}^{\infty }(U),$ there holds%
\begin{equation}
\int_{U\backslash \gamma }\frac{\nabla \tilde{u}+\vec{F}}{|\nabla \tilde{u}+%
\vec{F}|}\cdot \nabla \varphi +\int_{U}H\varphi =0.  \label{5.12.11}
\end{equation}

\noindent Take $p_{0}$ $\in $ $\gamma .$ Suppose that near $p_{0}$ we are in
the situation of Definition 5.1. Namely, we have the seed curves $\gamma
_{\pm }$ parametrized by $\tau ^{\pm }$ and the characteristic curves
issuing from $\gamma _{\pm }$ hit $\gamma .$ By adding the $u$-variable to (%
\ref{5.10.1}), we consider 
\begin{eqnarray}
\frac{dx}{d\sigma } &=&\sin \theta ,\text{ }\frac{dy}{d\sigma }=-\cos \theta
\label{5.12.12} \\
\frac{d\theta }{d\sigma } &=&-H,\text{ }\frac{du}{d\sigma }=-F_{1}\sin
\theta +F_{2}\cos \theta  \notag
\end{eqnarray}

\noindent where $\vec{F}$ $=$ $(F_{1},$ $F_{2}).$ Let ($x_{\pm }(\sigma
,\tau ^{\pm })$, $y_{\pm }(\sigma ,\tau ^{\pm }),$ $\theta _{\pm }(\sigma
,\tau ^{\pm })$, $u_{\pm }(\sigma ,\tau ^{\pm }))$ be the solution to (\ref%
{5.12.12}) with the initial data $(x_{\pm }(0,\tau ^{\pm }),$ $y_{\pm
}(0,\tau ^{\pm }))$ $=$ $\gamma _{\pm }(\tau ^{\pm })$, $\theta _{\pm
}(0,\tau ^{\pm })$ $=$ $\tilde{\theta}(\gamma _{\pm }(\tau ^{\pm }))$, and $%
u_{\pm }(0,\tau ^{\pm })$ $=$ $u(\gamma _{\pm }(\tau ^{\pm })).$ Here we
write $\frac{\nabla \tilde{u}+\vec{F}}{|\nabla \tilde{u}+\vec{F}|}$ $=$ $%
(\cos \tilde{\theta},$ $\sin \tilde{\theta})$ near $p_{0}.$ Suppose that $%
p_{0}$ is nondegenerate. Then there is a diffeomorphism between $(\sigma
,\tau ^{\pm })$ and $(x,y)$ near $p_{0}$ as shown before. By the basic ODE
theory, the solution ($x_{\pm }(\sigma ,\tau ^{\pm })$, $y_{\pm }(\sigma
,\tau ^{\pm }),$ $\theta _{\pm }(\sigma ,\tau ^{\pm })$, $u_{\pm }(\sigma
,\tau ^{\pm }))$ is $C^{1}$ smooth in ($\sigma $, $\tau ^{\pm })$ and is
defined near $p_{0}$ on the $xy$-plane. Recall that $N_{\pm }$ is defined to
be $(\cos \theta _{\pm },$ $\sin \theta _{\pm })$ and we call $p_{0}$ a
crack point if $N_{+}(p_{0})$ $\neq $ $N_{-}(p_{0}).$ Note that $dz$ $+$ $%
F_{1}dx$ $+$ $F_{2}dy$ is a contact form in $\boldsymbol{H}_{1}$ due to the
condition $\func{curl}\vec{F}$ $\neq $ $0.$ From Theorem A in \cite{CHY2}
and $du$ $+$ $F_{1}dx$ $+$ $F_{2}dy$ $=$ $0$ along the characteristic
curves, we have $\theta _{\pm }$ $=$ $\tilde{\theta},$ $u_{\pm }$ $=$ $%
\tilde{u}$ on $U^{\pm }$ by the uniqueness of solutions to the ODE system (%
\ref{5.12.12}) with the same initial data. It follows that $N_{\pm }$ $=$ $%
\frac{\nabla u_{\pm }+\vec{F}}{D_{\pm }}$ where $D_{\pm }$ $:=$ $|\nabla
u_{\pm }+\vec{F}|.$

\bigskip

\textbf{Theorem 5.4.} \textit{Let }$\gamma $\textit{\ be a }$C^{1}$\textit{\
smooth curve dividing a planar domain }$U$\textit{\ into two connected
regions }$U^{+}$\textit{\ and }$U^{-}.$\textit{\ Let }$\tilde{u}$\textit{\ }$%
\in $\textit{\ }$C^{1}(U\backslash \gamma )$\textit{\ }$\cap $\textit{\ }$%
C^{0}(U)$\textit{\ be such that }$U\backslash \gamma $\textit{\ is a
nonsingular domain. Moreover, }$\tilde{u}$\textit{\ is a weak solution to (}$%
\ref{1.1})$\textit{\ with }$\vec{F}$\textit{\ }$\in $\textit{\ }$C^{1}(U),$%
\textit{\ }$\func{curl}\vec{F}$\textit{\ }$\neq $\textit{\ }$0$\textit{\ and 
}$H$\textit{\ }$\in $\textit{\ }$C^{1}(U)$.\textit{\ Let }$p_{0}$\textit{\ }$%
\in $ $\gamma $ \textit{be a nondegenerate crack point. Then }$\tilde{u}$%
\textit{\ }$\in $\textit{\ }$C^{1}(V)$\textit{\ for some neighborhood }$V$ $%
\subset $ $U$ of $p_{0}$ \textit{and }$p_{0}$\textit{\ is a singular point
of }$\tilde{u}$\textit{. That is, }$\nabla \tilde{u}+\vec{F}$\textit{\ }$=$%
\textit{\ }$0$\textit{\ at }$p_{0}.$

\bigskip

\proof
Observe that along a characteristic curve, there holds $\frac{d\tilde{u}%
(x(\sigma ),y(\sigma ))}{d\sigma }$ $+$ $F_{1}\frac{dx}{d\sigma }$ $+$ $F_{2}%
\frac{dy}{d\sigma }$ $=$ $0$ on $U\backslash \gamma $ $=$ $U^{+}$ $\cup $ $%
U^{-}$ for $\tilde{u}$ $\in $ $C^{1}(U\backslash \gamma ).$ It follows from
the uniqueness of solutions to the ODE system (\ref{5.12.12}) with the same
initial data and $\tilde{u}$ $\in $ $C^{0}(U)$ that

\begin{eqnarray}
u_{+} &=&\tilde{u}\text{ }on\text{ }U^{+}\cup \gamma  \label{5.12.12a} \\
u_{-} &=&\tilde{u}\text{ }on\text{ }U^{-}\cup \gamma .  \notag
\end{eqnarray}%
\noindent So $u_{+}$ $=$ $u_{-}$ on $\gamma .$ Take a unit-speed parameter $%
\varsigma $ for $\gamma $ ($\in $ $C^{1})$ and note that both $u_{+}$ and $%
u_{-}$ are defined and $C^{1}$ smooth in a small neighborhood $V$ of $p_{0}$%
. We can choose $V$ such that each $p$ $\in $ $\gamma $ $\cap $ $V$ is a
nondegenerate crack point. Therefore $\frac{du_{\pm }}{d\varsigma }$ exists
and 
\begin{equation}
\frac{du_{+}}{d\varsigma }=\frac{du_{-}}{d\varsigma }  \label{5.12.13}
\end{equation}

\noindent on $\gamma $ $\cap $ $V.$ Compute%
\begin{eqnarray}
&&\frac{du_{\pm }}{d\varsigma }+F_{1}\frac{dx}{d\varsigma }+F_{2}\frac{dy}{%
d\varsigma }  \notag \\
&=&(\nabla u_{\pm }+\vec{F})\cdot \frac{d\gamma }{d\varsigma }
\label{5.12.14}
\end{eqnarray}

\noindent where $\frac{d\gamma }{d\varsigma }$ $=$ $(\frac{dx}{d\varsigma },$
$\frac{dy}{d\varsigma }).$ From (\ref{5.12.13}) and (\ref{5.12.14}) we obtain%
\begin{equation}
D_{+}N_{+}\cdot \frac{d\gamma }{d\varsigma }=D_{-}N_{-}\cdot \frac{d\gamma }{%
d\varsigma }.  \label{5.12.15}
\end{equation}

\noindent Here we recall that $D_{\pm }$ $:=$ $|\nabla u_{\pm }+\vec{F}|$
and $N_{\pm }$ $=$ $\frac{\nabla u_{\pm }+\vec{F}}{D_{\pm }}.$ From Theorem G%
$^{\prime }$ (b) (equal-angle condition) (see (\ref{5.12.7}) and (\ref%
{5.12.9})) and $N_{+}(p)$ $\neq $ $N_{-}(p)$ for each $p$ $\in $ $\gamma $ $%
\cap $ $V$, we obtain 
\begin{equation}
-N_{+}+N_{-}=(N_{+}^{\perp }-N_{-}^{\perp })^{\perp }\parallel \frac{d\gamma 
}{d\varsigma }.  \label{5.12.15a}
\end{equation}%
\noindent It follows from the identity ($N_{+}+N_{-})\cdot (N_{+}-N_{-})$ $=$
$0$ and (\ref{5.12.15a}) that 
\begin{equation*}
(N_{+}+N_{-})\cdot \frac{d\gamma }{d\varsigma }=0
\end{equation*}%
\noindent and hence%
\begin{equation}
N_{+}\cdot \frac{d\gamma }{d\varsigma }\text{ \textit{and} }N_{-}\cdot \frac{%
d\gamma }{d\varsigma }\text{ \textit{have different sign}}  \label{5.12.16}
\end{equation}

\noindent at each $p$ $\in $ $\gamma $ $\cap $ $V.$ Therefore we have $%
D_{\pm }$ $=$ $0$ on $\gamma $ $\cap $ $V$ in view of $D_{\pm }$ $\geq $ $0,$
(\ref{5.12.15}) and (\ref{5.12.16}). It follows that $\nabla u_{\pm }$ $+$ $%
\vec{F}$ $=$ $0$, and hence $\nabla u_{+}$ $=$ $\nabla u_{-}$ at each $p$ $%
\in $ $\gamma $ $\cap $ $V.$ This implies that $\nabla \tilde{u}$ exists and
equals $\nabla u_{+}$ $=$ $\nabla u_{-}$ at each $p$ $\in $ $\gamma $ $\cap $
$V$ by (\ref{5.12.12a})$.$ Since $u_{+}$, $u_{-}$ $\in $ $C^{1}(V)$, we see
that $\nabla \tilde{u}$ is continuous at each $p$ $\in $ $\gamma $ $\cap $ $%
V,$ and hence $\tilde{u}$ $\in $ $C^{1}(V).$ From $\nabla \tilde{u}$ $+$ $%
\vec{F}$ $=$ $\nabla u_{+}$ $+$ $\vec{F}$ $=$ $0$ on $\gamma $ $\cap $ $V,$
we learn that in particular, $p_{0}$ is a singular point of $\tilde{u}.$

\endproof%

\bigskip

To illustrate Theorem G, we consider the case of a $p$-minimal graph over a
planar domain $\Omega ,$ defined by $u$ $\in $ $C^{1}(\Omega )$. In this
case we have some interesting formulas such as (\ref{5.18.1}) and (\ref{5.26}%
). From (\ref{1.5'}) we can deduce that if $D^{\prime }$ $>$ $1$ at an
initial point, then%
\begin{equation}
D^{\prime }=\frac{4+cD^{2}\pm \sqrt{(4+cD^{2})cD^{2}}}{2}.  \label{5.1}
\end{equation}

\noindent Meanwhile, if $D^{\prime }$ $<$ $1$ at an initial point, we obtain%
\begin{equation}
D^{\prime }=\frac{4-cD^{2}\pm \sqrt{(4-cD^{2})(-cD^{2})}}{2}.  \label{5.2}
\end{equation}

\noindent Since $c$ $>$ $0,$ we must have $4-cD^{2}$ $\leq $ $0$ in this
case. It follows that $\sqrt{\frac{4}{c}}$ $\leq $ $D.$ So $D$ never reaches 
$0$ along a characteristic line if $D^{\prime }$ $<$ $1$ at an initial
point. Observe that%
\begin{equation}
\frac{2}{4+cD^{2}\pm \sqrt{(4+cD^{2})cD^{2}}}=\frac{1}{2}\mp \frac{1}{2}%
\frac{\sqrt{c}D}{\sqrt{4+cD^{2}}}.  \label{5.3}
\end{equation}

\noindent When $H$ $=$ $0,$ we can take $f$ $\equiv $ $1$ so that the
parameter $s$ in Theorem C of \cite{CHY2} is just the arc length $\sigma $
along characteristic lines. From (\ref{5.1}) and (\ref{5.3}) we can
integrate and obtain%
\begin{equation}
\frac{1}{2}(D\mp \sqrt{D^{2}+\frac{4}{c}})\mid _{\sigma _{0}}^{\sigma
_{1}}=\sigma _{1}-\sigma _{0}.  \label{5.4}
\end{equation}

Suppose $D(0)$ $=$ $0.$ Taking $\sigma _{0}$ $=$ $0,$ $\sigma _{1}$ $=$ $s$
in (\ref{5.4}) we obtain%
\begin{equation}
D(s)=s+a-\frac{a^{2}}{s+a}  \label{5.4.1}
\end{equation}%
\noindent where $a$ $=$ $\mp \sqrt{\frac{1}{c}}.$ By (\ref{5.4.1}) we compute%
\begin{equation}
D^{\prime }(s)=1+\frac{a^{2}}{(s+a)^{2}}.  \label{5.4.2}
\end{equation}

\noindent It follows from (\ref{5.4.2}) that $D^{\prime }$ $\rightarrow $ $2$
as $s$ $\rightarrow $ $0$ for $a$ $\neq $ $0$ while $D^{\prime }$ $\equiv $ $%
1$ for $a$ $=$ $0.$ This verifies Theorem B (a) for the case of $p$-minimal
graphs.

If we start with points on a $C^{1}$ smooth curve $\beta (\tau )$ transverse
to the characteristic lines, described by $\sigma _{0}$ $=$ $0,$ say, then
we can describe the first (singular) points where the characteristic lines
hit the singular set by $\sigma _{1}$ $=$ $\sigma _{1}(\tau ).$ Noting that $%
D$ $=$ $0$ at the points described by ($\tau ,$ $\sigma _{1}(\tau )),$ we
have%
\begin{equation}
\sigma _{1}(\tau )=\mp \frac{1}{\sqrt{c(\tau )}}-\frac{1}{2}D(\beta (\tau
))\pm \frac{1}{2}\sqrt{D(\beta (\tau ))^{2}+\frac{4}{c(\tau )}}  \label{5.5}
\end{equation}

\noindent by (\ref{5.4}). Here we have written $c$ $=$ $c(\tau ).$ Since $D$
and $D^{\prime }$ are continuous by Theorem D of \cite{CHY2}, we conclude
that $c$ $=$ $c(\tau )$ is continuous in $\tau $ by (\ref{1.5'}). It follows
from (\ref{5.5}) that $\sigma _{1}(\tau )$ is continuous in $\tau .$ For
general $H$ $\in $ $C^{1}$ and $\vec{F}$ $\in $ $C^{1},$ we cannot expect to
have an explicit formula for $\sigma _{1}(\tau )$ (replace characteristic
lines by characteristic curves in the definition) like (\ref{5.5}). But
still $\sigma _{1}(\tau )$ is $C^{0}$ in $\tau $ by Lemma 3.2$^{\prime }.$

Let $p$ $\in $ $S(u)$ be a singular point. Suppose we are in the situation
described in Definition 5.1 with $K$ replaced by $S(u)$. So we have $\gamma
_{\pm },$ parametrized by $\tau ^{\pm },$ $X_{\pm }(\tau ^{\pm }),$ and $%
X_{+}(\tau _{0}^{+})$ $=$ $X_{\_}(\tau _{0}^{-})$ $=$ $p_{0}.$ Since the
characteristic curves are straight lines in this case ($\theta $ being
constant along a characteristic curve due to $\frac{d\theta }{d\sigma }$ $=$ 
$-H$ $=$ $0$ by (\ref{5.10.1})), we can write%
\begin{eqnarray}
X_{+}(\tau ^{+}) &=&\gamma _{+}(\tau ^{+})+\sigma _{+}(\tau
^{+})N_{+}^{\perp }(\tau ^{+})  \label{5.12.17} \\
X_{\_}(\tau ^{-}) &=&\gamma _{-}(\tau ^{-})+\sigma _{-}(\tau
^{-})N_{-}^{\perp }(\tau ^{-})  \notag
\end{eqnarray}%
\noindent for some real functions $\sigma _{\pm }$ ($\in C^{0}$ since $%
X_{\pm },$ $\gamma _{\pm },$ and $N_{\pm }^{\perp }$ are $C^{0}),$ where $%
N_{\pm }^{\perp }(\tau ^{\pm })$ :$=$ $(\sin \theta _{\pm }(\tau ^{\pm }),$ $%
-\cos \theta _{\pm }(\tau ^{\pm }))$ for some angular functions $\theta
_{\pm }$ (cf. (\ref{5.11}))$.$ It follows from (\ref{5.12.17}) and (\ref%
{5.11}) that%
\begin{eqnarray}
\frac{\partial P_{\pm }}{\partial \tau ^{\pm }} &=&\frac{\partial \gamma
_{\pm }}{\partial \tau ^{\pm }}+\sigma _{\pm }\frac{\partial N_{\pm }^{\perp
}}{\partial \tau ^{\pm }}  \label{5.12.18} \\
&=&N_{\pm }+\sigma _{\pm }\theta _{\pm }^{\prime }N_{\pm }  \notag
\end{eqnarray}%
\noindent (noting that $\theta _{\pm }$ $\in $ $C^{1}$ by the ODE theory and 
$\theta $ $\in $ $C^{1}$ by Theorem D in \cite{CHY2}). From (\ref{5.12.10})
and (\ref{5.12.18}) we have%
\begin{eqnarray}
\lambda _{\pm }(p_{0}) &=&\frac{\partial P_{\pm }}{\partial \tau ^{\pm }}%
(p_{0})\cdot N_{\pm }(p_{0})  \label{5.18.1} \\
&=&1+\sigma _{\pm }(\tau _{0}^{\pm })\theta _{\pm }^{\prime }(\tau _{0}^{\pm
}).  \notag
\end{eqnarray}

Recall that $p_{0}$ is nondegenerate if and only if both $\lambda
_{+}(p_{0}) $ $\neq $ $0$ and $\lambda _{-}(p_{0})$ $\neq $ $0.$ Since $%
\frac{\partial P_{\pm }}{\partial \sigma _{\pm }}$ $=$ $N_{\pm }^{\perp }$
and hence $\frac{\partial P_{\pm }}{\partial \tau ^{\pm }}\cdot N_{\pm }$ $=$
$\frac{\partial P_{\pm }}{\partial \sigma _{\pm }}\wedge \frac{\partial
P_{\pm }}{\partial \tau ^{\pm }},$ we learn that $\lambda _{\pm }(p_{0})$ $%
\neq $ $0$ is independent of the choice of seed curves by (\ref{5.11.1a}).

\bigskip

\textbf{Proposition 5.5.} \textit{Consider a }$p$\textit{-minimal graph over
a planar domain }$\Omega ,$\textit{\ defined by }$u$\textit{\ }$\in $\textit{%
\ }$C^{1}(\Omega )$\textit{. Suppose we are in the situation described in
Definition 5.1 with }$K$\textit{\ replaced by }$S(u)$\textit{. Suppose
further }$X_{+}(\tau ^{+})$ $=$ $X_{-}(\tau ^{-}).$ \textit{Then we have}

\begin{equation}
\sigma _{+}(\tau ^{+})-\sigma _{+}(\tau _{0}^{+})=\sigma _{-}(\tau
^{-})-\sigma _{-}(\tau _{0}^{-}).  \label{5.26}
\end{equation}

\bigskip

\proof
Let $\Gamma _{\tau ^{+}}^{+}$ ($\Gamma _{\tau _{0}^{+}}^{+},$ $\Gamma _{\tau
^{-}}^{-},$ $\Gamma _{\tau _{0}^{-}}^{-},$ resp.) denote the characteristic
line (segment) which connect $\gamma _{+}(\tau ^{+})$ ($\gamma _{+}(\tau
_{0}^{+}),$ $\gamma _{-}(\tau ^{-}),$ $\gamma _{-}(\tau _{0}^{-}),$ resp.)
with $X_{+}(\tau ^{+})$ ($X_{+}(\tau _{0}^{+}),$ $X_{-}(\tau ^{-}),$ $%
X_{-}(\tau _{0}^{-}),$ resp.) (see Figure 5.1).

\bigskip

\begin{figure}[th]
\begin{center}
\includegraphics[width=6cm]{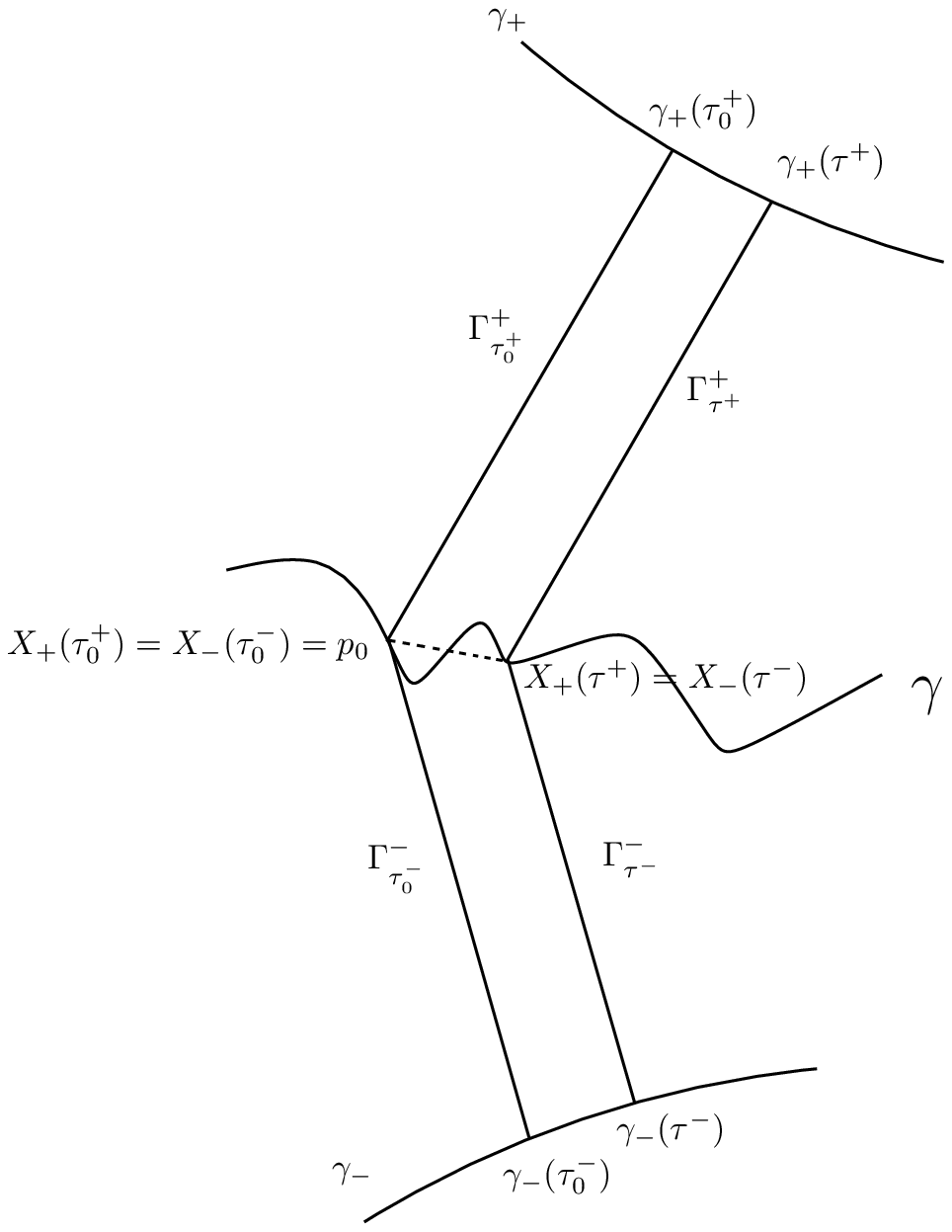}\\[0pt]
Figure 5.1
\end{center}
\par
\end{figure}

\bigskip

Let $\hat{\Omega}$ denote the region surrounded by $\Gamma _{\tau
_{0}^{+}}^{+},$ $\Gamma _{\tau _{0}^{-}}^{-},$ segment($\gamma _{-}(\tau
_{0}^{-})\gamma _{-}(\tau ^{-})),$ $\Gamma _{\tau ^{-}}^{-},$ $\Gamma _{\tau
^{+}}^{+}$, and segment($\gamma _{+}(\tau ^{+})\gamma _{+}(\tau _{0}^{+}))$.
By Theorem B (b) (the unit outward normal) $\nu $ $=$ $+N$ along $\Gamma
_{\tau _{0}^{+}}^{+}$ and $\Gamma _{\tau ^{-}}^{-}$ while $\nu $ $=$ $-N$
along $\Gamma _{\tau _{0}^{-}}^{-}$ and $\Gamma _{\tau ^{+}}^{+}.$ We can
now apply Lemma 5.1 with $\vec{F}$ $=$ $(-y,x)$ and $H$ $=$ $0$ to obtain%
\begin{eqnarray}
0 &=&\doint\limits_{\mathcal{\partial }\hat{\Omega}}N\cdot \nu =\int_{\Gamma
_{\tau _{0}^{+}}^{+}}N\cdot N  \label{5.28} \\
&&+\int_{\Gamma _{\tau _{0}^{-}}^{-}}N\cdot (-N)+\int_{\Gamma _{\tau
^{-}}^{-}}N\cdot N+\int_{\Gamma _{\tau ^{+}}^{+}}N\cdot (-N)  \notag \\
&=&\sigma _{+}(\tau _{0}^{+})-\sigma _{-}(\tau _{0}^{-})+\sigma _{-}(\tau
^{-})-\sigma _{+}(\tau ^{+}).  \notag
\end{eqnarray}

\noindent Note that $\nu $ $=$ $\pm N^{\perp }$ and hence $N\cdot \nu $ $=$ $%
0$ along the seed curves $\gamma _{+}$ and $\gamma _{-}.$ Now (\ref{5.26})
follows from (\ref{5.28}).

\endproof%

\bigskip

We can interpret the quantity $\lambda _{\pm }(p_{0})$ $=$ $1$ $+$ $\sigma
_{\pm }(\tau _{0}^{\pm })\theta _{\pm }^{\prime }(\tau _{0}^{\pm })$ in (\ref%
{5.18.1}) for a nonsingular point $p_{0}$ in terms of an integrating factor
for solving the $t$-coordinate in Theorem C of \cite{CHY2}. Recall that the
integrating factor $gD$ for solving $t$ such that $\nabla t$ $=$ $gDN$
satisfies%
\begin{equation}
N^{\perp }g+\frac{(curl\vec{F})g}{D}=0.  \label{5.40}
\end{equation}

\noindent Let $\tilde{g}$ $:=$ $gD.$ A direct computation shows that 
\begin{equation}
N^{\perp }\tilde{g}=\tilde{g}(\frac{N^{\perp }D-curl\vec{F}}{D}).
\label{5.41}
\end{equation}

\noindent By (1.13) and $N^{\perp }$ $=$ $f\frac{\partial }{\partial s}$ in 
\cite{CHY2} we can reduce (\ref{5.41}) to%
\begin{equation}
f\frac{\partial \tilde{g}}{\partial s}+\tilde{g}^{2}\frac{\partial \theta }{%
\partial t}=0.  \label{5.42}
\end{equation}

\noindent We can take $f$ $\equiv $ $1$ (satisfying the first equation of
(1.10) in \cite{CHY2}) for the case $H$ $=$ $0.$ Observe that $\theta $ $\in 
$ $C^{1}$ by Theorem D in \cite{CHY2} and $\partial _{t}(\partial _{s}\theta
)$ exists and is continuous by (1.12) of \cite{CHY2} for $H$ $\in $ $C^{1},$
say ($f$ is known to be $C^{1}$ smooth in $t).$ By Lemma 5.4 in \cite{CHY2}
we have the existence of $\partial _{s}(\partial _{t}\theta )$ and%
\begin{eqnarray}
\partial _{s}(\partial _{t}\theta ) &=&\partial _{t}(\partial _{s}\theta
)=\partial _{t}(-\frac{H}{f})  \label{5.43} \\
( &=&0\text{ \textit{if} }H=0).  \notag
\end{eqnarray}

\noindent So $\partial _{t}\theta $ is independent of $s$ for the case $H$ $%
= $ $0$ by (\ref{5.43})$.$ We can now solve $\tilde{g}$ for (\ref{5.42}) in
the case of $H$ $=$ $0$ (in which $f$ $\equiv $ $1)$ to get%
\begin{equation*}
\tilde{g}^{-1}(s,t)=1+s\frac{\partial \theta }{\partial t}(0,t).
\end{equation*}

\noindent Here we have taken $\tilde{g}(0,t)$ $=$ $1$ for which $\partial
_{t}$ $=$ $\partial _{\tau }$ at $s$ $=$ $0.$

\bigskip

\section{The local theory of surfaces with prescribed $p$-mean curvature}

We are going to prove Theorem H. We have in mind that $V$ plays the role of $%
N^{\perp }.$ We need to find $N$ in the $\xi ,$ $\eta $ coordinates. Recall
that on the $xy$-plane, $\func{div}DN^{\perp }$ $=$ $\func{div}$ $(u_{y}+x,$ 
$-(u_{x}-y))$ $=$ $2$ for $N$ $=$ $\frac{(u_{x}-y,u_{y}+x)}{D}$ and $D$ $=$ [%
$(u_{x}-y)^{2}$ $+$ $(u_{y}+x)^{2}]^{1/2}.$ It follows that $N^{\perp }(D)$ $%
=$ $2-N(\theta )D$ if we write $N$ $=$ $(\cos \theta ,$ $\sin \theta )$
while $N^{\perp }$ $=$ $(\sin \theta ,$ $-\cos \theta ).$ Therefore we have%
\begin{equation}
N(\theta )=\frac{2-N^{\perp }(D)}{D}.  \label{6.1}
\end{equation}

We compute the commutator of $N$ and $N^{\perp }:$%
\begin{eqnarray}
&&[N,N^{\perp }]  \label{6.2} \\
&=&(\theta _{x}\cos ^{2}\theta +\theta _{y}\sin \theta \cos \theta )\partial
_{x}+(\theta _{x}\cos \theta \sin \theta +\theta _{y}\sin ^{2}\theta
)\partial _{y}  \notag \\
&&-(\theta _{x}\sin \theta (-\sin \theta )+\theta _{y}\cos \theta \sin
\theta )\partial _{x}-(\theta _{x}\cos \theta \sin \theta -\theta _{y}\cos
^{2}\theta )\partial _{y}  \notag \\
&=&\theta _{x}\partial _{x}+\theta _{y}\partial _{y}=\nabla \theta .  \notag
\end{eqnarray}

\noindent We can express%
\begin{eqnarray}
\nabla \theta &=&N(\theta )N+N^{\perp }(\theta )N^{\perp }  \label{6.3} \\
&=&(\frac{2-N^{\perp }(D)}{D})N-HN^{\perp }  \notag
\end{eqnarray}

\noindent by (\ref{6.1}) and (2.23) in \cite{CHMY04}.

\bigskip

\proof
\textbf{(of Theorem H) }Without\textbf{\ }specifying the regularity, we mean 
$C^{\infty }$ smoothness for each quantity in the following argument. Note
that $L_{V}P$ $=$ $[V,P]$ and we (having $P$ $=$ $N$ in mind) obtain
equation (\ref{1.9.1}) in view of (\ref{6.2}) and (\ref{6.3}). Next we can
find a solution $f$ $>$ $0$ ($g$ $>$ $0,$ resp.) to the equation 
\begin{equation}
P(f)+Hf=0\text{ }(V(g)+\frac{2g}{D}=0,\text{ resp}.)  \label{6.5}
\end{equation}%
\noindent in a small neighborhood of $p$ in the $\xi \eta $-plane by
assigning any positive initial value of $f$ ($g,$ resp.) along an integral
curve of $V$ ($P,$ resp.) through $p.$ Now we want to solve in $s$ and $t$
for the following equations:%
\begin{equation}
V(s)=f,\text{ }P(s)=0;  \label{6.6.1}
\end{equation}

\begin{equation}
V(t)=0,\text{ }P(t)=gD.  \label{6.6.2}
\end{equation}

\noindent From (\ref{1.9.1}) we need to check whether the integrability
condition for (\ref{6.6.1}) ((\ref{6.6.2}), resp.) $-(\frac{2-V(D)}{D})P(s)$ 
$+$ $HV(s)$ $=$ $VP(s)$ $-$ $PV(s)$ $=$ $0$ $-$ $P(f)$ ($-(\frac{2-V(D)}{D}%
)P(t)$ $+$ $HV(t)$ $=$ $VP(t)$ $-$ $PV(t)$ $=$ $V(gD),$ resp.) holds by
Frobenius' theorem. A direct computation shows that (\ref{6.5}) makes these
integrability conditions hold. Therefore (\ref{6.6.1}) and (\ref{6.6.2}) are
solvable. Since $V$ and $P$ are transversal, $s$ and $t$ form a local
coordinate system by (\ref{6.6.1}) and (\ref{6.6.2}) (note that $f,$ $g,$
and $D$ are all positive). Moreover, $(V,$ $P)$ has the same orientation as $%
(\partial _{s},$ $\partial _{t}).$

Next we want to find local coordinate functions $x$ and $y$ such that 
\begin{equation}
\frac{ds^{2}}{f^{2}}+\frac{dt^{2}}{g^{2}D^{2}}=dx^{2}+dy^{2}  \label{6.7}
\end{equation}

\noindent (in view of (1.11) in \cite{CHY2}). By the fundamental theorem of
Riemannian geometry, we need to check if the Gaussian curvature $K$ of the
metric $\frac{ds^{2}}{f^{2}}$ $+$ $\frac{dt^{2}}{g^{2}D^{2}}$ equals zero.
For a metric of the form $Eds^{2}$ $+$ $Gdt^{2}$ (orthogonal parametrization)%
$,$ we have%
\begin{equation}
K=-\frac{1}{2A}[\partial _{s}(\frac{G_{s}}{A})+\partial _{t}(\frac{E_{t}}{A}%
)]  \label{6.8}
\end{equation}

\noindent where $A$ $=$ $\sqrt{EG}.$ Substituting $E$ $=$ $\frac{1}{f^{2}}$
and $G$ $=$ $\frac{1}{g^{2}D^{2}}$ into (\ref{6.8}) we have $A$ $=$ $\frac{1%
}{fgD}$ and%
\begin{eqnarray}
K &=&-\frac{fgD}{2}[\partial _{s}(\frac{(\frac{1}{g^{2}D^{2}})_{s}}{\frac{1}{%
fgD}})+\partial _{t}(\frac{(\frac{1}{f^{2}})_{t}}{\frac{1}{fgD}})]
\label{6.9} \\
&=&-\frac{fgD}{2}[\partial
_{s}(-2g^{-2}D^{-1}fg_{s}-2g^{-1}D^{-2}fD_{s})+\partial
_{t}(-2f^{-2}gDf_{t})].  \notag
\end{eqnarray}

\noindent From (\ref{6.6.1}) and (\ref{6.6.2}) we can easily relate $%
\partial _{s},$ $\partial _{t}$ to $V,$ $P$ as follows:%
\begin{equation}
V=f\frac{\partial }{\partial s},\text{ }P=gD\frac{\partial }{\partial t}.
\label{6.10}
\end{equation}

\noindent It follows from (\ref{6.10}) and (\ref{6.5}) that 
\begin{eqnarray}
fg_{s} &=&Vg=-\frac{2g}{D},\text{ }fD_{s}=V(D),\text{ and}  \label{6.11} \\
gDf_{t} &=&Pf=-Hf.  \notag
\end{eqnarray}

\noindent Substituting (\ref{6.11}) into (\ref{6.9}), we obtain%
\begin{eqnarray}
K &=&-\frac{fgD}{2}[\partial _{s}(4g^{-1}D^{-2}-2g^{-1}D^{-2}D^{\prime
})+\partial _{t}(2Hf^{-1})]  \label{6.12} \\
&=&-\frac{gD}{2}\{[4(-1)g^{-2}V(g)D^{-2}+4g^{-1}(-2)D^{-3}D^{\prime }](1-%
\frac{1}{2}D^{\prime })  \notag \\
&&+4g^{-1}D^{-2}(-\frac{1}{2}D^{\prime \prime })\}-\frac{f}{2}%
(2H(-1)f^{-2}N(f))-P(H)  \notag \\
&=&D^{-1}D^{\prime \prime }-2D^{-2}(D^{\prime }-1)(D^{\prime }-2)-H^{2}-P(H)
\notag
\end{eqnarray}

\noindent by (\ref{6.10}) and (\ref{6.11}). Comparing (\ref{6.12}) with the
condition (\ref{1.9}), we can finally conclude that $K$ $=$ $0.$ So we have
proved the existence of local coordinates $x$ and $y$ such that (\ref{6.7})
holds. Moreover, we can find $x,$ $y$ such that $(\partial _{x},$ $\partial
_{y})$ has the same orientation as $(\partial _{\xi },$ $\partial _{\eta }).$

Observe that both $V$ and $P$ are unit vectors and orthogonal with respect
to the metric (\ref{6.7}). So we can write 
\begin{equation}
V(x)=\sin \theta ,\text{ }V(y)=-\cos \theta  \label{6.13}
\end{equation}

\noindent for some function $\theta $ locally near $p.$ It follows from the
orthonormality and the arrangement of orientation that%
\begin{equation}
P(x)=\cos \theta ,\text{ }P(y)=\sin \theta .  \label{6.14}
\end{equation}

\noindent So we have $N$ $=$ $P$ and $N^{\perp }$ $=$ $V.$ With $N^{\perp }$
replaced by $V$ in (\ref{6.2}) and (\ref{6.3})$,$ we get 
\begin{equation}
\lbrack V,N]=-N(\theta )N-V(\theta )V  \label{6.15}
\end{equation}

\noindent by (\ref{6.13}) and (\ref{6.14}). Comparing (\ref{6.15}) with (\ref%
{1.9.1}) ($P$ replaced by $N)$, we obtain%
\begin{eqnarray}
N(\theta ) &=&\frac{2-V(D)}{D}\text{ and}  \label{6.15.1} \\
V(\theta ) &=&-H.  \notag
\end{eqnarray}

\noindent We have proved (\ref{1.11}) and (\ref{1.10}).

Next we are going to prove (2). Take a local integral curve $\ell $ of $N$
through $p$ (may assume $x(p)$ $=$ $y(p)$ $=$ $0).$ Choose a $C^{\infty }$
smooth function $u_{0}$ along $\ell $ with $N(u_{0})+(-y,x)\cdot N=D$ where "%
$\cdot "$ denotes the standard planar inner product$.$ For any point $q$ $%
\in $ $\ell ,$ there passes an integral (characteristic) curve $\Gamma _{q}$
of $V.$ Define the value of $u$ on $\Gamma _{q}$ by integrating the contact
form 
\begin{equation}
du+xdy-ydx=0  \label{6.16}
\end{equation}%
\noindent along $\Gamma _{q}.$ That is, at $\zeta $ $\in $ $\Gamma _{q},$ we
define%
\begin{equation*}
u(\zeta )=u_{0}(q)+\int_{q}^{\zeta }(ydx-xdy)
\end{equation*}

\noindent where the integral means the line integral from $q$ to $\zeta $
along $\Gamma _{q}.$ Thus $u$ $=$ $u(x,y)$ is defined in an open
neighborhood of $p$ and is a $C^{\infty }$ smooth function. Writing (\ref%
{6.16}) as ($u_{x}-y)dx$ $+$ $(u_{y}+x)dy$ $=$ $0,$ we obtain%
\begin{equation*}
\lbrack \nabla u+(-y,x)]\cdot V=0.
\end{equation*}

\noindent It follows that 
\begin{equation}
\nabla u+(-y,x)=\tilde{D}N  \label{6.17}
\end{equation}%
\noindent for some function $\tilde{D}$ by the orthonormality of $V$ and $N.$
Taking the inner product of (\ref{6.17}) and $N$ gives%
\begin{equation}
\tilde{D}=N(u)+(-y,x)\cdot N=D\text{ along }\ell .  \label{6.18}
\end{equation}

\noindent Since $D$ $>$ $0,$ we may assume that $\tilde{D}$ $>$ $0$ as well
in a small neighborhood of $p.$ From (\ref{6.17}) we have $(u_{y},-u_{x})$ $%
+ $ $(x,y)$ $=$ ($\nabla u$ $+$ $(-y,x))^{\perp }$ $=$ $\tilde{D}V$ (recall
that $\vec{G}^{\perp }$ :$=$ $(G_{2},$ $-G_{1})$ for $\vec{G}$ $=$ $(G_{1},$ 
$G_{2})).$ Applying the divergence operator to both sides and expressing $V$ 
$=$ $(\sin \theta ,$ $-\cos \theta ),$ we obtain%
\begin{eqnarray}
2 &=&\func{div}(\tilde{D}V)  \label{6.18.1} \\
&=&V(\tilde{D})+\tilde{D}N(\theta ).  \notag
\end{eqnarray}

\noindent It follows that 
\begin{equation}
\frac{2-V(\tilde{D})}{\tilde{D}}=N(\theta ).  \label{6.19}
\end{equation}

\noindent Comparing (\ref{6.19}) with (\ref{6.15.1}) we get%
\begin{equation}
V(\tilde{D})=V(D)\text{ along }\ell  \label{6.20}
\end{equation}

\noindent in view of (\ref{6.18}). Now applying $V$ to (\ref{6.19})$,$ we
compute%
\begin{eqnarray}
V(\frac{2-V(\tilde{D})}{\tilde{D}}) &=&V(N(\theta ))  \label{6.21} \\
&=&[V,N](\theta )+N(V(\theta ))  \notag \\
&=&-N(\theta )^{2}-V(\theta )^{2}-N(H)  \notag
\end{eqnarray}

\noindent by (\ref{6.15}) and $V(\theta )$ $=$ $-H$, the second equation of (%
\ref{6.15.1}). Substituting (\ref{6.19}) and $V(\theta )$ $=$ $-H$ into (\ref%
{6.21}) we finally obtain 
\begin{equation*}
\tilde{D}\tilde{D}^{\prime \prime }=2(\tilde{D}^{\prime }-1)(\tilde{D}%
^{\prime }-2)+(H^{2}+N(H))\tilde{D}^{2}
\end{equation*}

\noindent in which we denote $V(\tilde{D}),$ $V(V(\tilde{D}))$ by $\tilde{D}%
^{\prime },$ $\tilde{D}^{\prime \prime },$ resp.. That is, $\tilde{D}$
satisfies the same equation as $D$ does in (\ref{1.9}) (noting that $P$ $=$ $%
N)$. Since $\tilde{D}$ and $D$ satisfy the same initial data by (\ref{6.18})
and (\ref{6.20}), we can conclude that $\tilde{D}$ $=$ $D$ in a small
neighborhood of $p$ by the uniqueness of the solution to an ordinary
differential equation of second order. Substituting $\tilde{D}$ $=$ $D$ into
(\ref{6.17}) we have%
\begin{equation*}
N=\frac{\nabla u+(-y,x)}{D}
\end{equation*}

\noindent and hence $D$ $=$ $|\nabla u+(-y,x)|$ $=$ $\sqrt{%
(u_{x}-y)^{2}+(u_{y}+x)^{2}}.$ (\ref{1.12.1}) follows from $V(\theta )$ $=$ $%
-H$ and observing that $\func{div}N$ $=$ $-V(\theta ).$ (\ref{1.12.2}) is
simply (\ref{6.18.1}) with $\tilde{D}$ $=$ $D.$

\endproof%

\bigskip

\section{Index of the singular set}

In this section we are going to show some global results about the singular
set.

\bigskip

\textbf{Lemma 7.1}. \textit{Let }$\Omega $\textit{\ be a bounded planar
domain.} \textit{Let }$u$\textit{\ }$\in $ $C^{1}(\Omega )$\textit{\ be a
weak solution to (}$\ref{1.1})$\textit{\ with }$\vec{F}$ $\in $ $%
C^{1}(\Omega )$ \textit{and }$H$\textit{\ }$\in $\textit{\ }$C^{0}(\Omega ).$%
\textit{\ Assume further }$N^{\perp }(curl\vec{F})$ \textit{and} $N(H)$%
\textit{\ exist and are continuous (extended over singular points) in }$%
\Omega $ \textit{and} $curl\vec{F}$\textit{\ }$\neq $\textit{\ }$0.$ \textit{%
Suppose the singular set }$S_{\vec{F}}(u)$\textit{\ }$\subset $ $\Omega $ 
\textit{is compact. Then \# }$\pi _{0}(S_{\vec{F}}(u))$\textit{\ }$<$\textit{%
\ }$\infty ,$\textit{\ i.e. the number of the connected components of }$S_{%
\vec{F}}(u)$\textit{\ is finite.}

\bigskip

\proof
Suppose \textit{\# }$\pi _{0}(S_{\vec{F}}(u))$\textit{\ }$=$\textit{\ }$%
\infty .$ Since $\Omega $ is bounded, there exist a sequence of distinct
connected components $S_{j},$ $j$ $=$ $1,$ $2,$ ...and a sequence of points $%
q_{j}$ $\in $ $S_{j}$ converging to $q_{\infty }$ $\in $ $\bar{\Omega}.$ It
follows from the compactness of $S_{\vec{F}}(u)$ that $q_{\infty }$ $\in $ $%
\Omega $ and hence $q_{\infty }$ $\in $ $S_{\vec{F}}(u)$ by closedness. From
Theorem C (b) there exists a neighborhood $U$ $\subset $ $\Omega $ of $%
q_{\infty }$ such that $U\cap S_{\vec{F}}(u)$ is path-connected. This
implies that $S_{K}$ $=$ $S_{K+1}$ $=$ $S_{K+2}$ $=$ ...for some large $K.$
We have reached a contradiction.

\endproof%

\bigskip

In the following we want to use "step functions" to approximate a $C^{0}$
singular curve. Let $\beta $ $:$ $[0,\tilde{\tau}]$ $\rightarrow $ $\Omega $
be a nonsingular $C^{1}$ smooth curve which is transverse to the
characteristic curves $\Gamma (\beta (\tau ))$ issuing from $\beta (\tau )$
for all $\tau $ $\in $ $[0,\tilde{\tau}]$ ($\beta (0)$ $=$ $\beta (\tilde{%
\tau})$ if $\beta $ is closed). Suppose each $\Gamma (\beta (\tau ))$ hits a
singular point $\mathfrak{s}(\tau )$ (see Figure 7.1)$.$

\begin{figure}[th]
\begin{center}
\includegraphics[width=6cm]{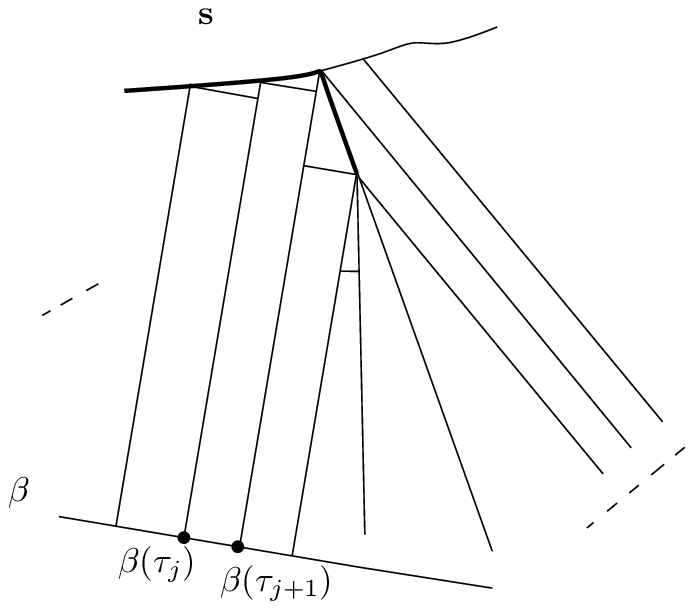}\\[0pt]
Figure 7.1
\end{center}
\par
\end{figure}

\bigskip

\textbf{Lemma 7.2.} \textit{Let }$\Omega $\textit{\ be a bounded planar
domain.} \textit{Let }$u$\textit{\ }$\in $ $C^{1}(\Omega )$\textit{\ be a
weak solution to (}$\ref{1.1})$\textit{\ with }$\vec{F}$ $\in $ $%
C^{1}(\Omega )$ \textit{and }$H$\textit{\ }$\in $\textit{\ }$C^{0}(\Omega ).$%
\textit{\ Assume further }$N^{\perp }(curl\vec{F})$ \textit{and} $N(H)$%
\textit{\ exist and are continuous (extended over singular points) in }$%
\Omega $ \textit{and} $curl\vec{F}$\textit{\ }$\neq $\textit{\ }$0.$ \textit{%
Suppose each of the characteristic curves passing through a nonsingular }$%
C^{1}$\textit{\ smooth curve }$\beta (\tau )$\textit{\ hits a singular point 
}$\mathfrak{s}(\tau )$\textit{\ as above. Then }$\mathfrak{s}$\textit{\ is
continuous.}

\bigskip

\proof
Suppose $\mathfrak{s}$ is not continuous at $\tau _{0}$ $\in $ $[0,\tilde{%
\tau}]$. Then we can find a sequence $\tau _{j}$ $\in $ $[0,\tilde{\tau}],$
converging to $\tau _{0},$ such that the characteristic curves $\Gamma
(\beta (\tau _{j}))$ converge ($C^{2}$ for any compact parameter interval)
to a characteristic curve $\Gamma _{\infty }$ and either $\Gamma _{\infty }$
contains $\mathfrak{s}(\tau _{0})$ or $\Gamma (\beta (\tau _{0}))$ $\cup $ $%
\{\beta (\tau _{0})\}$ contains $\bar{\Gamma}_{\infty }.$ Both cases
contradict Lemma $3.2^{\prime }$ (b) by noting that $\bar{\Gamma}_{\infty }$
contains a singular end point.

\endproof%

\bigskip

We remark that the result in Lemma 7.2 has been used in proving Lemma $%
3.5^{\prime }$. Recall that in Section 5 the map $P$ $:$ $(\sigma ,\tau )$ $%
\rightarrow $ $(x,y)$ describes $\Gamma (\beta (\tau ))$ with $P(0,\tau )$ $%
= $ $\beta (\tau ).$ Take a partition of $[0,\tilde{\tau}],$ $0$ $=$ $\tau
_{0} $ $\leq $ $\tau _{1}$ $\leq $ ...$\leq $ $\tau _{\kappa }$ $=$ $\tilde{%
\tau}. $ Let $\bar{\sigma}_{j}$ be such a number that either $P(\bar{\sigma}%
_{j},$ $\tau _{j})$ meets $\mathfrak{s}(\tau _{j})$ or $P(\bar{\sigma}_{j},$ 
$\tau _{j+1})$ meets $\mathfrak{s}(\tau _{j+1})$ for $j$ $=$ $0,$ $1,$..., $%
\kappa -1.$ Let $\sigma _{1}(\tau )$ denote the length of $\Gamma (\beta
(\tau ))$ from $\beta (\tau )$ to $\mathfrak{s}(\tau ).$ Note that $%
\mathfrak{s}$ and hence $\sigma _{1}$ is $C^{0}$ by Lemma 7.2 under the
assumption there and $P $ is $C^{1}$ smooth for $\tau $ $\in $ $[0,\tilde{%
\tau}]$ and $\sigma $ $\in $ $[0,\sigma _{1}(\tau ))$ originally, but in
fact $P$ extends $C^{1}$ smoothly over $\{(\sigma _{1}(\tau ),$ $\tau )\}$
if $H$\textit{\ }$\in $\textit{\ }$C^{1}(\Omega )$ besides the conditions in
Lemma 7.2$.$ See the argument in the proof of Theorem F in Section 5. So
there is small $\delta $ $>$ $0$ such that

\begin{eqnarray}
&&\sup_{\text{max}_{j}\text{\TEXTsymbol{\vert}}\tau _{j+1}-\tau _{j}|\leq
\delta }\sum_{j=0}^{\kappa -1}|P(\bar{\sigma}_{j},\tau _{j+1})-P(\bar{\sigma}%
_{j},\tau _{j})|  \label{7.2} \\
&\leq &\tilde{\tau}\sup_{(\sigma ,\tau )\in K}|\frac{\partial P}{\partial
\tau }(\sigma ,\tau )|  \notag
\end{eqnarray}

\noindent by the mean value theorem, where $K$ $=$ $K(\delta )$ is a compact
set in the $(\sigma ,\tau )$-plane, containing $\{(\sigma _{1}(\tau ),$ $%
\tau ):$ $\tau $ $\in $ $[0,\tilde{\tau}]\}.$

We remark that although the $C^{0}$ singular curve $\mathfrak{s}(\tau )$ may
not be rectifiable, the sum of the length of "steps" that approximate $%
\mathfrak{s}$ is bounded by (\ref{7.2}). Let $\Theta _{\vec{F}}$ $:=$ $du$ $%
+ $ $F_{1}dx$ $+$ $F_{2}dy$ $=$ ($u_{x}+F_{1})$ $dx$ $+$ $(u_{y}+F_{2})$ $%
dy. $ Let $Area(\beta ,\mathfrak{s})$ denote the area of region $R(\beta ,%
\mathfrak{s})$ surrounded by $\Gamma (\beta (\tilde{\tau})),$ $\mathfrak{s},$
$\Gamma (\beta (0)),$ and $\beta $ (see Figure 7.1 in which $\int_{\beta }$
means the integration from $\beta (0)$ to $\beta (\tilde{\tau})$).

\bigskip

\textbf{Lemma 7.3.} \textit{Suppose we are in the situation of Lemma 7.2.\
Assume further }$H$\textit{\ }$\in $\textit{\ }$C^{1}(\Omega ).$ \textit{%
Then if }$0$ $<$ $C_{1}$ $\leq $ $|curl\vec{F}|$ $\leq $ $C_{2}$ on region $%
R(\beta ,\mathfrak{s})$ \textit{we have}%
\begin{equation}
C_{1}Area(\beta ,\mathfrak{s})\leq |\int_{\beta }\Theta _{\vec{F}}|\leq
C_{2}Area(\beta ,\mathfrak{s}).  \label{7.5}
\end{equation}

\bigskip

\proof
Let $L_{j}$ denote the line segment from $P(\bar{\sigma}_{j},\tau _{j+1})$
to $P(\bar{\sigma}_{j},\tau _{j}).$ Let $\omega _{j}$ denote the region
surrounded by $\beta ([\tau _{j},\tau _{j+1}]),$ $P([0,\bar{\sigma}%
_{j}],\tau _{j+1}),$ $L_{j}$, and $P([0,\bar{\sigma}_{j}],\tau _{j}).$ By
Green's theorem we have%
\begin{eqnarray}
\int_{\omega _{j}}curl\vec{F}dx\wedge dy &=&\doint\limits_{\partial \omega
_{j}}\Theta _{\vec{F}}  \label{7.6} \\
&=&\int_{L_{j}}\Theta _{\vec{F}}+\int_{\beta ([\tau _{j},\tau
_{j+1}])}\Theta _{\vec{F}}  \notag
\end{eqnarray}

\noindent since $\Theta _{\vec{F}}$ $=$ $0$ when acting on $N^{\perp }(u)$.
Summing (\ref{7.6}) over $j$ we get%
\begin{equation}
\sum_{j=0}^{\kappa -1}\int_{\omega _{j}}curl\vec{F}dx\wedge dy-\int_{\beta
}\Theta _{\vec{F}}=\sum_{j=0}^{\kappa -1}\int_{L_{j}}\Theta _{\vec{F}}.
\label{7.7}
\end{equation}

\noindent Note that $D$ (:= $\sqrt{(u_{x}+F_{1})^{2}+(u_{y}+F_{2})^{2}})$ $=$
$0$ on $\mathfrak{s}$ (consisting of singular points) and $D$ is continuous$%
. $ Therefore given $\varepsilon $ $>$ $0,$ we have $D(q)$ $<$ $\varepsilon $
for $q$ $\in $ $\cup _{j=0}^{\kappa -1}L_{j}$ when $\kappa $ is large and $%
\max_{0\leq j\leq \kappa -1}$ $|\tau _{j+1}-\tau _{j}|$ is small enough. We
can then estimate%
\begin{eqnarray}
|\sum_{j=0}^{\kappa -1}\int_{L_{j}}\Theta _{\vec{F}}| &=&|\int_{\cup
_{j=0}^{\kappa -1}L_{j}}\Theta _{\vec{F}}|  \label{7.8} \\
&\leq &\int_{\cup _{j=0}^{\kappa -1}L_{j}}Dd\bar{s}\text{ \ (}\bar{s}:\text{%
arc-length parameter)}  \notag \\
&\leq &\varepsilon \tilde{\tau}\sup_{(\sigma ,\tau )\in K}|\frac{\partial P}{%
\partial \tau }(\sigma ,\tau )|  \notag
\end{eqnarray}

\noindent by (\ref{7.2}). Now (\ref{7.5}) follows from (\ref{7.7}) and (\ref%
{7.8}) as $\varepsilon $ $\rightarrow $ $0.$

\endproof%

\bigskip

\textbf{Lemma 7.4.} \textit{Let }$\Omega $\textit{\ be a bounded domain of }$%
R^{2}.$ \textit{Let }$u$\textit{\ }$\in $\textit{\ }$C^{1}(\Omega )$\textit{%
\ be a weak solution to (}$\ref{1.1})$\textit{\ with }$\vec{F}$\textit{\ }$%
\in $\textit{\ }$C^{1}(\Omega )$\textit{\ and }$H$\textit{\ }$\in $\textit{\ 
}$C^{1}(\Omega ).$\textit{\ Assume further }$N^{\perp }(curl\vec{F})$\textit{%
\ and }$N(H)$\textit{\ exist and are continuous (extended over singular
points) in }$\Omega .$\textit{\ Let }$\hat{S}$\textit{\ be a connected
component of }$S_{\vec{F}}(u).$ \textit{Suppose }$curl\vec{F}$ $>$ $0$%
\textit{\ (}$curl\vec{F}$ $<$ $0,$\textit{\ resp.) and }$\hat{S}$ \textit{is
compact and path-connected. Then for any }$\varepsilon $ $>$ $0,$ \textit{%
there exists a simple closed }$C^{0}$ \textit{curve }$\gamma _{\varepsilon }$%
\textit{\ in }$\Omega $\textit{\ such that }

(a) $dist(\gamma _{\varepsilon },\hat{S})$ $<$ $\varepsilon $ where $%
dist(\gamma _{\varepsilon },\hat{S})$ denotes the distance between $\gamma
_{\varepsilon }$ and $\hat{S}.$

(b) $N^{\perp }(u)$\textit{\ (}$-N^{\perp }(u),$ \textit{resp.) points
outward along }$\gamma _{\varepsilon },$\textit{\ i.e., the characteristic
curve} \textit{issuing from }$q$\textit{\ }$\in $\textit{\ }$\gamma
_{\varepsilon }$\textit{\ with tangent }$-N^{\perp }(u)$\textit{\ (}$%
N^{\perp }(u),$ \textit{resp.)} \textit{lies in the bounded domain, denoted
as }$\Omega _{\gamma _{\varepsilon }},$\textit{\ surrounded by }$\gamma
_{\varepsilon }.$

(c)\textit{\ }$\Omega _{\gamma _{\varepsilon }}$\textit{\ }$\subset $\textit{%
\ }$\Omega $ and $Area(\Omega _{\gamma _{\varepsilon }})$ (the area of $%
\Omega _{\gamma _{\varepsilon }})\rightarrow 0$ as $\varepsilon \rightarrow
0.$

\bigskip

\proof
Let $r_{0}$ :$=$ (2$\max_{p\in \hat{S}}|H(p)|)^{-1}.$ Let $\hat{\Omega}%
_{r_{1}}$ :$=$ $\{$ $p$ $\in $ $\Omega $ $|$ $0$ $<$ $dist(p,\hat{S})$ $<$ $%
r_{1}$ $\}.$ We claim the existence of $r_{1},$ $0$ $<$ $r_{1}$ $<$ $r_{0},$
such that for any $r_{2},$ $0$ $<$ $r_{2}$ $\leq $ $r_{1},$ there hold

(1) $\hat{\Omega}_{r_{2}}$ $\subset \subset $ $\Omega ,$

(2) $\hat{\Omega}_{r_{2}}\cap S_{\vec{F}}(u)$ $=$ $\phi $ (Note that $\hat{%
\Omega}_{r_{2}}$ does not contain $\hat{S}$ by definition)$,$ and

(3) the characteristic curve $\Gamma _{p}$ passing through $p$ hits $\hat{S}$
for any $p$ $\in $ $\hat{\Omega}_{r_{2}}.$

Conditions (1) and (2) are easily achieved. Suppose condition (3) fails.
That means we can find a sequence $p_{j}$ such that $dist(p_{j},\hat{S})$ $%
\rightarrow $ $0$ while $\Gamma _{p_{j}}$ does not hit $\hat{S}$ for each $%
j. $ Since $\hat{S}$ is compact, we end up finding a subsequence, still
denoted as $p_{j},$ converging to $p_{\infty }$ $\in $ $\hat{S}.$ On the
other hand, $\Gamma _{p_{j}}$ converges to a curve $\Gamma _{\infty }$ by
Lemma $3.2^{\prime }$ (a) and $\Gamma _{\infty }$ contains no singular
points by Lemma $3.2^{\prime }$ (b). But it is clear that the singular point 
$p_{\infty }$ $\in $ $\Gamma _{\infty }$. We have reached a contradiction.

Let $\mathfrak{s}(p)$ $\in $ $\hat{S}$ be the point at which $\Gamma _{p}$
hits $\hat{S}.$ Let $\mathring{\Gamma}_{p}$ $\subset $ $\Gamma _{p}$ denote
the part from $p$ to $\mathfrak{s}(p),$ not including end points $p$ and $%
\mathfrak{s}(p).$ By (3) and the choice of $r_{0}$ (recall that the
curvature of a characteristic curve is $-H),$ we conclude%
\begin{equation*}
\hat{\Omega}_{r_{2}}=\cup _{p\in l(r_{2})}\mathring{\Gamma}_{p}
\end{equation*}

\noindent where $l(r_{2})$ $:=$ $\{$ $p$ $\in $ $\Omega $ $|$ $dist(p,\hat{S}%
)$ $=$ $r_{2}$ $\}.$ Let $\gamma ^{r_{2}}$ $:=$ $\{$ $q$ $\in $ $\hat{\Omega}%
_{r_{2}}$ $|$ $\sigma _{1}(q)$ $:=$ the length of $\mathring{\Gamma}_{q}$ $=$
$\frac{r_{2}}{2}$ $\}.$ We claim that $\gamma ^{r_{2}}$ is a $C^{0}$ curve.
Take $q_{0}$ $\in $ $\gamma ^{r_{2}}.$ There is a point $p_{0}$ $\in $ $%
l(r_{2})$ such that $\mathring{\Gamma}_{p_{0}}$ $\supset $ $\mathring{\Gamma}%
_{q_{0}}.$ Take a $C^{1}$ smooth nonsingular curve $\beta $ $=$ $\beta (\tau
)$ $\subset $ $\hat{\Omega}_{r_{2}}$ $\cup $ $l(r_{2})$ for $\tau $ $\in $ $%
(\tau _{0}-\delta ,$ $\tau _{0}+\delta )$ with $\beta (\tau _{0})$ $=$ $%
p_{0},$ which is transversal to $\Gamma _{p_{0}}$ (a circular arc $\subset $ 
$\partial B_{r_{2}}(\mathfrak{s}_{0})$ where $\mathfrak{s}_{0}$ $\in $ $\hat{%
S}$ and $dist(p_{0},\mathfrak{s}_{0})$ $=$ $r_{2},$ passing through $p_{0},$
will serve as such $\beta $)$.$ Let $\mathfrak{s}(\tau )$ $\in $ $\hat{S}$
be the point which $\Gamma _{\beta (\tau )}$ hits. It follows that $%
\mathfrak{s}$ and hence $\sigma _{1}$ is $C^{0}$ by Lemma 7.2. Now any point
in $\gamma ^{r_{2}}$ near $q_{0}$ is the intersection of $\mathring{\Gamma}%
_{\beta (\tau )}$ and $\gamma ^{r_{2}}$ for a unique $\tau $ near $\tau
_{0}. $ We can therefore parametrize $\gamma ^{r_{2}}$ near $q_{0}$ by $\tau
,$ denoted as $\gamma ^{r_{2}}(\tau ).$ We compute%
\begin{eqnarray*}
\frac{r_{2}}{2} &=&\lim_{\tau _{j}\rightarrow \tau _{0}}\sigma _{1}(\gamma
^{r_{2}}(\tau _{j})) \\
&=&\sigma _{1}(q_{0}^{\prime })
\end{eqnarray*}

\noindent by the continuity of $\sigma _{1}$ if $\gamma ^{r_{2}}(\tau
_{j})\rightarrow q_{0}^{\prime }$ as $\tau _{j}\rightarrow \tau _{0}.$
Observe that $q_{0}^{\prime }$ $\in $ $\Gamma _{p_{0}}$ and hence $%
q_{0}^{\prime }$ $=$ $q_{0}$ since $\sigma _{1}(q_{0})$ $=$ $\frac{r_{2}}{2}$
too. So $\tau $ $\rightarrow $ $\gamma ^{r_{2}}(\tau )$ is continuous. In
fact, it is a homeomorphism near $\tau _{0}.$ We have shown that $\gamma
^{r_{2}}$ is a $C^{0}$ curve (in $R^{2}$ as a submanifold). Since $\hat{S}$
is compact, $\gamma ^{r_{2}}$ is bounded and hence any connected component
(still denoted as $\gamma ^{r_{2}})$ of $\gamma ^{r_{2}}$ must be simple
closed. Given $\varepsilon $ $>$ $0,$ we take $\gamma _{\varepsilon }$ $=$ $%
\gamma ^{\min \{\varepsilon ,r_{1}\}}.$ (a) follows from the fact that $%
\sigma _{1}(q)$ $\geq $ $dist(q,\mathfrak{s}(q)).$ Let $\Omega _{\gamma
_{\varepsilon }}$ denote the bounded domain surrounded by (a connected
component of) $\gamma _{\varepsilon }.$ Suppose the characteristic curve
issuing from $q$\ $\in $\ $\gamma _{\varepsilon }$\ with tangent $N^{\perp
}(u)$ lies in $\Omega _{\gamma _{\varepsilon }}$ (we are assuming $curl\vec{F%
}$ $>$ $0).$ Then the characteristic curve issuing from $q$\ $\in $\ $\gamma
_{\varepsilon }$\ with tangent $-N^{\perp }(u)$ hits a singular point $%
\mathfrak{s}(q)$ $\in $ $\hat{S}$ (we are assuming $curl\vec{F}$ $>$ $0,$ so 
$N^{\perp }(u)$ points away from $\hat{S}$ according to Theorem B (b)). From
(\ref{7.5}) with $\beta $ = $\gamma _{\varepsilon }$ in Lemma 7.3, we have%
\begin{equation}
C_{1}Area(\gamma _{\varepsilon },\mathfrak{s})\leq |\doint\limits_{\gamma
_{\varepsilon }}\Theta _{\vec{F}}|\leq C_{2}Area(\gamma _{\varepsilon },%
\mathfrak{s}).  \label{7.8.1}
\end{equation}

\noindent Note that we can take $C_{1}$ and $C_{2}$ to be independent of $%
\varepsilon .$

On the other hand, we compute%
\begin{eqnarray}
\doint\limits_{\gamma _{\varepsilon }}\Theta _{\vec{F}} &=&\int_{\Omega
_{\gamma _{\varepsilon }}}d\Theta _{\vec{F}}  \label{7.8.2} \\
&=&\int_{\Omega _{\gamma _{\varepsilon }}}curl\vec{F}dx\wedge dy\geq
C_{1}Area(\Omega _{\gamma _{\varepsilon }})  \notag
\end{eqnarray}

\noindent by Green's theorem. Choosing $\varepsilon $ small so that $%
Area(\gamma _{\varepsilon },\mathfrak{s})$ (close to 0) $<<$ $\frac{C_{1}}{%
C_{2}}Area(\Omega _{\gamma _{\varepsilon }})$ (note that $Area(\Omega
_{\gamma _{\varepsilon }})$ is nonincreasing in $\varepsilon ),$ we reach a
contradiction by substituting (\ref{7.8.2}) into (\ref{7.8.1}). So $\hat{S}$
must lie inside $\Omega _{\gamma _{\varepsilon }}$ and it is impossible for $%
\gamma _{\varepsilon }$ to have more than one connected components$.$ We
have proved (b). Since $Area(\gamma _{\varepsilon },\mathfrak{s})$ $%
\rightarrow $ $0$ as $\varepsilon $ $\rightarrow $ $0,$ we have 
\begin{equation*}
\lim_{\varepsilon \rightarrow 0}Area(\Omega _{\gamma _{\varepsilon }})=0
\end{equation*}%
\noindent by (\ref{7.8.1}) and (\ref{7.8.2}). That $\Omega _{\gamma
_{\varepsilon }}$\textit{\ }$\subset $\textit{\ }$\Omega $ follows by
observing that $\Omega _{\gamma _{\varepsilon }}$ $\subset $ the closure of $%
\hat{\Omega}_{^{\min \{\varepsilon ,r_{1}\}}}$\textit{\ }$\subset $ $\Omega .
$ We have proved (c). We can deal with the case $curl\vec{F}$ $<$ $0$
similarly.

\endproof%

\bigskip

\proof
\textbf{(of Theorem I) }By Lemma 7.1 we have \textit{\# }$\pi _{0}(S_{\vec{F}%
}(u))$ $<$ $\infty .$ Let $\hat{S}$ be a connected component of $S_{\vec{F}%
}(u).$ We claim that $\hat{S}$ is compact and path-connected. By Theorem C
(b) $\hat{S}$ is closed (in $\Omega $ and $R^{2})$, and hence compact since $%
S_{\vec{F}}(u)$ is compact. Suppose that $\hat{S}$ is decomposed as the
union of path-connected components $\hat{S}_{j},$ $j$ $=$ $1,$ $2,$.... Each 
$\hat{S}_{j}$ is closed by Theorem C (b), and hence compact since $\hat{S}$
is compact. Similarly we show that $\hat{S}^{\prime }$ $:=$ $\cup _{j\geq 2}%
\hat{S}_{j}$ is compact. Since we can separate two nonintersecting compact
sets by two nonintersecting open sets containing them respectively, $\hat{S}%
^{\prime }$ must be empty (otherwise contradicting the connectedness of $%
\hat{S})$ and hence $\hat{S}$ $=$ $\hat{S}_{1}$ is path-connected. Let $%
\gamma _{\varepsilon }$ be a simple closed $C^{0}$ curve around $\hat{S}$ in
Lemma 7.4.

Recall (see Section 1) that we denote the line field (1-dimensional
distribution) defined by the tangent lines (the lines having the direction $%
\pm N^{\perp }(u))$ of the characteristic curves by $\mathcal{D}.$ Let $%
P^{1} $ denote the projective line consisting of all lines through the
origin of $R^{2}.$ We define a homeomorphism $\varsigma $: $P^{1}$ $%
\rightarrow $ $S^{1} $ by noting that $P^{1}$ is the same as a semi-circle
with end points identified. Let $\approx $ denote a homeomorphism. We define 
$index(\gamma _{\varepsilon },\mathcal{D})$ to be half the degree of the map 
$S^{1}$ $\approx $ $\gamma _{\varepsilon }$ $\rightarrow $ $P^{1},$ defined
by $q$ $\in $ $\gamma _{\varepsilon }$ $\rightarrow $ $\mathcal{D(}q)$ $\in $
$P^{1}, $ composed with $\varsigma $ (cf. p.325 in \cite{Sp} for the index
of a line field at a point)$.$ It follows from Lemma 7.4 that%
\begin{equation}
index(\gamma _{\varepsilon },\mathcal{D})=1  \label{7.10}
\end{equation}

\noindent Namely, the connected component $\hat{S}$ of the singular set,
surrounded by $\gamma _{\varepsilon },$ has the index contribution $1.$

Next let us investigate the index contribution of the boundary curves. We
consider another copy of the domain $\Omega ,$ denoted as $\Omega ^{\prime
}. $ Denote the corresponding line field, boundary curves of $\Omega
^{\prime }$ by $\mathcal{D}^{\prime },$ $C_{j}^{\prime },$ $1$ $\leq $ $j$ $%
\leq $ $l,$ resp.$.$ We glue $C_{j}^{\prime }$ with $C_{j}$ for all $j$ to
get a closed surface $\Sigma $. Consider the line field $\mathcal{\tilde{D}}$
on $\Sigma , $ obtained from $\mathcal{D}$ and $\mathcal{D}^{\prime }$
(smoothing it along $\partial \Omega $ $=$ $\partial \Omega ^{\prime }$ so
that the topological type of the line field.does not change). Therefore the
index sum of $C_{j}^{\prime }$ $=$ $C_{j}$ $\subset $ $\Sigma $ with respect
to $\mathcal{\tilde{D}}$ equals $2$ times $index(C_{j};u)$ according to (\ref%
{1.12.3}) and (\ref{1.12.4})$.$ Together with (\ref{7.10}) we have 
\begin{equation}
\chi (\Sigma )=2\text{ }\#\text{ }\pi _{0}(S_{\vec{F}}(u))+2%
\sum_{j=1}^{l}index(C_{j};u)  \label{7.11}
\end{equation}

\noindent by the Hopf index theorem for a closed surface. Now (\ref{1.13})
follows from $\chi (\Omega )$ $=$ $\frac{1}{2}\chi (\Sigma )$ and (\ref{7.11}%
)$.$

\endproof%

\bigskip

We remark that (\ref{7.10}) is the key to the proof of Theorem I.

\bigskip

\proof
\textbf{(of Corollary J)} Observe that $\partial \Omega $ contains no
singular point by assumption. So either the characteristic line field is
transversal at each point of $\partial \Omega $ or there is a boundary point
at which the characteristic line is tangent to $\partial \Omega .$ In the
latter case, $\Omega $ must be foliated by the characteristic lines in view
of Lemma 3.2 (here we use the convexity of $\Omega )$. This contradicts the
assumption that $S(u)$ is nonempty in $\Omega .$ In the former case, the
index of $\partial \Omega $ is zero. Observe that a convex domain is
contractible to a point, and hence $\chi (\Omega )$ $=$ $1.$ Now $\#$\textit{%
\ }$\pi _{0}(S(u))$\textit{\ }$=$\textit{\ }$1$ follows from (\ref{1.13}).

\endproof%

\bigskip

Note that in Corollary J, if the singular set touches the boundary, we can
have more than one connected components of the singular set. For instance,
take a circular disc covering more than one singular line segments in Figure
4.5 (b). The boundary having no singular point can be obtained by the graph
restricted to it being a nonlegendrian (nonhorizontal) curve in $\boldsymbol{%
H}_{1}.$

\bigskip

\textbf{Example 7.1. }Look at Figure 7.2 (note that $\Omega $ is the
unbounded region):

\begin{figure}[ht]
\begin{center}
\includegraphics[width=6cm]{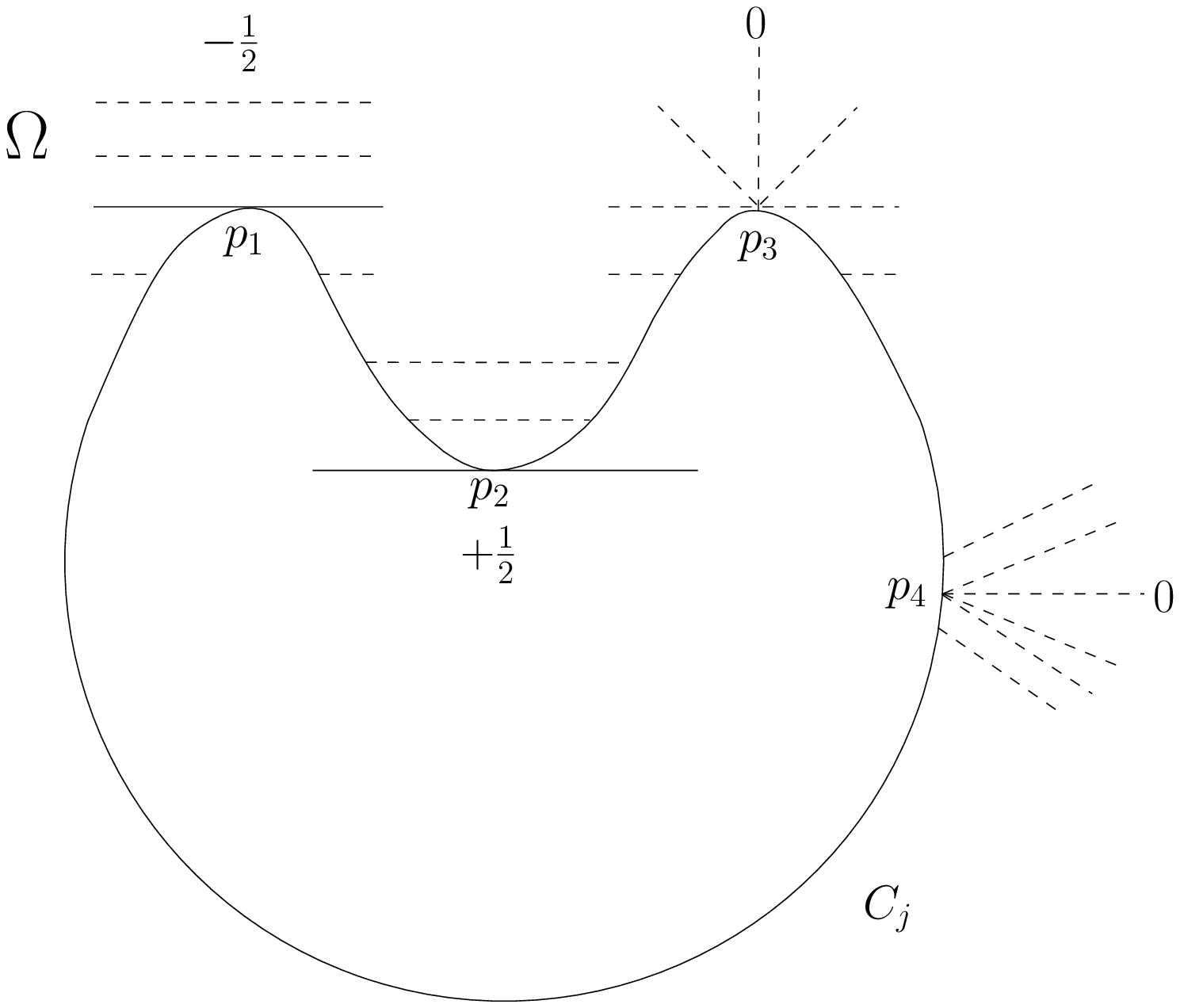}\\[0pt]
Figure 7.2
\end{center}
\par
\end{figure}

\noindent Let $\mathcal{D}_{U_{j}}$ denote the line field in a small
neighborhood $U_{j}$ of $p_{j}$ as shown in Figure 7.2 for $j$ $=$ $1,$ $2,$ 
$3,$ $4.$ According to the definition (\ref{1.12.3}) (see Figure 7.3),

\begin{figure}[th]
\begin{center}
\includegraphics[width=9cm]{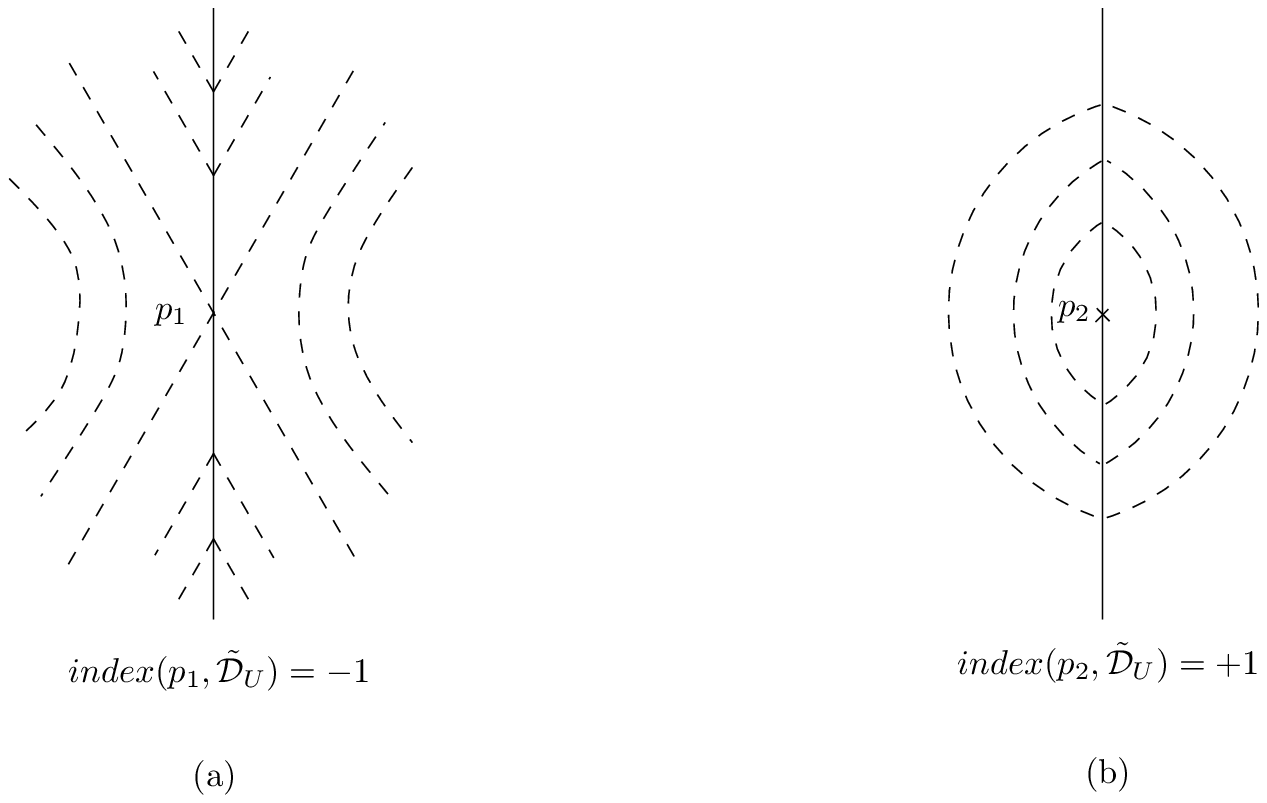}\\[0pt]
Figure 7.3
\end{center}
\par
\end{figure}

\noindent we have%
\begin{eqnarray*}
index(p_{1},\mathcal{D}_{U_{1}}) &=&-\frac{1}{2},\text{ }index(p_{2},%
\mathcal{D}_{U_{2}})=+\frac{1}{2}, \\
index(p_{3},\mathcal{D}_{U_{3}}) &=&index(p_{4},\mathcal{D}_{U_{4}})=0.
\end{eqnarray*}

\bigskip

\textbf{Example 7.2. }Let\textbf{\ }$\Omega $ be a $C^{1}$ smooth, bounded
planar domain with boundary curves $C_{1},$ $C_{2},$ and $C_{3}$ as shown in
Figure 7.4.

\begin{figure}[th]
\begin{center}
\includegraphics[width=7cm]{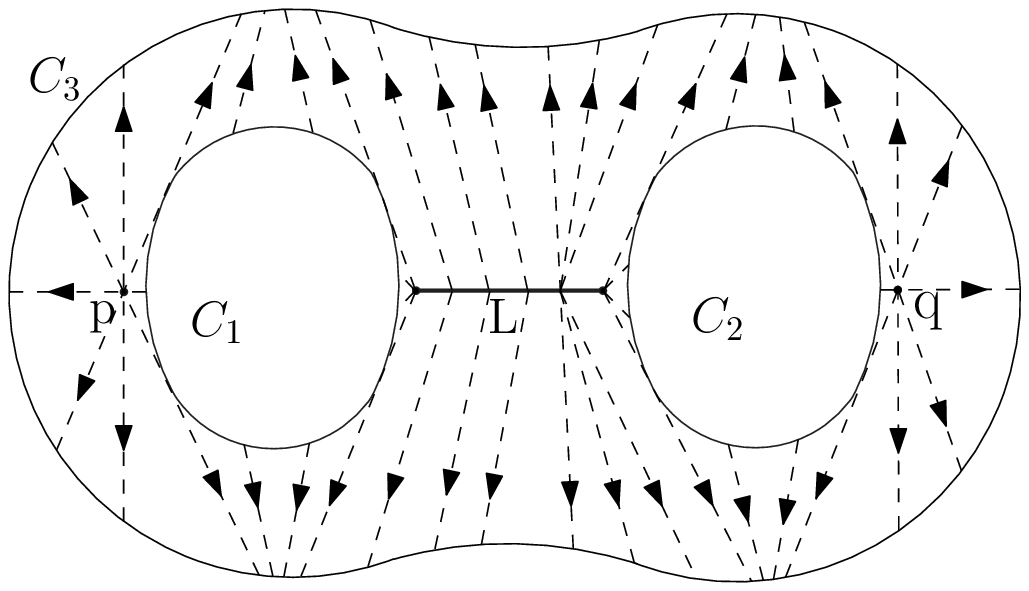}\\[0pt]
Figure 7.4
\end{center}
\par
\end{figure}

Suppose\textbf{\ }that a $C^{1}$\ smooth $p$-minimal graph defined by $u$\
over $\Omega $ has its singular set $S(u)$ $=$ $\{p\}$ $\cup $ $L$ $\cup $ $%
\{q\}$ where $L$ is a closed line segment. The dotted lines in Figure 7.4
denote the characteristic lines which are tangent to $\partial \Omega $ at
finitely many points, satisfying the situation to which we can apply Theorem
I. There are four points on $C_{1}$ each of which has index contribution $-%
\frac{1}{2}.$ Hence by (\ref{1.12.4}) we have%
\begin{equation}
index(C_{1};u)=4\cdot (-\frac{1}{2})=-2.  \label{7.12}
\end{equation}

\noindent Similarly we can easily compute%
\begin{eqnarray}
index(C_{2};u) &=&4\cdot (-\frac{1}{2})=-2,  \label{7.13} \\
index(C_{3};u) &=&0.  \notag
\end{eqnarray}

On the other hand, we have $\#$ $\pi _{0}(S(u))$ $=$ $3.$ So together with (%
\ref{7.12}) and (\ref{7.13}) we have%
\begin{eqnarray*}
&&\#\pi _{0}(S(u))+\sum_{j=1}^{3}index(C_{j};u) \\
&=&3+(-2)+(-2)+0=-1
\end{eqnarray*}

\noindent which equals the Euler characteristic number $\chi (\Omega )$ of $%
\Omega .$ We have verified Theorem I for this specific example.

\bigskip

\textbf{Example 7.3}. We consider closed surfaces of bounded $p$-mean
curvature in the Heisenberg group $\boldsymbol{H}_{1}.$ Let $\mathfrak{S}%
^{2} $ denote a Pansu sphere (\cite{Pan82}), a sphere of nonzero constant $p$%
-mean curvature in $\boldsymbol{H}_{1}$. There are exactly two singular
points on $\mathfrak{S}^{2},$ denoted as {\small N,S.} Take a short slit
along each of $g+1$ characteristic curves joining {\small N} to {\small S}
(see Figure 7.5 for $g$ $=$ 2).

\bigskip

\begin{figure}[th]
\begin{center}
\includegraphics[width=5cm]{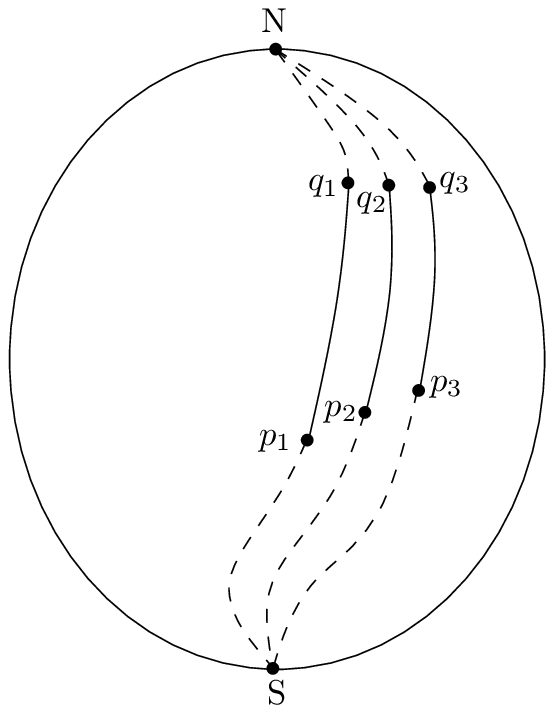}\\[0pt]
Figure 7.5
\end{center}
\par
\end{figure}

\bigskip

Take the 2-fold branched cover $\Sigma _{g}$ of $\mathfrak{S}^{2}$ with
branch locus consisting of the $2(g+1)$ points, $p_{j},$ $q_{j},$ $j$ $=$ $%
1, $ $...,$ $g+1,$ which are end points of the slits. From standard
topological arguments we learn that $\Sigma _{g}$ is a closed surface of
genus $g.$ Let $\hat{\varphi}$ : $\Sigma _{g}$ $\rightarrow $ $\mathfrak{S}%
^{2}$ $\subset $ $\boldsymbol{H}_{1}$ denote this 2-fold branched covering
map. Consider the characteristic line field on $\Sigma _{g}.$ Observe that
the index at a branch point is $-1$ while the index at {\small N} or {\small %
S} (two copies) is $+1.$ The total index count is correct since ($-1)\cdot
2(g+1)$ $+ $ $(+1)\cdot 4$ $=$ $2-2g,$ the Euler characteristic of $\Sigma
_{g}.$ We can deform the surface by moving the top copy of the Pansu sphere
a little bit away from the bottom copy (only do this away from the branch
cuts) to keep the $p$-mean curvature bounded while the only self
intersection of this new surface is the original branch cuts.

We can consider the problem of minimizing the $p$-area among mappings from
genus $g$ surfaces to $\boldsymbol{H}_{1}$ to enclose a fixed volume.
Observe that $\hat{\varphi}$ is not the minimizer, but only a critical point.

\bigskip

\section{Generalization to pseudohermitian manifolds}

Let $\Sigma $ be a (say, $C^{\infty })$ surface of a pseudohermitian
3-manifold $(M,J,\Theta )$ (see, e.g., \cite{CHMY04} for the definition of
pseudohermitian 3-manifolds). Let $\psi $ be a defining function of $\Sigma
. $ Let $\hat{e}_{1},$ $\hat{e}_{2}$ ($\hat{e}^{1},$ $\hat{e}^{2},$ resp.)
denote a local orthonormal basis (dual coframe, resp.) of the contact bundle
with respect to the Levi metric $\frac{1}{2}d\Theta (.,J.).$ Write%
\begin{equation}
\frac{\nabla _{b}\psi }{|\nabla _{b}\psi |}=\frac{\hat{e}_{1}\psi }{|\nabla
_{b}\psi |}\hat{e}_{1}+\frac{\hat{e}_{2}\psi }{|\nabla _{b}\psi |}\hat{e}%
_{2}=(\cos \theta )\hat{e}_{1}+(\sin \theta )\hat{e}_{2}  \label{O.1}
\end{equation}

\noindent for some angular function $\theta ,$ where $|\nabla _{b}\psi |=%
\sqrt{(\hat{e}_{1}\psi )^{2}+(\hat{e}_{2}\psi )^{2}}.$ Let $D:=|\nabla
_{b}\psi |.$ Recall (\cite{CHMY04}) that associated to the nonsingular part (%
$D\neq 0)$ of $\Sigma ,$ we have the characteristic vector field $e_{1}$,
the Legendrian normal $e_{2},$ and the dual coframe $e^{1},$ $e^{2}.$ Since $%
e_{1}\psi =0,$ we then have%
\begin{eqnarray}
D &=&|\nabla _{b}\psi |=\sqrt{(\hat{e}_{1}\psi )^{2}+(\hat{e}_{2}\psi )^{2}}
\label{O.2} \\
&=&\sqrt{(e_{1}\psi )^{2}+(e_{2}\psi )^{2}}=|e_{2}\psi |=e_{2}\psi  \notag
\end{eqnarray}%
\noindent by changing the sign of $\psi $ if necessary. Let $T$ denote the
Reeb vector field with respect to $\Theta .$ On $\Sigma $ we compute%
\begin{eqnarray}
0 &=&d\psi =(e_{1}\psi )e^{1}+(e_{2}\psi )e^{2}+(T\psi )\Theta  \label{O.3}
\\
&=&(e_{2}\psi )e^{2}+(T\psi )\Theta  \notag \\
&=&De^{2}+(T\psi )\Theta  \notag
\end{eqnarray}

\noindent by (\ref{O.2}). Taking the exterior differentiation of (\ref{O.3})
gives%
\begin{equation}
e_{1}(D)e^{1}\wedge e^{2}+Dde^{2}=-e_{1}(T\psi )e^{1}\wedge \Theta -T\psi
d\Theta .  \label{O.4}
\end{equation}

\noindent Here we have used the fact that $e^{2}\wedge \Theta =0$ on $\Sigma
.$ Formulas (A.3r) in \cite{CHMY04} when restricted to $\Sigma $ read%
\begin{eqnarray}
de^{1} &=&(H\alpha -\func{Re}A_{\bar{1}}^{1})e^{1}\wedge \Theta  \label{O.5}
\\
de^{2} &=&(\omega (\alpha e_{2}+T)-\func{Im}A_{\bar{1}}^{1})e^{1}\wedge
\Theta  \notag
\end{eqnarray}

\noindent (recall that $\omega _{1}^{1}=i\omega ,$ $\omega (e_{1})=H,$ the
p-mean curvature, $\alpha $ is defined so that $\alpha e_{2}+T\in T\Sigma ,$
and $e^{1}\wedge e^{2}=\alpha e^{1}\wedge \Theta $ on $\Sigma ).$
Substituting the second equality of (\ref{O.5}) into (\ref{O.4}) and noting
that $d\Theta =2e^{1}\wedge e^{2}=2\alpha e^{1}\wedge \Theta $ on $\Sigma ,$
we obtain%
\begin{equation}
e_{1}(D)\alpha +D(\omega (\alpha e_{2}+T)-\func{Im}A_{\bar{1}%
}^{1})=-e_{1}(T\psi )-2\alpha T\psi .  \label{O.6}
\end{equation}

\noindent If $T\psi =0,$ $T$ $\in $ $T\Sigma ,$ so it follows that $\alpha
=0.$ Let us assume $T\psi \neq 0.$ We may suppose $T\psi >0$ (otherwise
change $\Theta $ to $-\Theta ).$ Adjust $\psi $ such that $T\psi =1$ on $%
\Sigma $ while $e_{2}\psi $ is positive (choose for example $\frac{\psi }{%
T\psi }$). Under the condition $T\psi =1$ we can reduce (\ref{O.6}) to%
\begin{equation}
e_{1}(D)\alpha +D(\omega (\alpha e_{2}+T)-\func{Im}A_{\bar{1}}^{1})=-2\alpha
.  \label{O.7}
\end{equation}

Since ($\alpha e_{2}+T)\psi =0$ and $T\psi =1,$ we get%
\begin{equation}
\alpha =\frac{-1}{e_{2}\psi }=-\frac{1}{D}.  \label{O.8}
\end{equation}

\noindent Write%
\begin{eqnarray}
\nabla _{b}\psi &=&(\hat{e}_{1}\psi )\hat{e}_{1}+(\hat{e}_{2}\psi )\hat{e}%
_{2}=D(\cos \theta \hat{e}_{1}+\sin \theta \hat{e}_{2})  \label{O.9} \\
&=&(e_{1}\psi )e_{1}+(e_{2}\psi )e_{2}=(e_{2}\psi )e_{2}=De_{2}.  \notag
\end{eqnarray}

\noindent It follows that 
\begin{eqnarray}
e_{2} &=&\cos \theta \hat{e}_{1}+\sin \theta \hat{e}_{2},\text{ and }
\label{O.10} \\
e_{1} &=&\sin \theta \hat{e}_{1}-\cos \theta \hat{e}_{2},\text{ }e^{1}=\sin
\theta \hat{e}^{1}-\cos \theta \hat{e}^{2}  \notag \\
e^{2} &=&\cos \theta \hat{e}^{1}+\sin \theta \hat{e}^{2}.  \notag
\end{eqnarray}

\noindent From (\ref{O.10}) we learn that 
\begin{equation}
\theta ^{1}=e^{i(\frac{\pi }{2}-\theta )}\hat{\theta}^{1}.  \label{O.10'}
\end{equation}
\noindent Therefore $A_{11}=e^{-2i(\frac{\pi }{2}-\theta )}\hat{A}%
_{11}=-e^{2i\theta }\hat{A}_{11}$ and hence%
\begin{eqnarray}
\func{Re}A_{11} &=&\sin 2\theta \func{Im}\hat{A}_{11}-\cos 2\theta \func{Re}%
\hat{A}_{11}  \label{O.11} \\
\func{Im}A_{11} &=&-\sin 2\theta \func{Re}\hat{A}_{11}-\cos 2\theta \func{Im}%
\hat{A}_{11}.  \notag
\end{eqnarray}

\noindent The connection forms change according to $\omega _{1}^{1}=\hat{%
\omega}_{1}^{1}-id(\frac{\pi }{2}-\theta ).$ It follows that%
\begin{equation}
\omega =\hat{\omega}+d\theta .  \label{O.12}
\end{equation}

\noindent Let%
\begin{equation}
v_{2}:=\alpha e_{2}+T.  \label{O.13}
\end{equation}

\noindent From (\ref{O.8}), (\ref{O.12}), and \ref{O.13}, we can rewrite (%
\ref{O.7}) as%
\begin{equation}
e_{1}(D)=D^{2}(v_{2}(\theta )+\hat{\omega}(v_{2})-\func{Im}A_{\bar{1}%
}^{1})-2.  \label{O.14}
\end{equation}

\noindent Taking the derivative of (\ref{O.14}) in the direction of $e_{1}$,
we have%
\begin{eqnarray}
e_{1}^{2}(D) &=&2De_{1}(D)(v_{2}(\theta )+\hat{\omega}(v_{2})-\func{Im}A_{%
\bar{1}}^{1})  \label{O.15} \\
&&+D^{2}\{e_{1}(v_{2}(\theta ))+e_{1}(\hat{\omega}(v_{2})-\func{Im}A_{\bar{1}%
}^{1})\}.  \notag
\end{eqnarray}

On the other hand, we write%
\begin{equation}
e_{1}(v_{2}(\theta ))=v_{2}(e_{1}(\theta ))+[e_{1},v_{2}](\theta )
\label{O.16}
\end{equation}

\noindent and observe that 
\begin{equation}
e_{1}(\theta )=\omega (e_{1})-\hat{\omega}(e_{1})=H-\hat{\omega}(e_{1})
\label{O.17}
\end{equation}

\noindent by applying (\ref{O.12}) to $e_{1}.$ Making use of (A.6r) and
(A.7r) in \cite{CHMY04}, we have%
\begin{eqnarray}
\lbrack e_{1},v_{2}] &=&\alpha \lbrack e_{1},e_{2}]+e_{1}(\alpha
)e_{2}+[e_{1},T]  \label{O.18} \\
&=&(-\alpha H+\func{Re}A_{11})e_{1}+(-\omega (v_{2})+e_{1}(\alpha )-\func{Im}%
A_{11})e_{2}-2\alpha T  \notag \\
&=&(-\alpha H+\func{Re}A_{11})e_{1}-2\alpha v_{2}.  \notag
\end{eqnarray}

\noindent For the last equality of (\ref{O.18}) we have used the
integrability condition for $\alpha :$%
\begin{equation}
-\omega (v_{2})+e_{1}(\alpha )-\func{Im}A_{11}=-2\alpha ^{2}  \label{O.19}
\end{equation}

\noindent (in order to make $[e_{1},v_{2}]\in T\Sigma ).$ By (\ref{O.16}), (%
\ref{O.17}), (\ref{O.18}), and (\ref{O.14}), we obtain%
\begin{eqnarray}
e_{1}(v_{2}(\theta )) &=&v_{2}(H-\hat{\omega}(e_{1}))  \label{O.20} \\
&&+(-\alpha H+\func{Re}A_{11})e_{1}(\theta )-2\alpha v_{2}(\theta )  \notag
\\
&=&(v_{2}-\alpha H+\func{Re}A_{11})(H-\hat{\omega}(e_{1}))  \notag \\
&&-2\alpha \lbrack \frac{e_{1}(D)+2}{D^{2}}-\hat{\omega}(v_{2})-\func{Im}%
A_{11}].  \notag
\end{eqnarray}

\noindent Note that $\func{Im}A_{\bar{1}}^{1}=\func{Im}A_{\bar{1}\bar{1}}=-%
\func{Im}A_{11}.$ Next we observe from (\ref{O.18}), (A.5r), and (A.5) in 
\cite{CHMY04} that 
\begin{eqnarray}
&&e_{1}(\hat{\omega}(v_{2}))-v_{2}(\hat{\omega}(e_{1})-(-\alpha H+\func{Re}%
A_{11})\hat{\omega}(e_{1})+2\alpha \hat{\omega}(v_{2})  \label{O.21} \\
&=&e_{1}(\hat{\omega}(v_{2}))-v_{2}(\hat{\omega}(e_{1})-\hat{\omega}%
([e_{1},v_{2}])  \notag \\
&=&d\hat{\omega}(e_{1},v_{2})  \notag \\
&=&-2\alpha W+2\func{Im}A_{11,\bar{1}}.  \notag
\end{eqnarray}

\noindent For the last equality of (\ref{O.21}) we have used the
transformation law%
\begin{eqnarray*}
\func{Im}A_{11,\bar{1}} &=&\sin \theta \func{Im}\hat{A}_{11,\bar{1}}-\cos
\theta \func{Re}\hat{A}_{11,\bar{1}} \\
(\func{Re}A_{11,\bar{1}} &=&\sin \theta \func{Re}\hat{A}_{11,\bar{1}}+\cos
\theta \func{Im}\hat{A}_{11,\bar{1}})
\end{eqnarray*}

\noindent under the coframe change (\ref{O.10'}). Denote $\frac{v_{2}}{%
\alpha }=e_{2}+\frac{T}{\alpha }$ by $\tilde{e}_{2}.$ Substituting 
\begin{equation*}
v_{2}(\theta )+\hat{\omega}(v_{2})-\func{Im}A_{\bar{1}}^{1}=\frac{e_{1}(D)+2%
}{D^{2}}
\end{equation*}%
\noindent from (\ref{O.14}) and (\ref{O.20}) into (\ref{O.15}), we finally
reach%
\begin{eqnarray}
&&e_{1}^{2}(D)  \label{O.23} \\
&=&\frac{2(e_{1}(D)+1)(e_{1}(D)+2)}{D}+D[-\tilde{e}_{2}(H)+H^{2}]  \notag \\
&&+D^{2}\{2D^{-1}W+2\func{Im}A_{11,\bar{1}}+(e_{1}-2D^{-1})(\func{Im}%
A_{11})+(\func{Re}A_{11})H\}  \notag
\end{eqnarray}

\noindent in view of (\ref{O.21}) and $\alpha =-\frac{1}{D}$ by (\ref{O.8}).
Note that for a graph over the $xy$-plane in the 3-dimensional Heisenberg
group, $e_{1}(D)=-N^{\perp }(D)$ and $\tilde{e}_{2}(H)=-N(H)$ for $H$ being
a function of $x$ and $y.$ So (\ref{O.23}) is reduced to previous formula (%
\ref{1.2}) with $\vec{F}$ $=$ $(-y,x)$.

\bigskip

\textbf{Example 8.1}. We want to show that $\func{div}\frac{\nabla u+\vec{F}%
}{|\nabla u+\vec{F}|}$ $(u$ $\in $ $C^{2},$ say) is proportional to the $p$%
-mean curvature of a certain pseudohermitian structure. Let $\Theta _{\vec{F}%
}$ :$=$ $dz$ $+$ $F_{1}dx$ $+$ $F_{2}dy$ where $\vec{F}$ $=$ $(F_{1},$ $%
F_{2}).$ Then we have%
\begin{equation*}
d\Theta _{\vec{F}}=(curl\vec{F})dx\wedge dy
\end{equation*}

\noindent where $curl\vec{F}$ $:=$ $\frac{\partial F_{2}}{\partial x}-\frac{%
\partial F_{1}}{\partial y}.$ We require that $curl\vec{F}$ $\neq $ $0.$
Without loss of generality we may assume $curl\vec{F}$ $>$ $0.$ We have $%
\Theta _{\vec{F}}\wedge d\Theta _{\vec{F}}$ $\neq $ $0.$ That is to say, $%
\Theta _{\vec{F}}$ is a contact form. Let 
\begin{eqnarray}
\hat{e}_{1} &=&\frac{1}{\sqrt{curl\vec{F}/2}}(\frac{\partial }{\partial x}%
-F_{1}\frac{\partial }{\partial z})  \label{5.1.1} \\
\hat{e}_{2} &=&\frac{1}{\sqrt{curl\vec{F}/2}}(\frac{\partial }{\partial y}%
-F_{2}\frac{\partial }{\partial z}).  \notag
\end{eqnarray}

\noindent It is clear that $\hat{e}_{1},$ $\hat{e}_{2}$ $\in $ $Ker\Theta _{%
\vec{F}}$ form a basis$.$ Define the $CR$ structure $J_{\vec{F}}$ on this
basis by $J_{\vec{F}}(\hat{e}_{1})$ $:=$ $\hat{e}_{2}$ and $J_{\vec{F}}(\hat{%
e}_{2})$ $:=$ $-\hat{e}_{1}.$ It follows that $\hat{e}_{1},$ $\hat{e}_{2}$
form an orthonormal basis with respect to the Levi metric $\frac{1}{2}%
d\Theta _{\vec{F}}(\cdot ,J_{\vec{F}}\cdot ).$ Let $\psi $ $:=$ $z-u(x,y)$
be a defining function. Then we have%
\begin{eqnarray}
\hat{e}_{1}\psi &=&\frac{1}{\sqrt{curl\vec{F}/2}}(-u_{x}-F_{1})
\label{5.1.2} \\
\hat{e}_{2}\psi &=&\frac{1}{\sqrt{curl\vec{F}/2}}(-u_{y}-F_{2}).  \notag
\end{eqnarray}

\noindent So from $\nabla _{b}\psi $ $=$ ($\hat{e}_{1}\psi )\hat{e}_{1}$ $+$
($\hat{e}_{2}\psi )\hat{e}_{2},$ we have%
\begin{eqnarray}
|\nabla _{b}\psi | &=&\sqrt{(\hat{e}_{1}\psi )^{2}+(\hat{e}_{2}\psi )^{2}}
\label{5.1.3} \\
&=&\frac{1}{\sqrt{curl\vec{F}/2}}|\nabla u+\vec{F}|  \notag
\end{eqnarray}

\noindent by (\ref{5.1.2}). Now we can compute the $p$-mean curvature $H_{J_{%
\vec{F}},\Theta _{\vec{F}}}$ of the graph $z$ $=$ $u(x,y)$ with respect to
the pseudohermitian structure ($J_{\vec{F}},$ $\Theta _{\vec{F}})$ according
to formula $(pMCE)$ in \cite{CHMY04}:%
\begin{eqnarray*}
H_{J_{\vec{F}},\Theta _{\vec{F}}} &=&-\func{div}_{b}\frac{\nabla _{b}\psi }{%
|\nabla _{b}\psi |} \\
&=&-\hat{e}_{1}(\frac{\hat{e}_{1}\psi }{|\nabla _{b}\psi |})-\hat{e}_{2}(%
\frac{\hat{e}_{2}\psi }{|\nabla _{b}\psi |}) \\
&=&\frac{1}{\sqrt{curl\vec{F}/2}}\func{div}\frac{\nabla u+\vec{F}}{|\nabla u+%
\vec{F}|}
\end{eqnarray*}

\noindent by (\ref{5.1.2}) and (\ref{5.1.3}).

\bigskip

We may write equation (\ref{O.23}) in the form of (\ref{A1}) (see the
Appendix). By Theorem A.1 we can then conclude

\bigskip

\textbf{Theorem B}$^{\prime }.$ \textit{Let }$\Sigma $\textit{\ be a (say, }$%
C^{\infty })$\textit{\ surface of a pseudohermitian 3-manifold }$(M,J,\Theta
),$\textit{\ defined by }$\{\psi $\textit{\ }$=$\textit{\ }$0\}.$\textit{\
Suppose that }$T\psi $\textit{\ }$\neq $\textit{\ }$0$\textit{. Then either }%
$e_{1}(D)$\textit{\ tends to }$-1$\textit{\ or }$e_{1}(D)$\textit{\ tends to 
}$-2$\textit{\ as the argument tends to a singular point along a
characteristic curve}

\bigskip

\section{Appendix: a generalized ODE}

We generalize the result (\ref{1.6}) for equation (\ref{1.2}) in this
section.

\bigskip

\textbf{Theorem A.1. }\textit{Let }$v$ $\in $\textit{\ }$C^{2}(0,\rho _{0}),$%
\textit{\ }$E_{1},$\textit{\ }$l,$\textit{\ }$m$\textit{\ }$\in $\textit{\ }$%
C^{0}[0,\rho _{0})$\textit{\ (}$\rho _{0}$\textit{\ }$>$\textit{\ }$0$%
\textit{) be real functions such that }$v(\rho )$\textit{\ }$>$\textit{\ }$0$%
\textit{\ and }$v$ \textit{is bounded on }$(0,\rho _{0}),$\textit{\ }$%
E_{1}(0)$\textit{\ }$>$\textit{\ }$0,$\textit{\ }$l(0)$\textit{\ }$<$\textit{%
\ }$m(0).$\textit{\ Let }$E_{2}$\textit{\ }$=$\textit{\ }$E_{2}(\rho ,v)$%
\textit{\ be a real function continuous in }$(\rho ,v)$\textit{\ }$\in $%
\textit{\ }$[0,\rho _{0})$\textit{\ }$\times $\textit{\ }$\mathit{[0,\infty )%
}$\textit{\ with }$E_{2}(0,0)$\textit{\ }$=$\textit{\ }$0.$\textit{\
Consider the following ODE (which generalizes (\ref{1.2})):}%
\begin{equation}
v(\rho )v^{\prime \prime }(\rho )=E_{1}(\rho )(v^{\prime }(\rho )-l(\rho
))(v^{\prime }(\rho )-m(\rho ))+E_{2}(\rho ,v(\rho )).  \label{A1}
\end{equation}

\textit{\noindent Then there hold}

\textit{(a) The limit of }$v(\rho )$\textit{\ and }$v^{\prime }(\rho ),$%
\textit{\ resp. exists as }$\rho $\textit{\ }$\rightarrow $\textit{\ }$0.$%
\textit{\ }

\textit{(b) In case }$\lim_{\rho \rightarrow 0}v(\rho )$\textit{\ }$=$%
\textit{\ }$0,$ \textit{we have either}%
\begin{equation}
\lim_{\rho \rightarrow 0}v^{\prime }(\rho )=l(0)\text{ \ or \ }\lim_{\rho
\rightarrow 0}v^{\prime }(\rho )=m(0).  \label{A1'}
\end{equation}

\bigskip

\textbf{Lemma A.2. }\textit{Suppose we are in the situation of Theorem A.1
(excluding the condition }$v(\rho )$ $>$ $0$ for $0$ $<$ $\rho $ $<$ $\rho
_{0}$ \textit{and replacing }$E_{1}(0)$\textit{\ }$>$\textit{\ }$0$\textit{\
by }$E_{1}(0)$\textit{\ }$\neq $\textit{\ }$0)$\textit{. Then for any }$%
0<\varepsilon <\frac{m(0)-l(0)}{3}$\textit{\ there exists }$\delta $\textit{%
\ }$=$\textit{\ }$\delta (\varepsilon )$\textit{\ (}$<\rho _{0})$\textit{\
such that for any }$0<\rho <\delta $\textit{\ there holds}%
\begin{equation}
|v(\rho )v^{\prime \prime }(\rho )|\geq \frac{|E_{1}(0)|}{16}|(v^{\prime
}(\rho )-l(0))(v^{\prime }(\rho )-m(0))|.  \label{A2'}
\end{equation}

\textit{\noindent for }$v^{\prime }(\rho )$\textit{\ }$\in $\textit{\ }$%
(-\infty ,$\textit{\ }$l(0)-\varepsilon ]$\textit{\ }$\cup $\textit{\ }$%
[l(0)+\varepsilon ,$\textit{\ }$m(0)-\varepsilon ]$\textit{\ }$\cup $\textit{%
\ }$[m(0)+\varepsilon ,$\textit{\ }$\infty ).$

\textit{\bigskip }

\proof
Write $E_{1}(\rho )$ $=$ $E_{1}(0)$ $+$ $E_{1}(\rho )$ $-$ $E_{1}(0).$ There
exists $\delta _{1}$ $>$ $0$ such that $|E_{1}(\rho )$ $-$ $E_{1}(0)|$ $\leq 
$ $\frac{|E_{1}(0)|}{2}$ for $0$ $<$ $\rho $ $<$ $\delta _{1}$ by the
continuity of $E_{1}$ and the assumption $E_{1}(0)$ $\neq $ $0.$ It follows
that%
\begin{equation}
|E_{1}(\rho )|\geq \frac{|E_{1}(0)|}{2}  \label{A3}
\end{equation}

\noindent for $0$ $<$ $\rho $ $<$ $\delta _{1}.$ Choose $\delta _{2}$ $>$ $0$
such that $|l(\rho )-l(0)|$ $\leq $ $\frac{\varepsilon }{2}$ for $0$ $<$ $%
\rho $ $<$ $\delta _{2}.$ So we have%
\begin{eqnarray}
|v^{\prime }(\rho )-l(\rho )| &=&|v^{\prime }(\rho )-l(0)+l(0)-l(\rho )|
\label{A4} \\
&\geq &\frac{|v^{\prime }(\rho )-l(0)|}{2}+\frac{|v^{\prime }(\rho )-l(0)|}{2%
}-|l(\rho )-l(0)|  \notag \\
&\geq &\frac{|v^{\prime }(\rho )-l(0)|}{2}  \notag
\end{eqnarray}%
\noindent for $0$ $<$ $\rho $ $<$ $\delta _{2}$ since $|v^{\prime }(\rho
)-l(0)|$ $\geq $ $\varepsilon $ by assumption. Similarly we can choose $%
\delta _{3}$ $>$ $0$ such that 
\begin{equation}
|v^{\prime }(\rho )-m(\rho )|\geq \frac{|v^{\prime }(\rho )-m(0)|}{2}
\label{A5}
\end{equation}

\noindent for $0$ $<$ $\rho $ $<$ $\delta _{3}.$ By continuity and $%
E_{2}(0,0)$ $=$ $0,$ we can find $\delta _{4}$ $>$ $0$ such that 
\begin{equation}
|E_{2}(\rho ,v(\rho ))|\leq \frac{|E_{1}(0)|}{16}|(v^{\prime }(\rho
)-l(0))(v^{\prime }(\rho )-m(0))|  \label{A6}
\end{equation}

\noindent for $0$ $<$ $\rho $ $<$ $\delta _{4}.$ Note that the right-hand
side of (\ref{A6}) is greater or equal to $\frac{|E_{1}(0)|}{16}\varepsilon
^{2}$ $>$ $0.$ From (\ref{A1}), (\ref{A3}), (\ref{A4}), (\ref{A5}), and (\ref%
{A6}), we have%
\begin{eqnarray*}
|v(\rho )v^{\prime \prime }(\rho )| &=&|E_{1}(\rho )(v^{\prime }(\rho
)-l(\rho ))(v^{\prime }(\rho )-m(\rho ))+E_{2}(\rho ,v(\rho ))| \\
&\geq &|E_{1}(\rho )||(v^{\prime }(\rho )-l(\rho ))||(v^{\prime }(\rho
)-m(\rho ))|-|E_{2}(\rho ,v(\rho ))| \\
&\geq &\frac{|E_{1}(0)|}{16}|(v^{\prime }(\rho )-l(0))(v^{\prime }(\rho
)-m(0))|
\end{eqnarray*}

\noindent for $0$ $<$ $\rho $ $<$ $\delta =\min \{\delta _{1},\delta
_{2},\delta _{3},\delta _{4},\rho _{0}\}.$ We have proved (\ref{A2'}).

\endproof%

\bigskip

\proof%
\textbf{\ (of Theorem A.1)} Suppose that $\lim_{\rho \rightarrow 0}v(\rho )$
does not exist. Then $a$ $:=$ $\lim \inf_{\rho \rightarrow 0}v(\rho )$ $<$ $%
\lim \sup_{\rho \rightarrow 0}v(\rho )$ $:=$ $b.$ $a,$ $b$ $\in $ $R$ since $%
v$ is bounded by assumption. There exist sequences $\bar{\rho}_{j}$ and $%
\rho _{j}$ such that $v(\bar{\rho}_{j})$ $\rightarrow $ $b$ as $\bar{\rho}%
_{j}$ $\rightarrow $ $0$ while $v(\rho _{j})$ $\rightarrow $ $a$ as $\rho
_{j}$ $\rightarrow $ $0.$ Now as both $\bar{\rho}_{j}$ and $\rho _{k}$ tend
to $0,$ we have%
\begin{eqnarray*}
v^{\prime }(\xi _{jk}) &=&\frac{v(\bar{\rho}_{j})-v(\rho _{k})}{\bar{\rho}%
_{j}-\rho _{k}}\text{ \textit{for} }\mathit{\xi }_{jk}\text{ \textit{between}
}\mathit{\bar{\rho}}_{j}\text{ \textit{and} }\mathit{\rho }_{k}\text{ } \\
&&\text{\textit{approximates} }\frac{b-a}{\bar{\rho}_{j}-\rho _{k}}
\end{eqnarray*}

\noindent which goes to $+\infty $ or $-\infty $ depending on $\bar{\rho}%
_{j}-\rho _{k}$ is positive or negative, resp.. On the other hand, we always
have $v^{\prime \prime }$ $>$ $0$ according to (\ref{A1}) for \TEXTsymbol{%
\vert}$v^{\prime }|$ large. We claim that it is impossible for $v^{\prime }$
to change drastically from positive large to negative large while $v^{\prime
\prime }$ is positive. Suppose that $v^{\prime }(\xi _{j_{1}k_{1}})$ is
positively large while $v^{\prime }(\xi _{j_{2}k_{2}})$ is negatively large
for $0$ $<$ $\xi _{j_{1}k_{1}}$ $<$ $\xi _{j_{2}k_{2}}.$ Then there exists $%
\hat{\xi},$ $\xi _{j_{1}k_{1}}$ $<$ $\hat{\xi}$ $<$ $\xi _{j_{2}k_{2}}$,
such that $v^{\prime }(\hat{\xi})$ $>$ $v^{\prime }(\xi _{j_{1}k_{1}})$ and $%
v^{\prime \prime }(\hat{\xi})$ $=$ $0$ since $v^{\prime \prime }(\xi
_{j_{1}k_{1}})$ $>$ $0.$ But applying (\ref{A1}) to $\rho $ $=$ $\hat{\xi},$
we get $0$ $=$ \textit{a positive large number}, a contradiction. We have
proved the existence of $\lim_{\rho \rightarrow 0}v(\rho )$.

Suppose $\lim_{\rho \rightarrow 0}v(\rho )$ $=$ $0$ and (\ref{A1'}) does not
hold. Then there exists a sequence of $\rho _{j}\rightarrow 0$ such that 
\begin{equation}
v^{\prime }(\rho _{j})\mathit{\ }\in \mathit{\ }(-\infty ,\mathit{\ }%
l(0)-\varepsilon ]\mathit{\ }\cup \mathit{\ }[l(0)+\varepsilon ,\mathit{\ }%
m(0)-\varepsilon ]\mathit{\ }\cup \mathit{\ }[m(0)+\varepsilon ,\mathit{\ }%
\infty )  \label{A6'}
\end{equation}%
\noindent for a given $0<\varepsilon <\frac{m(0)-l(0)}{3}.$ By Lemma A.2
(may assume $0$ $<$ $\rho _{j}$ $<$ $\delta ),$ we have%
\begin{equation}
|v(\rho _{j})v^{\prime \prime }(\rho _{j})|\geq \frac{|E_{1}(0)|}{16}%
|(v^{\prime }(\rho _{j})-l(0))(v^{\prime }(\rho _{j})-m(0))|.  \label{A7}
\end{equation}

\noindent Since $v(\rho _{j})\rightarrow 0$ as $\rho _{j}\rightarrow 0,$ $%
|v^{\prime \prime }(\rho _{j})|$ must tend to infinity in view of (\ref{A7})
and (\ref{A6'}). We may assume either a subsequence of $v^{\prime \prime
}(\rho _{j})$ goes to $+\infty $ or a subsequence of $v^{\prime \prime
}(\rho _{j})$ goes to $-\infty .$ Still denote the subsequence by $v^{\prime
\prime }(\rho _{j}).$ Observe from (\ref{A1}) and $v(\rho )$\textit{\ }$>$%
\textit{\ }$0$\textit{\ }for ($0$ $<)$ $\rho $ small that 
\begin{eqnarray}
v^{\prime \prime }\text{(}\rho \text{) \textit{has the same sign (}}\mathit{>%
}\text{ }\mathit{0}\text{ }\mathit{or}\text{ }\mathit{<}\text{ }\mathit{0)}%
\text{ }\mathit{for}\text{ }\mathit{all}\text{ }\rho && \text{ }  \label{A8}
\\
\text{\textit{small enough so that} }\mathit{v}^{\prime }\mathit{(\rho )}%
\text{ }\mathit{\in }\text{ }\mathit{(-\infty ,l(0)-\varepsilon ]} &&  \notag
\end{eqnarray}%
($[l(0)+\varepsilon ,\mathit{\ }m(0)-\varepsilon ]$ or $[m(0)+\varepsilon ,%
\mathit{\ }\infty ),$ resp.). Now suppose $v^{\prime \prime }(\rho _{j})$ $%
\rightarrow $ $-\infty $ as $\rho _{j}\rightarrow 0.$ Clearly $v^{\prime }$
is increasing at $\rho _{j}.$ Let $\bar{c}$ $:=$ $\lim \sup v^{\prime }(\rho
_{j})$ $\in $ $(-\infty ,\mathit{\ }l(0)-\varepsilon ]\mathit{\ }\cup 
\mathit{\ }[l(0)+\varepsilon ,\mathit{\ }m(0)-\varepsilon ]\mathit{\ }\cup 
\mathit{\ }[m(0)+\varepsilon ,\mathit{\ }\infty ]$ by (\ref{A6'}). Then in
view of (\ref{A8}) we can easily show that as $\rho \rightarrow 0$ $%
v^{\prime }(\rho )$ increases and converges to $\bar{c}$ (in particular, $%
\bar{c}$ $\neq $ $l(0)+\varepsilon $ and $\bar{c}$ $\neq $ $m(0)+\varepsilon
).$ So there exists $\hat{\rho}$ $>$ $0$ (may depend on $\bar{c})$ such that 
$v^{\prime }(\hat{\rho})$ $\in $ $(-\infty ,\mathit{\ }l(0)-\varepsilon ]%
\mathit{\ }\cup \mathit{\ }[l(0)+\varepsilon ,\mathit{\ }m(0)-\varepsilon ]%
\mathit{\ }\cup \mathit{\ }[m(0)+\varepsilon ,\mathit{\ }\infty ),$ $%
v^{\prime \prime }(\rho )$ $<$ $0$ for all $0<\rho <\hat{\rho},$ and%
\begin{equation}
v^{\prime }(\rho )v^{\prime \prime }(\rho )\text{ has the same sign for all }%
0<\rho <\hat{\rho}.  \label{A9}
\end{equation}

\noindent By the assumption $E_{1}(0)$ $>$ $0$ and $u^{\prime \prime }$ $<$ $%
0$ for all $0<\rho <\hat{\rho},$ we have 
\begin{equation}
\bar{c}\notin (-\infty ,\mathit{\ }l(0)-\varepsilon ]\cup \lbrack
m(0)+\varepsilon ,\mathit{\ }\infty ]  \label{A9'}
\end{equation}%
in view of equation (\ref{A1}).

On the other hand, we deduce from (\ref{A2'}) that%
\begin{eqnarray}
&&|\frac{v^{\prime }(\rho )v^{\prime \prime }(\rho )}{(v^{\prime }(\rho
)-l(0))(v^{\prime }(\rho )-m(0))}|  \label{A10} \\
&\geq &\frac{|E_{1}(0)|}{16}|\frac{v^{\prime }(\rho )}{v(\rho )}|=\frac{%
|E_{1}(0)|}{16}|(\log v)^{\prime }(\rho )|.  \notag
\end{eqnarray}

\noindent We express the left side of (\ref{A10}) as follows:%
\begin{eqnarray}
&&\frac{v^{\prime }(\rho )v^{\prime \prime }(\rho )}{(v^{\prime }(\rho
)-l(0))(v^{\prime }(\rho )-m(0))}  \label{A11} \\
&=&\frac{\alpha (v^{\prime }(\rho )-l(0))^{\prime }(\rho )}{v^{\prime }(\rho
)-l(0)}+\frac{\beta (v^{\prime }(\rho )-m(0))^{\prime }(\rho )}{v^{\prime
}(\rho )-m(0)}  \notag \\
&=&(\log |v^{\prime }(\rho )-l(0)|^{\alpha }|v^{\prime }(\rho )-m(0)|^{\beta
})^{\prime }(\rho )  \notag
\end{eqnarray}

\noindent where $\alpha $ $:=$ $\frac{-l(0)}{m(0)-l(0)},$ $\beta $ $:=$ $%
\frac{m(0)}{m(0)-l(0)}.$ Now substituting (\ref{A11}) into (\ref{A10}) and
integrating (\ref{A10}) over $\rho $ $\in $ $(0,\hat{\rho}),$ we obtain%
\begin{eqnarray}
&&\int_{0}^{\hat{\rho}}|(\log |v^{\prime }(\rho )-l(0)|^{\alpha }|v^{\prime
}(\rho )-m(0)|^{\beta })^{\prime }(\rho )|d\rho  \label{A12} \\
&\geq &\frac{|E_{1}(0)|}{16}\int_{0}^{\hat{\rho}}|(\log v)^{\prime }(\rho
)|d\rho  \notag \\
&\geq &\frac{|E_{1}(0)|}{16}|\int_{0}^{\hat{\rho}}(\log v)^{\prime }(\rho
)d\rho |  \notag \\
&=&\frac{|E_{1}(0)|}{16}|(\log v(\hat{\rho})-\log v(0))|=+\infty .  \notag
\end{eqnarray}

\noindent On the other hand, either $(\log |v^{\prime }(\rho )-l(0)|^{\alpha
}|v^{\prime }(\rho )-m(0)|^{\beta })^{\prime }(\rho )$ is positive for all $%
0 $ $<$ $\rho $ $<$ $\hat{\rho}$ or negative for all $0$ $<$ $\rho $ $<$ $%
\hat{\rho}$ by (\ref{A9}) and (\ref{A11}). But we then have%
\begin{eqnarray}
&&\int_{0}^{\hat{\rho}}(\log |v^{\prime }(\rho )-l(0)|^{\alpha }|v^{\prime
}(\rho )-m(0)|^{\beta })^{\prime }(\rho )d\rho  \label{A13} \\
&=&\log (|v^{\prime }(\hat{\rho})-l(0)|^{\alpha }|v^{\prime }(\hat{\rho}%
)-m(0)|^{\beta })-\log (|c-l(0)|^{\alpha }|c-m(0)|^{\beta }),  \notag
\end{eqnarray}

\noindent a finite number ($\bar{c}$ $\neq $ $+\infty $ by (\ref{A9'})),
contradicting (\ref{A12}). For the situation that $v^{\prime \prime }(\rho
_{j})$ $\rightarrow $ $+\infty $ as $\rho _{j}\rightarrow 0,$ we have a
similar reasoning with $v^{\prime }$ being decreasing at $\rho _{j}$ and \b{c%
} $:=$ $\lim \inf v^{\prime }(\rho _{j})$ $\in $ $[-\infty ,\mathit{\ }%
l(0)-\varepsilon ]\mathit{\ }\cup \mathit{\ }[l(0)+\varepsilon ,\mathit{\ }%
m(0)-\varepsilon ]\mathit{\ }\cup \mathit{\ }[m(0)+\varepsilon ,\mathit{\ }%
\infty )$ by (\ref{A6'}). Then in view of (\ref{A8}) we can also show that
as $\rho \rightarrow 0$ $v^{\prime }(\rho )$ decreases and converges to \b{c}
(in particular, \b{c} $\neq $ $l(0)-\varepsilon $ and \b{c} $\neq $ $%
m(0)-\varepsilon ).$ Since $v(\rho )$ $>$ $0$ for $0$\textit{\ }$<$\textit{\ 
}$\rho $\textit{\ }$<$\textit{\ }$\rho _{0},$ we can find a sequence of $%
a_{j}\rightarrow 0$ such that $v^{\prime }(a_{j})$ $>$ $0.$ This property
implies that \b{c} $\geq $ $0$ (in particular, \b{c} $=$ $-\infty $ is
excluded). By a similar argument as in (\ref{A10})-(\ref{A13}) we finally
reach a contradiction again. We have proved (\ref{A1'}), hence (b).

Now suppose $\lim_{\rho \rightarrow 0}v(\rho )$ $:=$ $v(0)$ $\neq $ $0$ (so $%
v(0)$ $>$ $0$ since $v$ $>$ $0$ in $(0,\rho _{0})).$ We still want to prove
the existence of $\lim_{\rho \rightarrow 0}v^{\prime }(\rho ).$ If $%
v^{\prime }$ is bounded in $(0,\hat{\rho}_{0})$ for $0$ $<$ $\hat{\rho}_{0}$ 
$<$ $\rho _{0},$ then $v^{\prime \prime }$ is bounded in $(0,\bar{\rho}_{0})$
for $\bar{\rho}_{0}$ small, $0$ $<$ $\bar{\rho}_{0}$ $<$ $\hat{\rho}_{0}$ by
(\ref{A1}). It follows that $v^{\prime }$ is Cauchy in $(0,\bar{\rho}_{0})$
since%
\begin{eqnarray*}
|v^{\prime }(\check{\rho})-v^{\prime }(\tilde{\rho})| &=&|\int_{\tilde{\rho}%
}^{\check{\rho}}v^{\prime \prime }(\rho )d\rho | \\
&\leq &\int_{\tilde{\rho}}^{\check{\rho}}|v^{\prime \prime }(\rho )|d\rho
\leq C_{1}|\check{\rho}-\tilde{\rho}|
\end{eqnarray*}

\noindent for $0$ $<$ $\tilde{\rho}$ $<$ $\check{\rho}$ $<$ $\bar{\rho}_{0},$
where $|v^{\prime \prime }(\rho )|$ $\leq $ $C_{1},$ a positive constant
independent of $\rho ,$ for $\rho $ $\in $ $(0,\bar{\rho}_{0}).$ So $%
\lim_{\rho \rightarrow 0}v^{\prime }(\rho )$ exists. On the other hand, if $%
v^{\prime }$ is not bounded near $0,$ then there exists a sequence $\rho
_{j}\rightarrow 0$ such that $\lim_{\rho _{j}\rightarrow 0}v^{\prime }(\rho
_{j})$ $=$ $-\infty $ ($+\infty $ is impossible by a similar argument as in
the first paragraph of the proof since $v^{\prime \prime }(\rho _{j})$ $>$ $%
0 $). In fact we can easily show that $v^{\prime }(\rho )$ is monotonically
decreasing to $-\infty $ as $\rho $ $\rightarrow $ $0$ since $v^{\prime
\prime }$ $>$ $0$ for \TEXTsymbol{\vert}$v^{\prime }|$ large by (\ref{A1}).
We can find $\rho _{0}^{\prime }$ $>$ $0$ small so that $v^{\prime }$ is
negatively large in $(0,\rho _{0}^{\prime }]$ and there holds%
\begin{equation}
v(\rho )v^{\prime \prime }(\rho )\leq C_{2}(v^{\prime }(\rho ))^{2}
\label{A14}
\end{equation}

\noindent for $\rho $ $\in $ $(0,\rho _{0}^{\prime }]$ and some positive
constant $C_{2}$ independent of $\rho .$ Dividing (\ref{A14}) by $%
-vv^{\prime }$ $>$ $0$ (noting that $v^{\prime }$ $<$ $0$ in $(0,\rho
_{0}^{\prime }]$) we obtain%
\begin{eqnarray}
-(\log (-v^{\prime }))^{\prime } &=&-\frac{v^{\prime \prime }}{v^{\prime }}
\label{A15} \\
&\leq &C_{2}(-\frac{v^{\prime }}{v})=-C_{2}(\log v)^{\prime }.  \notag
\end{eqnarray}

\noindent Integrating (\ref{A15}) from $\varepsilon ,$ $0$ $<$ $\varepsilon $
$<$ $\rho _{0}^{\prime }$ to $\rho _{0}^{\prime }$ we get%
\begin{equation}
-\log (-v^{\prime }(\rho _{0}^{\prime }))+\log (-v^{\prime }(\varepsilon
))\leq C_{2}\{-\log v(\rho _{0}^{\prime })+\log v(\varepsilon )\}.
\label{A16}
\end{equation}

\noindent Letting $\varepsilon $ $\rightarrow $ $0$ in (\ref{A16}) we reach +%
$\infty $ $\leq $ $C_{2}\{-\log v(\rho _{0}^{\prime })+\log v(0)\}$ (note
that $v(0)$ $>$ $0),.$a contradiction. We have proved the existence of $%
\lim_{\rho \rightarrow 0}v^{\prime }(\rho )$ and completed the proof of (a).

\endproof%

\bigskip

We remark that Theorem A.1 can be applied to show a generalized version of
Theorem B (a) (see Theorem B$^{\prime }$ in Section 8). Since $v(0)$ is not
necessarily $0,$ we may also apply Theorem A.1 to the piecewise $C^{1}$ case
in which $\{\rho =0\}$ corresponds to a nonsmooth edge.

\end{document}